\documentclass[fleqn,11pt,oneside]{article}

\usepackage{hyperref}
\usepackage{graphicx}
\usepackage{bbm}  
\usepackage{bm}  
\usepackage{amsmath,amsthm,amssymb,amscd}
\usepackage{caption,subcaption}
\usepackage{booktabs,multirow}
\usepackage{setspace} 
\usepackage[top=1.1in, bottom=1.1in, left=1.6in, right=1.1in]{geometry}
\usepackage[inline,final]{showlabels}
\usepackage{dashbox}
\usepackage{microtype}
\usepackage[medium]{titlesec}
\usepackage{algpseudocode}
\usepackage[ruled]{algorithm}
\usepackage{xpatch}
\usepackage{accents}
\usepackage{bbding}
\usepackage[table]{xcolor}
\usepackage{srcltx}

\newtheorem{definition1}{Definition}[section]
\newtheorem{observe}{Observation}[section]
\newtheorem{remark1}[observe]{Remark}
\newtheorem{example1}{Example}[section]
\newtheorem{aside1}[observe]{Aside}

\newenvironment{remark}{\begin{remark1} \rm}{\end{remark1}}

\def\qed{\hfill$\blacksquare$\\} 

\algrenewcommand\algorithmicindent{3.0em}

\makeatletter
\xpatchcmd{\algorithmic}{\itemsep\z@}{\itemsep=2.3ex}{}{}
\makeatother

\newif\ifshowboxes \showboxestrue

\renewcommand{\hat}{\widehat}
\renewcommand{\tilde}{\widetilde}

\newcommand{\norm}[1]{\ensuremath{ {\lVert #1 \rVert} }}

\newcommand{\Bnorm}[1]{\ensuremath{ {\Bigl\lVert #1 \Bigr\rVert} }}

\def\1{\mathbbm{1}}

\newcommand{\e}[1]{\ensuremath{\cdot 10^{#1}}}
\newcommand\T{\rule{0pt}{2.6ex}}       

\begin{document}

\begin{center}
   \begin{minipage}[t]{6.0in}

In this paper, we describe an algorithm for fitting an 
analytic and bandlimited closed or open curve to interpolate an arbitrary collection of points in $\mathbb{R}^{2}$.
The main idea is to smooth the parametrization of the curve by iteratively filtering 
the Fourier or Chebyshev coefficients of both the derivative of the 
arc length function and the tangential angle of the curve, and applying smooth perturbations, after
each filtering step,
until the curve is represented by a reasonably
small number of coefficients. The algorithm produces a curve passing 
through the set of points to an accuracy of machine 
precision,
after a limited number of iterations. It costs O($N\log{N}$) operations at each iteration, provided that the number
of discretization nodes is $N$. 
The resulting curves are smooth and visually appealing, and do not exhibit any
ringing artifacts. The bandwidths of the constructed curves are much
smaller than those of curves constructed by previous methods.
We demonstrate the performance of our algorithm with several numerical
experiments.

\thispagestyle{empty}

  \vspace{ -100.0in}

  \end{minipage}
\end{center}

\vspace{ 2.60in}
\vspace{ 0.50in}

\begin{center}
  \begin{minipage}[t]{4.4in}
    \begin{center}

\textbf{A Continuation Method for Fitting a Bandlimited Curve to Points in the Plane}\\

  \vspace{ 0.50in}

Mohan Zhao$\mbox{}^{\dagger \, \diamond}$ 
and
Kirill Serkh$\mbox{}^{\ddagger \, \diamond}$  \\
University of Toronto NA Technical Report \\
              \today

    \end{center}
  \vspace{ -100.0in}
  \end{minipage}
\end{center}

\vspace{ 2.00in}

\vfill

\noindent 
$\mbox{}^{\diamond}$  This author's work  was supported in part by the NSERC
Discovery Grants RGPIN-2020-06022 and DGECR-2020-00356.
\\

\vspace{2mm}

\noindent
$\mbox{}^{\dagger}$ Dept.~of Computer Science, University of Toronto,
Toronto, ON M5S 2E4 \\
Corresponding author. Email: mohan.zhao@mail.utoronto.ca \\

\noindent
$\mbox{}^{\ddagger}$ Dept.~of Math. and Computer Science, University of Toronto,
Toronto, ON M5S 2E4 \\
Email: kserkh@math.toronto.edu \\

\vspace{2mm}



\vfill
\eject

\tableofcontents

\section{Introduction}
The construction of smooth curves passing through data points has 
uses in many areas of applied
science, including
boundary integral equation methods, computer graphics and geometric modeling. 
While much of the time, $C^{k}$ continuity is sufficient, there are certain applications
for which $C^{\infty}$ continuity is essential.
One such example is the high accuracy solution of partial differential equations
on general geometries. In CAD/CAM 
systems, $C^{\infty}$ smooth curves 
can be used as primitives to construct arbitrary
smooth objects. Solving partial differential equations on these smooth objects   
prevents the loss of  
accuracy due to imperfect smoothness of $C^k$ shapes.

Countless methods have been proposed for fitting a spline or a $C^{k}$ curve to a given
set of data points.
Most interpolation techniques use piecewise polynomials and impose 
constraints to ensure global $C^k$ smoothness of the curve 
(see, for example,~\cite{new}, \cite{stiff}, \cite{bica}, \cite{intro}).
In CAD/CAM systems, 
non-uniform rational B-splines (NURBS) are commonly used to
construct a curve which
approximates a set of control points,
by defining the curve as a linear combination of the control points multiplied
by $C^{k}$ and compactly supported B-spline basis functions.
The contribution of each control point
to the overall curve is 
determined by the corresponding weight, and the B-spline basis functions 
are normalized to 
ensure that the approximating curve remains affine invariant \cite{nurbs}.
A generalization of NURBS, called  
partition of unity parametrics (PUPs) was
first introduced by Runions and Samavati (\cite{pup}). The PUP curves
are constructed by
replacing the weighted B-spline
basis functions with arbitrary normalized weight functions (WFs), 
so that the resulting curves
still exhibit the desired properties, including compact support and $C^{k}$ smoothness.
In \cite{pup}, the authors specifically discuss uniform B-spline WFs,
to illustrate  
that
each WF can be adjusted independently to fine-tune 
various shape parameters of the curve.
Additionally, they observe that
it is possible to choose the WFs to generate a PUP curve 
that interpolates the control points without solving a system of equations.

Another method proposed by Zhang and Ma (\cite{zhang}) 
employs products of the sinc function and Gaussian functions as
basis functions for constructing
$C^{\infty}$ interpolating curves that pass through all the given data points
exactly.
The resulting curves
are almost affine invariant and almost compactly 
supported, and their shapes can be adjusted locally by directly adding or moving
control points.
Subsequently, Runions and Samavati (\cite{cinpact}) designed 
CINPACT-splines, 
by employing $C^{\infty}$ and compactly supported bump functions
as the WFs in a PUP curve, optionally multiplied by the normalized sinc function.
When the WFs are chosen to be products of bump functions and the normalized
sinc function, the resulting $C^{\infty}$ curve interpolates the control
points exactly, 
and when the WFs are bump functions,
the resulting $C^{\infty}$ curve approximates
a uniform B-spline with 
the given control
points. 
In addition to the properties
inherited from PUP curves,
CINPACT-splines possess $C^{\infty}$ smoothness
and the ability to specify tangents
at control points.
To increase the accuracy of the approximation to uniform B-splines,
Akram, Alim and Samavati (\cite{cinapact}) further proposed CINAPACT-splines,
by successively convolving a CINPACT-spline with B-splines of
order one, ensuring
any finite order of 
approximation to uniform B-splines,
as well as to other compactly supported kernels with maximal order 
and minimal support (\cite{moms}),
while preserving
$C^{\infty}$ smoothness and compact support.
Zhu (\cite{zhu}) proposed curves that share similarities with CINPACT-splines
in terms of
affine invariance, compact support, and $C^{\infty}$ 
smoothness.
In \cite{zhu},  
a class of non-negative blending functions is constructed 
by designing basis functions which
combine 
bump functions 
with the sinc function.
The resulting interpolating curves are defined by
three local shape parameters, 
with
one of the parameters determining whether the curve 
approximates or interpolates the given control points.

One notable distinction of the approach of Zhang and Ma (\cite{zhang}) from
the other methods we have discussed 
is that,
since Gaussian functions are utilized in the
basis functions, the 
interpolating curves produced by \cite{zhang} are not only $C^{\infty}$
smooth,
but also are analytic.
This paper mainly compares our method with \cite{zhang}, as
the interpolating curves 
in \cite{zhang}
have a smaller bandwidth, compared to methods based on 
$C^{\infty}$ compactly supported bump functions.
The approach in \cite{zhang} (as well as \cite{pup}, \cite{cinpact}, \cite{cinapact},
\cite{zhu}) necessitates a more 
specially chosen distribution of data points to achieve a visually
smooth curve, as it only guarantees smoothness in the curve parameter,
which does not necessarily correspond to smoothness of the curve
in $\mathbb{R}^2$.
However, our method directly smooths 
the tangential angle of the curve 
and the first derivative of the arc length function,
yielding a significantly
smoother curve which is also more visually appealing.

Among all the methods for constructing
a $C^{\infty}$ interpolating curve, the algorithm
described by Beylkin and Rokhlin (\cite{bandlimited})
bears the closest resemblance to our method, generating 
a bandlimited closed curve through a set of data points.
The bandlimited curve is constructed by filtering the Fourier coefficients 
of the tangential angle of the curve,
parametrized by arc length. However, the number of coefficients required to 
represent the curve can be large, which appears to be a major drawback of the algorithm 
in practical applications.

In this paper, we describe an algorithm for fitting a bandlimited closed or 
open curve to pass through a collection of points. The main idea is to iteratively filter 
the tangential angle and the first derivative of the arc length function of the curve,
and apply small corrections after each filtering step, until 
the desired bandwidth of the curve is reached, to the required precision.
Our algorithm produces an analytic and affine invariant curve 
with
far fewer coefficients, and the curve is visually appealing and free of
ringing artifacts.

The structure of this paper is as follows. Section~\ref{sec:prelim}
describes the mathematical preliminaries.
Section~\ref{sec:algorithm} describes the algorithm to construct the 
bandlimited approximation to a closed curve, and to an open curve. 
Finally, Section~\ref{sec:numerical} presents several numerical 
examples to show the performance of our algorithm, as well as 
some comparisons between our algorithm and the methods proposed in \cite{zhang}
and \cite{bandlimited}.

\section{Preliminaries}
  \label{sec:prelim}

In this section, we describe the mathematical and numerical preliminaries.

\subsection{Geometric properties of a curve}
  \label{sec:geometry}
Let $\gamma\colon[a,b]\rightarrow \mathbb{R}^{2}$ 
be a smooth curve parametrized by the curve parameter $t$,
such that
  \begin{align}
\gamma(t)=(x(t),y(t)), \qquad t \in [a,b],
  \end{align}
where $x(t)$ and $y(t)$ are the $x$ and $y$ coordinates.

Assuming $\gamma \in C^1([a,b])$, we define the tangent vector $T(t)$,
  \begin{align}
  T(t)=(x'(t),y'(t)),  \qquad t \in [a,b],
  \end{align}
and the arc length $s(t)$, which is the length of the curve from the point $(x(a)$, $y(a))$
to the point $(x(t)$, $y(t))$, 
  \begin{align}
  s(t)=\int_{a}^{t}\norm{T(\tau)} d\tau, \qquad t \in [a,b]. 
  \end{align}
It is obvious that 
  \begin{align}
  \label{eq:sprime} 
  s'(t)=\norm{T(t)}, \qquad t \in [a,b]. 
  \end{align}
Thus, we have 
  \begin{align}
  \label{eq:sprimecon} 
  s'(b)=s'(a)
  \end{align}
when the curve is closed.
The tangential angle $\theta(t)$ of the curve at the point $(x(t)$, $y(t))$ measures the angle between 
the tangent vector $T(t)$ at that point and the x-axis, defined by the formula 
  \begin{align}
  \label{eq:theta} 
  \theta(t)=\text{atan2}(y'(t),x'(t)), \qquad t \in [a,b], 
  \end{align}
where $\text{atan2}\colon\mathbb{R}^{2} \rightarrow (-\pi,\pi]$  is the arctangent at the point $(x(t),y(t))$. As
a result, $\theta(t) \in (-\pi,\pi]$.
Since the function $\text{atan2}$ has a branch cut at $\theta=-\pi$, it is possible for $\theta(t)$ to have $\omega$ jump discontinuities of size $2\pi$,
where $\omega \in \mathbb{Z}$ is the winding number.

The curve $(x(t)$, $y(t))$ can be constructed from $\theta(t)$ and $s'(t)$ by the formulas
  \begin{align}
  x(t) = \int_{a}^{t} s'(\tau)\cos{\theta(\tau)}\, d\tau + x(a), \qquad t \in [a,b], \label{eq:x}\\ 
  y(t) = \int_{a}^{t} s'(\tau)\sin{\theta(\tau)}\, d\tau + y(a), \qquad t \in [a,b], \label{eq:y}
  \end{align}
and $(x(a),y(a)) = \gamma(a)$.
If the curve is closed, we require $x(a) = x(b)$ and 
$y(a) = y(b)$, which means that
  \begin{align}
  \label{eq:cos} 
\int_{a}^{b} s'(\tau)\cos{\theta(\tau)}\, d\tau = 0
  \end{align}
and 
  \begin{align}
  \label{eq:sin} 
\int_{a}^{b} s'(\tau)\sin{\theta(\tau)}\, d\tau = 0.
  \end{align}

\subsection{Cubic B\'ezier Interpolation}

A B\'ezier curve is a function $\textbf{B}\colon[0,1]\rightarrow \mathbb{R}^{2}$ defined by a set of control points
$\textbf{P}_0$, \ldots, $\textbf{P}_{m} \in \mathbb{R}^{2}$.
The B\'ezier curve is designed to go through the first and
and the last control point ${\textbf{P}}_0$ and
${\textbf{P}}_m$, and the shape of the curve is determined by the intermediate 
control points $\textbf{P}_1$, \ldots, $\textbf{P}_{m-1}$. A $m$th order B\'ezier curve
is a polynomial of degree $m$, defined by
  \begin{align*}
\hspace*{-5em}\textbf{B}(t)& = \sum_{i=0}^{m} \binom{m}{i}(1-t)^{m-i}t^{i}\textbf{P}_i, \notag\\
\hspace*{-5em}&= (1-t)^{m}\textbf{P}_0 + \binom{m}{1}(1-t)^{m-1}t\textbf{P}_1 + 
                \dots + \binom{m}{m-1}(1-t)t^{m-1}\textbf{P}_{m-1} + t^{m}\textbf{P}_{m},
  \end{align*}
where $t \in [0,1]$.

A continuous B\'ezier spline connecting all the given points $\textbf{C}_0$, \ldots, $\textbf{C}_n$ 
can be constructed by combining $n$ cubic B\'ezier curves
$$
  \begin{aligned}
\textbf{B}_{i}(t) = (1-t)^{3}\textbf{P}_{i0} + 3(1-t)^{2}t\textbf{P}_{i1} +
     3(1-t)t^{2}\textbf{P}_{i2} + t^{3}\textbf{P}_{i3}, \qquad i = 1, \ldots, n,
  \end{aligned}
$$
where $t \in [0,1]$ and $\textbf{B}_{i}(t)$ is the $i$th B\'ezier curve, with controls points
\begin{align}
&\textbf{P}_{i0} = \textbf{C}_{i-1} \label{eq:con1},\\
&\textbf{P}_{i3} = \textbf{C}_{i} \label{eq:con2},
\end{align}
for $i = 1$, \ldots, $n$. 
We define the spline $\textbf{S}\colon[0,n]\rightarrow \mathbb{R}^{2}$ from the cubic B\'ezier curves 
$\textbf{B}_{i}$ by letting $\textbf{S}(t)$ = $\textbf{B}_{i}(t-i+1)$ for $t \in [i-1,i]$, for $i = 1,$ \ldots, $n$.
It is easy to see that $\textbf{S} \in {C}^0([0,n])$, however, in general,
$\textbf{S} \notin {C}^1([0,n])$. It is possible to ensure
$\textbf{S} \in {C}^2([0,n])$ by imposing additional conditions on the intermediate control points, which we derive
as follows.
Note that the following derivation is similar to the one presented in \cite{controlpoints}.
First, we observe that the first and second derivatives of a cubic B\'ezier curve are
  \begin{align}
\notag \textbf{B}_{i}'(t) &= -3(1-t)^{2}\textbf{P}_{i0} + 3(3t^{2}-4t+1)\textbf{P}_{i1} +
     3t(2-3t)\textbf{P}_{i2} + 3t^{2}\textbf{P}_{i3}, \\
\notag \textbf{B}_{i}''(t) &= 6(1-t)\textbf{P}_{i0} + 6(3t-2)\textbf{P}_{i1} +
     6(1-3t)\textbf{P}_{i2} + 6t\textbf{P}_{i3},
  \end{align}
for $i=1$, \ldots, $n$.
In order for $\textbf{S} \in {C}^2([0,n])$, we require that 
  \begin{align}
\textbf{B}_{i-1}'(1) = \textbf{B}_{i}'(0), \label{eq:con1}\qquad i = 1, \ldots, n,\\
\textbf{B}_{i-1}''(1) = \textbf{B}_{i}''(0), \label{eq:con2}\qquad i = 1, \ldots, n.
  \end{align}
Then, \eqref{eq:con1} implies that 
  \begin{align}
  \label{eq:res1} 
\textbf{P}_{(i-1)2} = 2\textbf{C}_{i-1} - \textbf{P}_{i1},\qquad i = 1, \ldots, n. 
  \end{align}
Likewise, it is possible to show that \eqref{eq:con2} implies that 
  \begin{align}
  \label{eq:res2} 
\textbf{P}_{(i-1)1} + 2\textbf{P}_{i1} = \textbf{P}_{i2} + 2\textbf{P}_{(i-1)2},
\qquad i = 1, \ldots, n. 
  \end{align}
Substituting \eqref{eq:res1} into \eqref{eq:res2},we get
  \begin{align}
  \label{eq:res3} 
\textbf{P}_{(i-1)1} + 4\textbf{P}_{i1} + \textbf{P}_{(i+1)1}= 
  2\textbf{C}_{i} + 4\textbf{C}_{i-1},\qquad i = 1, \ldots, n. 
  \end{align}
\subsubsection{Solving for control points for an open curve}
  \label{sec:bezieropen}

When the curve is open, we have \eqref{eq:res3} must hold for $i = 2$, \ldots, $n-1$, and we need  
two boundary conditions in order to solve a linear system of $n$ equations for the values of 
$\textbf{P}_{11}$, \ldots, $\textbf{P}_{n1}$. Assume that users specify the slope at  
two end points of the curve, $c_{\text{left}}$ and $c_{\text{right}}$, we have 
  \begin{align}
  \label{eq:bcon1} 
\textbf{B}_{1}'(0) = c_{\text{left}},
  \end{align}
and
  \begin{align}
  \label{eq:bcon2} 
\textbf{B}_{n}'(1) = c_{\text{right}}.
  \end{align}
It is possible to show that \eqref{eq:bcon1} implies that 
  \begin{align}
  \label{eq:bres1}
\textbf{P}_{11} = \frac{c_{\text{left}}+3\textbf{C}_{0}}{3} 
  \end{align}
and \eqref{eq:bcon2} implies that 
  \begin{align}
  \label{eq:bres2}
\textbf{P}_{n2} = \frac{3\textbf{C}_{n}-c_{\text{right}}}{3}. 
  \end{align}
Substituting \eqref{eq:res1} and \eqref{eq:bres2} into \eqref{eq:res2},we get
  \begin{align}
  \label{eq:bres3}
\textbf{P}_{(n-1)1} + 4\textbf{P}_{n1} = 4\textbf{C}_{n-1} + \textbf{C}_{n} - \frac{c_{\text{right}}}{3}. 
  \end{align}

With \eqref{eq:res3}, \eqref{eq:bres1} and \eqref{eq:bres3}, we build a system of $n$
equations
to calculate $\textbf{P}_{11},$ \ldots, $\textbf{P}_{n1}$ and use \eqref{eq:res1},  
$\eqref{eq:bres2}$ and the values of $\textbf{P}_{11}$, \ldots, $\textbf{P}_{n1}$  
to calculate $\textbf{P}_{12}$, \ldots, $\textbf{P}_{n2}$. This
system of equations is tridiagonal, and so can be solved in $O(n)$ operations.

\subsubsection{Solving for control points for a closed curve}
  \label{sec:bezierclosed}

When the curve is closed, we require $n+1$ cubic B\'ezier curves instead of $n$ cubic B\'ezier curves 
to connect the points $\textbf{C}_0$, \ldots, $\textbf{C}_n$, where the  
$(n+1)$th curve connects the points $\textbf{C}_{n}$ and $\textbf{C}_{0}$.
We have that the conditions \eqref{eq:res3} 
must hold for $i = 2$, \ldots, $n$, and we need the following two boundary conditions,
  \begin{align}
\textbf{B}_{1}'(0) &= \textbf{B}_{n+1}'(1), \label{eq:bcon3} \\
\textbf{B}_{1}''(0) &= \textbf{B}_{n+1}''(1)\label{eq:bcon4}, 
  \end{align}
to solve a linear system of
$(n+1)$ equations for the values of $\textbf{P}_{11}$, \ldots, $\textbf{P}_{(n+1)1}$.

It is possible to show that \eqref{eq:bcon3} implies that 
  \begin{align}
  \label{eq:bres4}
\textbf{P}_{11} + \textbf{P}_{(n+1)2} =  2\textbf{C}_{0} 
  \end{align}
and \eqref{eq:bcon4} implies that 
  \begin{align}
  \label{eq:bres5}
- 2\textbf{P}_{11} + \textbf{P}_{12} = \textbf{P}_{(n+1)1} - 2\textbf{P}_{(n+1)2}.
  \end{align}
Substituting \eqref{eq:res1} and \eqref{eq:res2} into \eqref{eq:bres4}, we get
  \begin{align}
  \label{eq:bres6}
\textbf{P}_{11} + \textbf{P}_{n1} + 4\textbf{P}_{(n+1)1} = 2\textbf{C}_{0} + 4\textbf{C}_{n}. 
  \end{align}
Substituting \eqref{eq:res1} and \eqref{eq:res2} into \eqref{eq:bres5}, we get
  \begin{align}
  \label{eq:bres7}
- 2\textbf{P}_{11} - \textbf{P}_{21} + 2\textbf{P}_{n1} + 7\textbf{P}_{(n+1)1} = -2\textbf{C}_{1} + 8 \textbf{C}_{n}. 
  \end{align}

Similarly, with \eqref{eq:res3}, \eqref{eq:bres6} and \eqref{eq:bres7}, we build a 
system of $(n+1)$ equations
to calculate $\textbf{P}_{11}$, \ldots, $\textbf{P}_{(n+1)1}$ and use \eqref{eq:res1},  
$\eqref{eq:bres4}$ and the values of $\textbf{P}_{11}$, \ldots, $\textbf{P}_{(n+1)1}$ 
to calculate $\textbf{P}_{12}$, \ldots, $\textbf{P}_{(n+1)2}$. This system of equations is
cyclic tridiagonal, and thus we can solve it in $O(n)$ operations.

\subsection{Chebyshev Polynomial Interpolation}
A smooth function $f(x)$ on the interval $[-1,1]$ can be approximated by a $(n-1)$th order Chebyshev expansion 
with the formula 
  \begin{align}
  \label{eq:cheb1}
f(x)\approx \sum_{k=0}^{n-1}\hat{f}_kT_k(x), 
  \end{align}
where $T_k(x)$ is the Chebyshev polynomial of the first kind of degree $k$, 
defined by
  \begin{align}
  T_{k}(x)=\cos(k\arccos{x}), \qquad x \in [-1,1].
  \end{align}
It is known that the Chebyshev coefficients $\{\hat{f}_k\}$ decay like $O(n^{-k+\frac{1}{2}})$
when $f \in {C}^{k}([-1,1])$,
when the coefficients $\hat{f}_k$ are chosen to satisfy the $n$ collocation equations 
  \begin{align}
  \label{eq:chebcoefsdef}
f(x_i)= \sum_{k=0}^{n-1}\hat{f}_kT_k(x_i),\qquad i = 0, \ldots, n-1,
  \end{align}
for the practical Chebyshev nodes $\{x_i\}$,
  \begin{align}
  \label{eq:chebnodes}
x_i=-\cos\bigl{(}\frac{i\pi}{n-1}\bigr{)}, \qquad i = 0, \ldots, n-1.
  \end{align}

Alternatively, one can compute $\hat{f}_k$ for $k = 0,\ldots,n-1$ using the Discrete Chebyshev Transform,
  \begin{align}
  \label{eq:chebcoefs0}
\hat{f}_0=\frac{1}{n-1}\bigl{(}\frac{1}{2}(f(x_0)+f(x_{n-1})\bigr{)}+
\sum_{i=1}^{n-2}f(x_i)T_0(x_i),
  \end{align}
and
  \begin{align}
  \label{eq:chebcoefs}
\hat{f}_k=\frac{2}{n-1}\bigl{(}\frac{1}{2}(f(x_0)(-1)^k+f(x_{n-1})\bigr{)}+
\sum_{i=1}^{n-2}f(x_i)T_k(x_i), 
  \end{align}
for $k=1$, \ldots, $n-1$.
Once the coefficients $\{\hat{f}_k\}$ are computed, we can use the expansion $\sum_{k=0}^{n-1}\hat{f}_kT_k(x)$ to evaluate 
$f(x)$ everywhere on the interval $[-1,1]$.

\subsubsection{Spectral Differentiation and Integration}
Assuming that $k \geq 1$ is an integer, the formula
  \begin{align}
2T_k(x)=\frac{T'_{k+1}(x)}{k+1}-\frac{T'_{k-1}(x)}{k-1},
  \end{align}
can be used to spectrally differentiate the Chebyshev expansion of $f(x)$, as follows. 
Suppose that
  \begin{align}
f(x)\approx\sum_{k=0}^{n-1}\hat{f}_kT_k(x) 
  \end{align}
and that
  \begin{align}
  \label{eq:chebdiff}
f'(x)\approx \sum_{k=0}^{n-1}\hat{f}'_kT_k(x). 
  \end{align}
The coefficients $\hat{f}'_k$ can be computed from $\hat{f}_k$ by iterating from $k=n-1,n-2$, \ldots, $2$ and, at each iteration,
assigning $\hat{f}'_{k-1}$ the value $2k\hat{f}_k$, and assigning 
$\hat{f}_{k-2}$ the value $\frac{k}{k-2}\hat{f}_k+\hat{f}_{k-2}$. 

Similarly, the formula
  \begin{align}
2\int_{-1}^{t}T_k(x) \,dx = \frac{T_{k+1}(t)}{k+1}-\frac{T_{k-1}(t)}{k-1}
-\frac{(-1)^{k+1}}{k+1}+\frac{(-1)^{k-1}}{k-1}
  \end{align}
can be used to spectrally integrate the Chebyshev expansion of $f(x)$. 
Suppose that
  \begin{align}
  \label{eq:chebint}
\int_{-1}^{t} f(x) \,dx \approx \sum_{k=0}^{n}\hat{\tilde{f}}_kT_k(t).
  \end{align}
Since
  \begin{align}
  \label{eq:chebint}
\int_{-1}^{t} f(x)dx&\approx \sum_{k=0}^{n-1}\hat{f}_k\int_{-1}^{t}T_k(x) \,dx \notag\\
                   &= \sum_{k=1}^{n-1}\hat{f}_k \frac{1}{2}
\Bigl(\frac{T_{k+1}(t)}{k+1}-\frac{T_{k-1}(t)}{k-1}
-\frac{(-1)^{k+1}}{k+1}+\frac{(-1)^{k-1}}{k-1}\Bigr) \notag \\
&  \ \ \ +\hat{f}_0 (t+1),
  \end{align}
one can compute the coefficients $\hat{\tilde{f}}_k$ from $\hat{f}_k$ by firstly
assigning $\hat{{f}}_1$ the value $\hat{f}_1+\hat{f}_0$, then iterating
from $k=n-1$, \ldots, $1$, and at each iteration,
assigning ${\hat{f}}_{k+1}$ the value $\hat{f}_{k+1} +
\frac{\hat{f}_k}{2(k+1)}$,
assigning $\hat{f}_{k-1}$ the value $\hat{f}_{k-1}-\frac{\hat{f}_k}{2(k-1)}$,
and assigning $\hat{f}_0$ the value
$\hat{f}_0-\hat{f}_k(\frac{(-1)^{k+1}}{2(k+1)}-\frac{(-1)^{k-1}}{2(k-1)})$.
Finally, $\hat{\tilde{f}}_k$ takes the value $\hat{f}_k$, for $k=n$, \ldots, $0$.

\subsection{The Discrete Fourier Transform (DFT)}
A periodic and smooth function $f(x)$ on the interval $[0,1]$ can be approximated by 
a $n$-term Fourier series using the Discrete Fourier Transform.
The Discrete Fourier Transform 
defines a transform from a sequence of $n$ complex numbers $f_0$, \ldots, $f_{n-1}$ 
to another sequence of $n$ complex numbers $\hat{f}_0$, \ldots, $\hat{f}_{n-1}$, by
  \begin{align}
  \label{eq:dft1}
\hat{f}_{k} = \sum_{j=0}^{n-1}f_je^{-\frac{2\pi i}{n}kj},\qquad k = 0, \ldots, n-1,
  \end{align}
The sequence $\{\hat{f}_k\}$ consists of the Fourier coefficients of $\{f_k\}$. 

The Inverse Discrete Fourier Transform (IDFT) is given by 
  \begin{align}
f_j = \frac{1}{n}\sum_{k=0}^{n-1}\hat{f}_{k}e^{\frac{2\pi i}{n}kj}, \qquad j = 0, \ldots, n-1.
  \end{align}
Another representation of the DFT which is usually used in applications 
is given by a shift in the index $k$, and a change in the placement of 
the scaling by $\frac{1}{n}$,
  \begin{align}
  \label{eq:rmk1}
\hat{f}_{k} = \frac{1}{n}\sum_{j=0}^{n-1}f_je^{-\frac{2\pi i}{n}kj},\qquad k = -\frac{n}{2}, \ldots, \frac{n}{2}-1.
  \end{align}
Thus, the corresponding IDFT is
  \begin{align}
  \label{eq:rmk2}
f_j = \sum_{k=-\frac{n}{2}}^{\frac{n}{2}-1}\hat{f}_{k}e^{\frac{2\pi i}{n}kj}, \qquad j = 0, \ldots, n-1.
  \end{align}
Suppose that $f\colon [0,1] \rightarrow \mathbb{C}$ is a smooth and periodic function,
and that $f_j = f(t_j)$ for $j=0$, \ldots, $n-1$, where $\{t_j\}$ are the equispaced points on $[0,1]$. 
Observing that
  \begin{align}
\hat{f}_{k} &=\frac{1}{n} \sum_{j=0}^{n-1}f_je^{-\frac{2\pi i}{n}kj}, \notag \\
            & \approx \int_{0}^{1}f(x)e^{-2 \pi ikx}\, dx ,
  \end{align}
for $k = -\frac{n}{2}$, \ldots, $\frac{n}{2}-1$,
we obtain the approximation to $f(x)$ by a truncated Fourier series, 
  \begin{align}
  \label{eq:dft3}
f(x) \approx \sum_{k=-\frac{n}{2}}^{\frac{n}{2}-1}\hat{f}_{k}e^{2\pi i kx}, \qquad x \in [0,1].
  \end{align}
It is known that the Fourier coefficients $\{\hat{f}_k\}$ decay
like $O(n^{-k+\frac{1}{2}})$
when $f \in {C}^{k}(S^1)$, where $S^1 = [0,1]$ is the circle.

\subsubsection{Spectral Differentiation and Integration}

The spectral differentiation of the truncated Fourier series approximation to $f(x)$ on
$[0,1]$ is as follows.
Suppose that $f(x)$ is given by \eqref{eq:dft3} and that
  \begin{align}
f'(x) \approx \frac{1}{n}\sum_{k=-\frac{n}{2}}^{\frac{n}{2}-1}\hat{f}'_{k}e^{2\pi ikx}.
  \end{align}
Since
  \begin{align}
  \label{eq:dftdiff}
f'(x) \approx \frac{1}{n}\sum_{k=-\frac{n}{2}}^{\frac{n}{2}-1}\hat{f}_{k}e^{{2\pi}ikx}\cdot {2\pi}ik,
  \end{align}
the coefficients $\hat{f}'_{k}$ can be computed from $\hat{f}_k$ by 
assigning $\hat{f}'_{k}$ the value $\hat{f}_{k}\cdot {2\pi}ik$ for $k=
-\frac{n}{2}$, \ldots, $\frac{n}{2}-1$.

Similarly, the spectral integration of the truncated Fourier series approximation to $f(x)$ is as follows.
Suppose that
  \begin{align}
\int_{0}^{t} f(x)dx \approx \frac{1}{n} \sum_{k=-\frac{n}{2}}^{\frac{n}{2}-1}\hat{\tilde{f}}_{k}e^{{2\pi}ikt}.
  \end{align}
Since
  \begin{align}
  \label{eq:dftint}
\int_{0}^{t} f(x) dx \approx \frac{1}{n}\sum_{k\neq 0}\frac{\hat{f}_k}{2\pi ik}e^{2\pi ikt}
  -\frac{1}{n}\sum_{k\neq 0}\frac{\hat{f}_k}{2\pi ik}+\frac{1}{n}\hat{f}_0t,
  \end{align}
it is easy to see that, for $\int_{0}^{t} f(x) dx$ to be periodic, it must be the 
case that $\hat{f}_0 =0$. Then, we have
  \begin{align}
  \label{eq:dftint}
\int_{0}^{t} f(x) dx \approx \frac{1}{n}\sum_{k\neq 0}\frac{\hat{f}_k}{2\pi ik}e^{2\pi ikt}
  -\frac{1}{n}\sum_{k\neq 0}\frac{\hat{f}_k}{2\pi ik}.
  \end{align}
We can compute the Fourier coefficients $\hat{\tilde{f}}_{k}$ from $\hat{f}_k$ by 
assigning $\hat{\tilde{f}}_{k}$ the value $\frac{\hat{f}_{k}}{2 \pi ik}$ for 
$k=-\frac{n}{2}$, \ldots, $\frac{n}{2}-1$, $k \neq 0$ and
assigning $\hat{\tilde{f}}_{0}$ the value $-\sum_{k\neq 0}\frac{\hat{f}_k}{2\pi ik}$.

\subsection{Gaussian filter}
A low-pass filter is commonly used in signal processing to construct a bandlimited 
function. In this paper, we use the Gaussian filter, which is a popular low-pass filter whose 
impulse response is a Gaussian function,
  \begin{align}
  \label{eq:filter1}
  g(x)=a e^{-\pi a^2 x^2},
  \end{align}
where $a$ determines the bandwidth of $g(x)$.

The Gaussian filter $g_0$, \ldots, $g_{n-1}$ is defined to be the IDFT of the sequence 
  \begin{align}
  \label{eq:filter2}
  \hat{g}_k=e^{-{\pi}\frac{k^2}{a^2}}, \qquad k =-\frac{n}{2},\ldots,\frac{n}{2}-1,
  \end{align}
and
coincides with the discrete values of $g(x)$ at the equispaced nodes 
$x_j=\frac{j}{n}$, $j=0$, \ldots, $n-1$.

To filter the Fourier coefficients $\hat{f}_0$, \ldots, $\hat{f}_{n-1}$ in \eqref{eq:rmk1}, we take
the product
  \begin{align}
  \label{eq:filterf1}
  \hat{h}_k=\hat{g}_k\hat{f}_k,\qquad k =-\frac{n}{2},\ldots,\frac{n}{2}-1.
  \end{align}
It is easily to obtain $h_0$, \ldots, $h_{n-1}$ by the IDFT,
  \begin{align}
  \label{eq:filterf2}
  h_j=\frac{1}{n}\sum_{k=0}^{j} g_kf_{j-k},\qquad j =0,\ldots,n-1.
  \end{align}
This can be considered to be a smoothing of $f(x)$ by a convolution of $f(x)$ with the Gaussian function $g(x)$.

Filtering the Chebyshev coefficients $\hat{f}_0$, \ldots, $\hat{f}_{n-1}$
defined in \eqref{eq:chebcoefsdef} is very similar to filtering
the Fourier coefficients, which we describe as follows. Substituting $x=\cos(\theta)$, where $x \in [-1,1]$, into \eqref{eq:cheb1}, we have
  \begin{align}
f(\cos(\theta)) &\approx \sum_{k=0}^{n-1}\hat{f}_kT_k(\cos(\theta))\notag\\ 
&=\sum_{k=0}^{n-1}\hat{f}_k\cos(k \theta),
  \end{align}
where $\theta \in [-\pi,\pi]$.
Letting $\hat{f}_{-k} = \hat{f}_k$, $k=1$, \ldots, $n-1$, we have
  \begin{align}
  \label{eq:filtercheb0}
f(\cos(\theta)) \approx \frac{1}{2}\sum_{k=-n+1}^{n-1}\hat{f}_k  e^{ik\theta} + 
 \frac{1}{2}\hat{f}_0  , \qquad \theta \in [-\pi,\pi]. 
  \end{align}
Hence, by defining $\phi$ by the formula $\theta=2 \pi \phi$,
  \begin{align}
  \label{eq:filtercheb}
f(\cos(2 \pi \phi)) \approx \frac{1}{2}\sum_{k=-n+1}^{n-1}\hat{f}_k e^{2\pi ik\phi}
+\frac{1}{2}\hat{f}_0 , \qquad \phi \in [-\tfrac{1}{2},\tfrac{1}{2}].
  \end{align}
Since \eqref{eq:filtercheb} can be viewed as a Fourier Transform in $\phi$ with the 
Fourier coefficients $\{\hat{f}_k\}$, we follow the equation \eqref{eq:filterf1} to filter
$\{\hat{f}_k\}$, and apply the IDFT to obtain the filtered 
values of $\{{f}_j\}$.
Therefore, $f(x)$ is smoothed by a convolution with the Gaussian function $g(\phi)$ in 
the $\phi$-domain, where $x=\cos(2 \pi \phi)$.

Alternatively, there are other low-pass filters that can be used, such 
as the Butterworth filter (see, for example, Chapter $14$ of~\cite{butterfly}) which resembles the Gaussian filter but is 
flatter in the passband.
The brick-wall filter also preserves signals
with lower frequencies and excludes signals with higher frequencies.
However, after applying the brick-wall filter, the resulting functions tend to 
oscillate at the cutoff frequency (this phenomenon is known as ringing).

\section{The Algorithm}
  \label{sec:algorithm}
In this section, we give an overview of our algorithm for fitting a $C^{\infty}$ curve to pass through a collection of points
$\textbf{C}_{0}$, \ldots, $\textbf{C}_{n}$. We begin with a $C^2$ cubic B\'ezier spline connecting the points 
$\textbf{C}_{0}$, \ldots, $\textbf{C}_{n}$. Given that the curve is at least $C^2$, we interpolate the tangential angle $\theta(t)$ and 
the first derivative of the arc length vector $s'(t)$, which are both $C^1$, using Chebyshev expansions when the curve is
open, or using truncated Fourier series when the curve is closed. We then iteratively filter the coefficients of 
$\theta(t)$ and $s'(t)$ by applying a Gaussian filter, whose bandwidth decreases with each iteration.  
If the curve is closed before filtering, we impose constraints on $\theta(t)$ and 
$s'(t)$ to ensure that 
the curve remains closed. We then reconstruct the curve with the
filtered values of $\theta(t)$ and $s'(t)$ at discretization nodes.

While filtering leads to small discrepancies between the reconstructed curve 
and the points $\textbf{C}_{0}$, \ldots, $\textbf{C}_{n}$,
it also improves the bandwidth of the curve. To fix the discrepancies after
each filtering step, we rotate and rescale the curve
to minimize the total distance between the curve and the points, and add small, smooth perturbations,
which do not negatively affect the smoothness of the curve. We 
stop filtering when the desired bandwidths of the Chebyshev or Fourier approximations to  
$\theta(t)$ and $s'(t)$ are achieved.
This algorithm gives us a $C^{\infty}$ smooth curve that can be represented by a reasonably 
small number of coefficients.

\subsection{Initial Approximation}
  \label{sec:initialinterp}
To initialize our algorithm, we require a $C^2$ curve, the reasons for which are
described in Section \ref{sec:alg_filtering}. 

Given a set of data points $\textbf{C}_{0}$, \ldots, $\textbf{C}_{n} 
\in \mathbb{R}^{2}$, we fit a cubic B\'ezier spline by solving for the intermediate control points $\{\textbf{P}_{i1}\}$ and $\{\textbf{P}_{i2}\}$
described in Section \ref{sec:bezieropen} for an open curve, or in Section \ref{sec:bezierclosed} for a closed curve.
We define the B\'ezier spline $S\colon [0,L] \rightarrow \mathbb{R}^{2}$
connecting all the points $\textbf{C}_{0}$, \ldots, $\textbf{C}_{n}$ by 
  \begin{align}
  \label{eq:initial1}
\textbf{S}(t) = \textbf{B}_{i}(t-i+1),\qquad t \in [0,n]\ \ \text{and}\ \ i =1,\ldots,n,
  \end{align}
if the curve is open, or
  \begin{align}
  \label{eq:initial2}
\textbf{S}(t) = \textbf{B}_{i}(t-i+1),\qquad t \in [0,n+1]\ \ \text{and}\ \ i =1,\ldots,n+1,
  \end{align}
if the curve is closed.

\subsection{Representations of the Curve}
In this section, we denote the curve by
  \begin{align}
\gamma(t)=(x(t),y(t)),
  \end{align}
where $\gamma \colon [0,L] \rightarrow \mathbb{R}^{2}$ 
 is at least $C^2$.

\subsubsection{Representation of an Open Curve}
  \label{sec:repopen}
When the curve is open, we discretize $x(t)$ and $y(t)$ at $N \gg n$
practical Chebyshev nodes $\{t_j\}$ on the interval $[0,L]$ (see formula \eqref{eq:chebnodes})
to obtain $\{x_j\}$ and $\{y_j\}$, where $x_j = x(t_j)$ and $y_j = y(t_j)$.
We use ($N-1$)th order Chebyshev expansions
to approximate $x(t)$ and $y(t)$, constructing the coefficients from 
$\{x_j\}$ and $\{y_j\}$ 
using the Discrete Chebyshev Transform,
and then spectrally differentiate $x(t)$ and $y(t)$ to derive    
the Chebyshev expansions approximating $x'(t)$ and $y'(t)$.
By \eqref{eq:sprime} and \eqref{eq:theta}, we can compute the values of $s'(t)$
and $\theta(t)$ sampled at nodes $\{t_j\}$, and then construct the corresponding 
Chebyshev expansions, 
again using the Discrete Chebyshev Transform.
However, performing the Chebyshev Transform on $\theta(t)$ requires $\theta(t)$ to be continuous,
and as discussed in Section \ref{sec:geometry}, $\theta(t)$ can have jump discontinuities of size $2\pi$.
These can be 
fixed by adding or subtracting multiples of $2\pi$ to $\theta(t)$ wherever 
a discontinuity is detected.

\subsubsection{Representation of a Closed Curve}
  \label{sec:repclosed}
When the curve is closed, 
we discretize $x(t)$ and $y(t)$ at $N \gg n$ equispaced nodes $\{t_j\}$ on
the interval $[0,L]$, where
  \begin{align}
  \label{eq:fournodes}
  {t}_j = \frac{j}{N}L, \qquad j = 0,\ldots,N-1,
  \end{align}
to obtain $\{x_j\}$ and $\{y_j\}$ by $x_j=x(t_j)$ and $y_j=y({t}_j)$.
We then approximate $x(t)$ and $y(t)$ by an $N$-term Fourier series, separately,
and spectrally differentiate $x(t)$ and $y(t)$
to approximate $x'(t)$ and $y'(t)$. 
Following the same procedures in Section \ref{sec:repopen},
we ensure that $\theta(t)$ is continuous, and approximate $s'(t)$ 
by a truncated Fourier series. Recall that, in order to approximate functions
by their Fourier series, the 
functions
must be both smooth and periodic.
The sequence $\{\theta_j\}$, which are the discrete values of $\theta(t)$ at $\{{t}_j\}$,
is not periodic after shifting  
by multiples of $2\pi$ to remove the discontinuities.
Defining $c$ by
  \begin{align}
  c = \theta(n+1)-\theta(0),
  \end{align}
we have that
  \begin{align}
  \tilde{\theta}_j = \theta_j - \frac{c}{L}{t}_j, \qquad {t}_j \in [0,L],
  \end{align}
transforms $\{\theta_j\}$ into a periodic sequence $\{\tilde{\theta}_j\}$ on the interval $[0,L]$, 
which can be approximated by a truncated Fourier series.
To recover the true values of $\{\theta_j\}$ after filtering, we can add $\frac{c}{L}{t_j}$ to 
$\tilde{\theta}_j$. In an abuse of notation, we denote $\{\tilde{\theta}_j\}$ by 
$\{\theta_j\}$ wherever the meaning is clear.

\subsection{Filtering the Curve}
  \label{sec:alg_filtering}
In this section, we describe the process of iteratively filtering $\theta(t)$ and $s'(t)$ using a Gaussian filter.
Given $\gamma(t) \in C^2$, we have $\theta(t) \in C^1$ and $s'(t) \in C^1$. It is known that the decay rate of 
the Chebyshev coefficients or the Fourier coefficients of a $C^1$ function is $O(N^{-\frac{1}{2}})$, where $N$
is the order of the expansion. By iteratively
decreasing the bandwidth of the Gaussian filter, we construct a sequence of bandlimited representations of $\theta(t)$ and $s'(t)$.
The decay rate of the Fourier coefficients or the Chebyshev coefficients in the expansions of $\theta(t)$ and $s'(t)$
increases with each iteration.
This filtering process smooths both the curve itself 
and the parameterization of the curve. 

\subsubsection{Filtering the Open Curve}
Let $\{{t}_j\}$ denote the practical Chebyshev nodes translated to the interval $[0,L]$ (see formula
\eqref{eq:chebnodes}).
Using the Chebyshev expansions of $\theta(t)$ and $s'(t)$ computed in Section \ref{sec:repopen},
we discretize $\theta(t)$ and $s'(t)$ at the points $\{{t}_j\}$ to obtain the sequences $\{\theta_j\}$ and 
$\{s'_j\}$, where $\theta_j = \theta({t}_j) \ \text{and} \ s'_j=s'({t}_j)$.
We filter the Chebyshev coefficients $\{\hat{\theta}_k\}$ of $\theta(t)$,
and $\{\hat{s'}_k\}$ of $s'(t)$ using the Gaussian filter in \eqref{eq:filter2},
and obtain the filtered coefficients $\{\hat{\theta}^{(f)}_k\}$ and $\{\hat{s'}^{(f)}_k\}$,
  \begin{align}
  \hat{\theta}^{(f)}_k=e^{-{\pi}\frac{k^2}{a^2}}\hat{\theta}_{k},\qquad k =0,\ldots,N-1,
  \end{align}
and
  \begin{align}
  \hat{s'}^{(f)}_k=e^{-{\pi}\frac{k^2}{a^2}}\hat{s'}_k,\qquad  k =0,\ldots,N-1.
  \end{align}
Applying the IDFT to $\{\hat{\theta}^{(f)}_k\}$ and $\{\hat{s'}^{(f)}_k\}$, 
we obtain 
  \begin{align}
\theta^{(f)}_j = \sum_{k=0}^{N-1}\hat{\theta}^{(f)}_k T_k(\bar{t}_j), \qquad \bar{t}_j \in [-1,1], 
  \end{align}
where $\bar{t}_j = \frac{2}{L}{t_j}-1$, $t_j \in [0,L]$, and
  \begin{align}
s'^{(f)}_j = \sum_{k=0}^{N-1} \hat{s'}^{(f)}_k T_k(\bar{t}_j).
  \end{align}
We can then use the values of $\{\theta^{(f)}_j\}$ and $\{s'^{(f)}_j\}$ to
recover $\{x^{(f)}_j\}$ and $\{y^{(f)}_j\}$
using \eqref{eq:x} and \eqref{eq:y}.

\subsubsection{Filtering the Closed Curve}
Assume that $\{\theta_j\}$
and $\{s'_j\}$ are the values of $\theta(t)$ and $s'(t)$ discretized
at the equispaced nodes $\{t_j\}$ in \eqref{eq:fournodes},
where $\theta_j = \theta({t}_j) \ \text{and} \ s'_j=s'({t}_j)$. We apply the DFT to derive the 
Fourier coefficients $\{\hat{\theta}_k\}$ of $\theta(t)$  and $\{\hat{s'}_k\}$ of $s'(t)$. 
Using the Gaussian filter, we filter the 
Fourier coefficients $\{\hat{\theta}_k\}$ and $\{\hat{s'}_k\}$ to obtain the filtered Fourier coeffcients
$\{\hat{\theta}^{(f)}_k\}$ and $\{\hat{s'}^{(f)}_k\}$, 
  \begin{align}
\hat{\theta}^{(f)}_k=e^{-{\pi}\frac{k^2}{a^2}}\hat{\theta}_{k},\qquad k =-\frac{N}{2},\ldots,\frac{N}{2}-1,
  \end{align}
and
  \begin{align}
\hat{s'}^{(f)}_k=e^{-{\pi}\frac{k^2}{a^2}}\hat{s'}_{k},\qquad k =-\frac{N}{2},\ldots,\frac{N}{2}-1.
  \end{align}
We recover the filtered sequences $\{\theta^{(f)}_j\}$ and $\{s'^{(f)}_j\}$ by applying the IDFT to the filtered Fourier coefficients $\{\hat{\theta}^{(f)}_{k}\}$
and $\{\hat{s'}^{(f)}_k\}$,
  \begin{align}
\theta^{(f)}_j= \sum_{k=-\frac{N}{2}}^{\frac{N}{2}-1}\hat{\theta}^{(f)}_k e^{\frac{2\pi i}{N}kj}+\frac{c}{L}t_j, \qquad j = 0, \ldots, N-1,
  \end{align}
and
  \begin{align}
s'^{(f)}_j= \sum_{k=-\frac{N}{2}}^{\frac{N}{2}-1} \hat{s'}^{(f)}_k e^{\frac{2\pi i}{N}kj}, \qquad j = 0, \ldots, N-1.
  \end{align}

Similarly, the curve can be reconstructed from
$\{\theta^{(f)}_j\}$ and $\{s'^{(f)}_j\}$, using equations \eqref{eq:x} and \eqref{eq:y}.

\subsection{Closing the Curve}
Applying a filter to $\theta(t)$ and $s'(t)$ for a closed curve, in general, makes the curve become open. To close the curve,
we require that
  \begin{align}
  \label{eq:close1} 
\int_{0}^{L} s'(t)\cos{\theta(t)}\, dt = 0,
  \end{align}
and 
  \begin{align}
  \label{eq:close2} 
\int_{0}^{L} s'(t)\sin{\theta(t)}\, dt = 0.
  \end{align}
The process of orthogonalizing $s'(t)$ to $\cos{\theta(t)}$ and $\sin{\theta(t)}$ using the
trapezoidal rule is as follows.
Supposing that we have 
the values $\{s'_j\}$, $\{\cos{\theta_j}\}$ and $\{\sin{\theta_j}\}$ 
of $s'(t)$, $\cos{\theta(t)}$ and $\sin{\theta(t)}$ sampled at 
the points $\{{t}_j\}$ defined in \eqref{eq:fournodes}. We ensure that $\{s'_j\}$ is orthogonal to $\cos{\theta_j}$
by setting $\{s'_j\}$ to the values
  \begin{align}
s'_j-\cos{\theta_j}\frac{\frac{1}{N}\sum_{j=0}^{N-1}s'_j\cos{\theta_j}}{\frac{1}{N}\sum_{j=0}^{N-1}{{\cos}^2{\theta_j}}},
\qquad j=0, \ldots, N-1.
  \end{align}
We let $\{\lambda_j\}$ be the vector defined by the formula
  \begin{align}
\lambda_j = \sin{\theta_j}-\cos{\theta_j}\frac{\frac{1}{N}\sum_{j=0}^{N-1}\sin{\theta_j}\cos{\theta_j}}{\frac{1}{N}\sum_{j=0}^{N-1}{{\cos}^2{\theta_j}}},
\qquad j=0, \ldots, N-1.
  \end{align}
Finally, we orthogonalize $\{s'_j\}$ to $\{\lambda_j\}$ by setting $\{s'_j\}$ 
to the values 
  \begin{align}
s'_j-\lambda_j\frac{\frac{1}{N}\sum_{j=0}^{N-1}s'_j\lambda_j}{\frac{1}{N}\sum_{j=0}^{N-1}{{\lambda}^2_j}},
\qquad j=0, \ldots, N-1.
  \end{align}
The sequence $\{s'_j\}$ is now orthogonal to both $\{\cos{\theta_j}\}$ and $\{\sin{\theta_j}\}$. Thus, the conditions 
\eqref{eq:close1} and \eqref{eq:close2} are satisfied to within the accuracy of the trapezoidal rule.

\subsection{Repositioning the Curve}
In general, the curve will not pass through the original data points after filtering. Moreover, filtering $\theta(t)$
and $s'(t)$ changes the tangential vector
$T(t)$, which results in changes in the orientation and position of the curve. In this section, we describe how to rotate the reconstructed curve
so that the sum of squares of the distances between the curve and the original data points is minimized.

Given the original data points $\{\textbf{C}_{i}\}$, where $\textbf{C}_{i}=(\textbf{C}_{ix}, \textbf{C}_{iy}), i=0,\ldots,n$, we find 
$\tilde{t}_{0}$, \ldots, $\tilde{t}_{n}$ $\in$ $[0,L]$ such that,
if $(\tilde{x}_i$, $\tilde{y}_i) = (x(\tilde{t}_i)$, $y(\tilde{t}_i))$,
then $(\tilde{x}_i,\tilde{y}_i)$ is the closest point on the curve to $(\textbf{C}_{ix}, \textbf{C}_{iy})$ for $i=0$, \ldots, $n$. 
We determine $\tilde{t}_{0},\dots,\tilde{t}_{n}$ only once, described in Remark~\ref{rem:perturb1}.
Suppose that $\{\phi_i\}$ are the 
values of the angle between $\{(\tilde{x}_i$, $\tilde{y}_i)\}$ and $(\bar{x},\bar{y})$,
and that $\{r_i\}$ are the distances between $\{(\tilde{x}_i,\tilde{y}_i)\}$ and $(\bar{x},\bar{y})$, where 
$(\bar{x},\bar{y})$ is the center of all the closest points $\{(\tilde{x}_i,\tilde{y}_i)\}$.
We shift the center of all the closest points $\{(\tilde{x}_i,\tilde{y}_i)\}$ by $(\Delta x,\Delta y)$,
and rotate the curve by an angle of $\psi$ around the
center. Observed that the sum of squares of the distances between the closest points
and the original data points is given by 
  \begin{align}
  \label{eq:repos2} 
\hspace {-2em} f(\psi,\Delta x, \Delta y) = \sum_{i=0}^{n}&(\bar{x} + \Delta x + r_i\cos{(\phi_i+\psi)}-\textbf{C}_{ix})^2+ \notag\\
&(\bar{y} + \Delta y +r_i\sin{(\phi_i+\psi)}-\textbf{C}_{iy})^2,
  \end{align}
where $(\bar{x},\bar{y})$ is the average of $\{(\tilde{x}_i,\tilde{y}_i)\}$. 
We use Newton's method to obtain the values of $\psi$, $\Delta x$ and $\Delta y$ which minimize $f(\psi, \Delta x, \Delta y)$.

\begin{remark}
  \label{rem:repos1}
One might also think to rescale the curve by multiplying $\{r_i\}$ by a constant $c$, since filtering $s'(t)$ changes the length of the curve.
However, rescaling the curve distorts the structure of the closest points $(\tilde{x}_i,\tilde{y}_i)$ on the curve.
Large perturbations, as described in Section~\ref{sec:perturb},
are sometimes needed as a result, 
and therefore the smoothness of the curve after adding perturbations can be reduced. 
\end{remark}

Notice that each point on the curve is, in some sense, equivalent.
The procedure of repositioning
ensures that the resulting curve is affine invariant.

\subsection{Adding Perturbations to the Curve}
  \label{sec:perturb}
Since the curve does not pass through the original data points $\{\textbf{C}_{i}\}$ after filtering $\theta(t)$ and $s'(t)$, we introduce a 
set of Gaussian functions, $\{g_i(t)\}$, which we use as 
smooth perturbations that can be added to the curve to ensure that the curve passes through the points $\{\textbf{C}_{i}\}$.
We define $g_i(t)$ by 
  \begin{align}
  g_i(t) = e^{-\sigma_i \bigl(\frac{t-\tilde{t}_i}{L}\bigr)^2}, \qquad i = 0, \ldots, n,
  \end{align}
for $t \in [0,L]$, where $\tilde{t}_i$ is the curve parameter of the closest point
$(\tilde{x}_i$, $\tilde{y}_i) = (x(\tilde{t}_i)$, $y(\tilde{t}_i))$ to $C_i$, and $\sigma_i$ determines the bandwidth of the 
perturbation. 
When the curve is closed, $g_i(t)$ is modified to be a periodic function with period $L$, given by the formula
  \begin{align}
g_i(t) = \sum_{k=-\infty}^{\infty}
 e^{-\sigma_i \bigl(\frac{t-\tilde{t}_i}{L}+ k \bigr)^2},\qquad i = 0,\ldots, n.
  \end{align}
It is obvious that $g_i(t)=g_i(t+L)$. 
We construct $\{(\bar{x}_j,\bar{y}_j)\}$ from $\{(x_j,y_j)\}$ by adding 
$g_i(t)$ at the discretized points $\{t_j\}$, 
  \begin{align}
  \bar{x}_j=x_j+\sum_{i=0}^{n}c_{ix}g_i(t_j), \qquad j=0,\ldots,N-1,
  \end{align}
and
  \begin{align}
  \bar{y}_j=y_j+\sum_{i=0}^{n}c_{iy}g_i(t_j),\qquad j=0,\ldots,N-1,
  \end{align}
where $\{c_{ix}\}$ and $\{c_{iy}\}$ are the coefficients of perturbations in $x$ and $y$, separately, 
which are reasonably small since the curve
is filtered slightly at each iteration. Let $\{(\tilde{\bar{x}}_j,\tilde{\bar{y}}_j)\}$ denote 
the points on the perturbed curve corresponding to $\tilde{t}_{0},\dots,\tilde{t}_{n}$.
We require 
  \begin{align}
  \tilde{\bar{x}}_i=\textbf{C}_{ix}, \qquad i=0,\ldots,n,
  \end{align}
and
  \begin{align}
  \tilde{\bar{y}}_i=\textbf{C}_{iy}, \qquad i=0,\ldots,n,
  \end{align}
and solve two linear systems of $n+1$ equations to compute the values of $\{c_{ix}\}$ and $\{c_{iy}\}$. 
We observe that, since the perturbations $g_i(t)$ are Gaussians, they are each, to finite precision, compactly
supported. Thus, the linear system that we solve is effectively banded, and the number of bands is determined by 
$\min_{i} \sigma_i$. An $O(n+1)$ solver can be used to 
speed up the computations.
\begin{remark}
  \label{rem:perturb1}
We only calculate $\{\tilde{t}_i\}$ once, at the first iteration before filtering, and use the same set of $\{\tilde{t}_i\}$ at each iteration.
Although it seems more natural to recalculate $\{\tilde{t}_i\}$ at each iteration, so that the discrepancies are fixed by smaller perturbations, 
the resulting perturbations are always orthogonal to the curve.
The effect of the changes in the length
of the curve due to filtering can not be eliminated by adding such perturbations, 
with the effect that the length of the curve grows if the points $\{\tilde{t}_i\}$
are calculated at each iteration.
By using the same set of closest points for all iterations, 
the perturbations can be oblique, which results in nice control over the total length of the curve during the filtering process. 
\end{remark}

\subsection{The Termination Criterion of the Algorithm}
Since the bandwidths of the coefficients of $\theta(t)$ and $s'(t)$ are reduced
at each iteration, and adding small, smooth perturbations has a negligible effect 
on the bandwidth of
the curve, one can expect to achieve 
the desired bandwidth of the representations 
of $\theta(t)$ and $s'(t)$ by iteratively filtering the coefficients. 
However, we note that there is a minimum number of coefficients that are necessary to represent a curve, as 
determined by the sample data points. When fewer than this number of 
coefficients are used, the curve reconstructed by these overfiltered coefficients
may deviate drastically from the sample data points. The resulting large perturbations
required to fix the discrepancies can harm the smoothness
of the curve.
The purpose of this section is to set up a termination criterion,
so that the algorithm will terminate if the coefficients of $\theta(t)$ and $s'(t)$, 
beyond a user-specified number of terms,
are filtered to zero, to the requested accuracy.

We denote the desired accuracy of the approximation by $\epsilon$, which 
is often set to be machine precision, and the number of
coefficients representing the curve that are larger than $\epsilon$ by 
$n_{\text{coefs}}$. 
Due to the potentially large condition number of 
spectral differentiation, some accuracy is  
lost when computing the coefficients of $x'(t)$ and $y'(t)$, and thus $\theta(t)$ and $s'(t)$, at each iteration.
Thus, we measure thresholds for the coefficients of $\theta(t)$ and $s'(t)$, below 
which they are considered to be zero, and denote them by $\delta_{\theta}$
and $\delta_{s'}$. 
We consider first the open curve case. Since the condition number of the Chebyshev
differentiation matrix 
is bounded by approximately $N^{\frac{3}{2}}$, where $N$ is the number of coefficients, the error 
induced by differentiating $x(t)$ and $y(t)$ is approximately
  \begin{align}
  \label{eq:delta0} 
  & \ \epsilon N^{\frac{3}{2}}\sqrt{\norm{x(t)}_{L^2[0,L]}^2+\norm{y(t)}_{L^2[0,L]}^2} \\
  \approx & \ \epsilon N^{\frac{3}{2}}\sqrt{\sum_j {x^2_jw_j}+\sum_j {y^2_jw_j}},
  \end{align}
where $\{w_j\}$ denotes
the Chebyshev weights on [$0,L$].
Considering the way $\theta(t)$ is calculated, the error
in $\theta(t)$ is proportional to the error in $x'(t)$ and $y'(t)$, divided by
the norm of the tangential vector $(x'(t),y'(t))$. Thus, we set 
  \begin{align}
  \label{eq:delta1} 
  \delta_{\theta} &= \epsilon N^{\frac{3}{2}}\sqrt{\norm{x(t)}_{L^2[0,L]}^2+\norm{y(t)}_{L^2[0,L]}^2}
  \cdot \Bnorm{\frac{1}{\sqrt{x'(t)^2+y'(t)^2}}}_{L^{\infty}[0,L]} \notag \\
  &\approx \frac{\epsilon N^{\frac{3}{2}}\sqrt{\sum_j {x^2_{j} w_j}+\sum_j {y^2_{j} w_j}}}
  {\min {\sqrt{x'^2_jw_j+y'^2_jw_j}}},
  \end{align}
where $x'_i$, $y'_i$ are the discretized values of $x'(t)$, $y'(t)$.
Similarly, 
the error in $s'(t)$
is proportional to the error in $x'(t)$ and $y'(t)$.
Thus, we set 
  \begin{align}
  \label{eq:delta2} 
  \delta_{s'} &= \epsilon N^{\frac{3}{2}}\sqrt{\norm{x(t)}_{L^2[0,L]}^2+\norm{y(t)}_{L^2[0,L]}^2}, \notag \\
  &\approx \epsilon N^{\frac{3}{2}}\sqrt{\sum_j {x^2_{j} w_j}+\sum_j {y^2_{j} w_j}}.
  \end{align}

The thresholds $\delta_{\theta}$ and $\delta_{s'}$ for the closed curve case
are almost identical, except that the condition number of spectral differentiation
matrix is approximately $N$, where $N$ is the number of coefficients, 
from which it follows that $N^{\frac{3}{2}}$ is replaced by $N$, and 
the weights $w_j$ 
are replaced by $\frac{L}{N}$.

Suppose that we have the desired accuracy of the approximation, $\epsilon$, 
the threshold, $\delta_{\theta}$, and the number of coefficients larger than $\epsilon$, $n_{\text{coefs}}$.
We consider first the coefficients of $\theta(t)$.
Our goal is to 
determine the number of coefficients, $n_{\text{coefs}}^{\delta_{\theta}}$, that we expect
to be larger than $\delta_{\theta}$, when there are only $n_{\text{coefs}}$ terms
larger than $\norm{\hat{\theta}}_{\infty} \epsilon$. 
In order to approximate $n_{\text{coefs}}^{\delta_{\theta}}$, 
we assume that the coefficients $\{\hat{\theta}_k \}$ decay exponentially, like
$\norm{\hat{\theta}}_{\infty} e^{-Ck}$, from the maximum value
$\norm{\hat{\theta}}_{\infty}$
to $\norm{\hat{\theta}}_{\infty} \epsilon$.
This implies that
$C=\frac{\log{(1/{\epsilon})}}{n_{\text{coefs}}}$.
Thus, 
  \begin{align}
  e^{-\log{(1/{\epsilon}}) \frac{n_{\text{coefs}}^{\delta_{\theta}}}
  {n_{\text{coefs}}}} =\delta_{\theta},
  \end{align}
so,
  \begin{align}
n_{\text{coefs}}^{\delta_{\theta}} = n_{\text{coefs}}
\frac{\log{(1/{\delta_{\theta}})}}{ \log{(1/{\epsilon})}}.
  \end{align}

We compute $n_{\text{coefs}}^{\delta_{s'}}$ in exactly the same way.
At each iteration, 
if only $n_{\text{coefs}}^{\delta_{\theta}}$ and $n_{\text{coefs}}^{\delta_{s'}}$
numbers of terms
are larger than $\delta_{\theta}$ and 
$\delta_{s'}$, respectively, then the algorithm terminates.
Eventually, $n_{\text{coefs}}$ coefficients are returned to 
the user to represent the curve, up to the precision 
$\epsilon$.

\begin{remark}
Since the values of $\delta_{\theta}$, $\delta_{s'}$, 
$n_{\text{coefs}}^{\delta_{\theta}}$ and $n_{\text{coefs}}^{\delta_{s'}}$
are fairly consistent in each iteration,
we only calculate these values 
once, at the first iteration. 
\end{remark}

\subsection{Summary and Cost of the Algorithm}

The algorithm can be summarized as follows:

\begin{enumerate}
\item Given $n+1$ points $\textbf{C}_0$, \ldots, $\textbf{C}_n$, 
fit a $C^2$ B\'ezier spline to connect the points.  

\item Discretize the curve at $N \gg n+1$ Chebyshev nodes if the curve is open, or $N \gg n+1$  
equispaced nodes if the curve is closed, and compute $\{\theta_j\}$
and $\{s'_j\}$.
\end{enumerate}

Repeat the steps $3$, \ldots, $9$ until a $C^{\infty}$ smooth curve
can be represented by the requested number of coefficients, $n_{\text{coefs}}$:  
\begin{enumerate}
\setcounter{enumi}{2}
\item Obtain the Chebyshev coefficients
or the Fourier coefficients of $\{\theta_j\}$ 
and $\{s'_j\}$.

\item Determine the number of coefficients of $\{\theta_j\}$ and $\{s'_j\}$
larger than $\delta_{\theta}$ and $\delta_{s'}$. If there are fewer than
$n_{\text{coefs}}^{\delta_{\theta}}$ and $n_{\text{coefs}}^{\delta_{s'}}$, respectively,
then return the first $n_{\text{coefs}}$ coefficients of $x(t)$ and $y(t)$.

\item Apply the filter to the coefficients 
of $\{\theta_j\}$ and $\{s'_j\}$ to compute the filtered values of
$\{\theta_j\}$ and $\{s'_j\}$.
\item In the case of a closed curve, modify $\{s'_j\}$ to satisfy 
the constraints \eqref{eq:close1} and \eqref{eq:close2} in order to close the curve after filtering.
\item Reconstruct the curve from $\{\theta_j\}$ and $\{s'_j\}$ by equations 
\eqref{eq:x} and \eqref{eq:y}.
\item Rotate the curve to minimize the sum of squares of the distances between the curve and the points 
$\textbf{C}_0$, \ldots, $\textbf{C}_n$.
\item Add smooth Gaussian perturbations to make the curve pass through 
the points $\textbf{C}_0$, \ldots, $\textbf{C}_n$.
\end{enumerate}

Solving for the control points of the B\'ezier spline in Step $1$ costs $O(n+1)$ operations, and discretizing the 
spline at N points in Step $2$ costs $O(N)$ operations. Step $3$ involves 
spectral differentiation and the Discrete Chebyshev
Transform in the open curve case, or the DFT in the closed curve case, where the Discrete Chebyshev
Transform can be replaced by the Fast Chebyshev Transform and the DFT can be replaced by
the FFT. The cost of step $3$ is thus reduced to $O(N \log{N})$. 
Checking the termination condition in Step $4$
costs approximately $O(N)$ operations.
Applying the filter and 
reconstructing $\{\theta_j\}$ and $\{s'_j\}$ in Step $5$ has the same cost as applying the inverse 
Fast Chebyshev Transform or the IFFT, which costs $O(N \log{N})$ operations.
If the curve is closed, we must modify $\{s'_j\}$ so that the curve remains closed. 
The cost of closing the curve by looping through $\{s'_j\}$ in Step $6$ is $O(N)$.  
Step $7$ involves spectral integration, and the inverse Fast Chebyshev Transform in the open curve case, or 
the IFFT in the closed curve case, which has
the same $O(N \log{N})$ cost as Step $3$. The cost of using Newton's method to rotate the curve in 
Step $8$ is $O(n+1)$,
and the cost of solving for the coefficients of the smooth perturbations added 
to the curve in Step $9$ is $O(n+1)$.
The total cost is thus $O(N \log{N})$ per iteration.

\section{Numerical Results}
  \label{sec:numerical}
In this section, we demonstrate the performance of our algorithm with several
numerical examples, 
and present both the analytic curves produced by
the algorithm and filtered coefficients of the functions $\theta(t)$ and $s'(t)$ representing the curves,
where $\theta(t)$ is the tangential angle and $s'(t)$ is the first derivative of the arc length. 
We implemented our algorithm in Fortran 77, and compiled it using the Gfortran 
Compiler, version 9.4.0, with -O3 flag. All experiments were conducted on a 
laptop with 16 GB of RAM and an Intel $11$th Gen Core i7-1185G7 CPU. Furthermore, we use FFTW library
(see~\cite{fftw}) for the implementations of the FFT and the Fast Cosine Transform. The latter 
is used to implement the Fast Chebyshev Transform.

The following variables appear in this section: 
\begin{itemize}
\item[$-$] $N$: the number of discretization nodes.
\item[$-$] $n$: the number of sample data points.
\item[$-$] $n_{\text{iters}}$: the maximum number of iterations.
\item[$-$] $n_{\text{stop}}$: the number of iterations needed for the algorithm to terminate. 
\item[$-$]  $h_{\text{filter}}$: the proportion of the coefficients that are filtered to zero at each iteration.
\item[$-$] $\epsilon$: the desired accuracy of the approximation to the curve. 
As $\epsilon$ is dependent on the size of the curve, for consistency,
we scale the sample data points, so that either the width or height of
the collection of data points,
whichever is closer to $1$, is $1$.
\item[$-$] $n_{\text{coefs}}$: the requested number of the coefficients representing the curve to precision $\epsilon$.
\item[$-$] $n_{\text{bands}}$: the bandwidth of the matrix describing the effect 
of the Gaussian perturbations centered at each sample point. 
\item[$-$] $x'_{\text{left}}$, $y'_{\text{left}}$: the derivative of the initial curve specified at the left end point,
in the $x$ coordinate and $y$ coordinate separately. 
This variable only exists in the open curve case. Notice that the filtering process can potentially
alter the value of this variable.
\item[$-$] $x'_{\text{right}}$, $y'_{\text{right}}$: the derivative of the initial curve specified at the right end point,
in the $x$ coordinate and $y$ coordinate separately.
This variable only exists in the open curve case. Notice that the filtering process can potentially
alter the value of this variable.
\item[$-$] $E_{\text{samp}}$: the maximum $l_2$ norm of the distance between the curve,
defined by $n_{\text{coefs}}$ Chebyshev or Fourier coefficients, and 
the sample data points.
\end{itemize}
While there is no strict rule on how to choose these variables, we assume that the users  
pick a reasonable combination of inputs, so that the algorithm terminates before reaching
the maximum number of iterations, $n_{\text{iters}}$.

\subsection{Open Curve Examples}
\begin{figure}[!h]
  \centering
  \begin{subfigure}[b]{0.28\linewidth}
    \includegraphics[width=\linewidth]{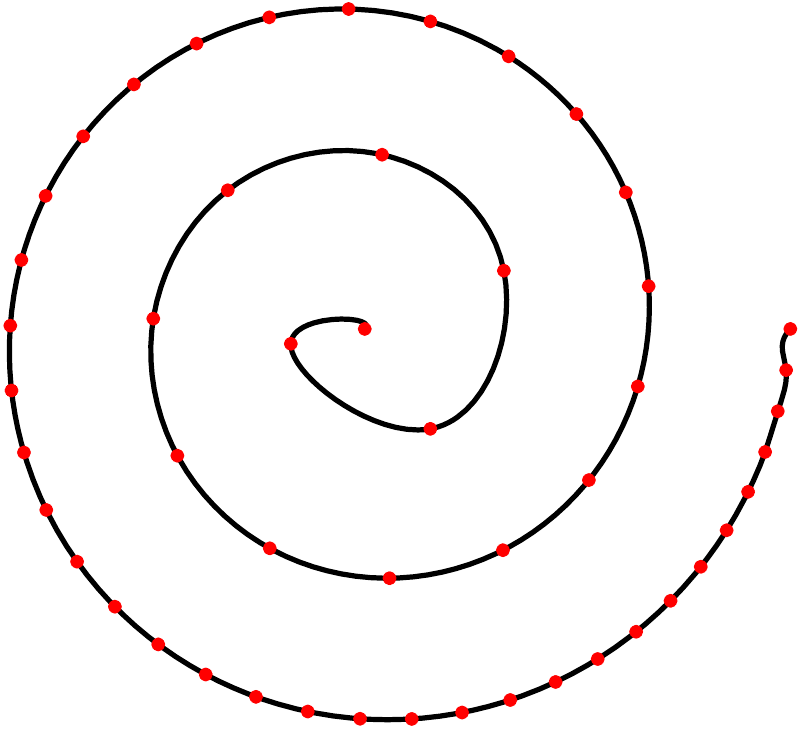}
    \caption{The curve before smoothing}
  \end{subfigure}\hspace{8mm}
  \begin{subfigure}[b]{0.28\linewidth}
    \includegraphics[width=\linewidth]{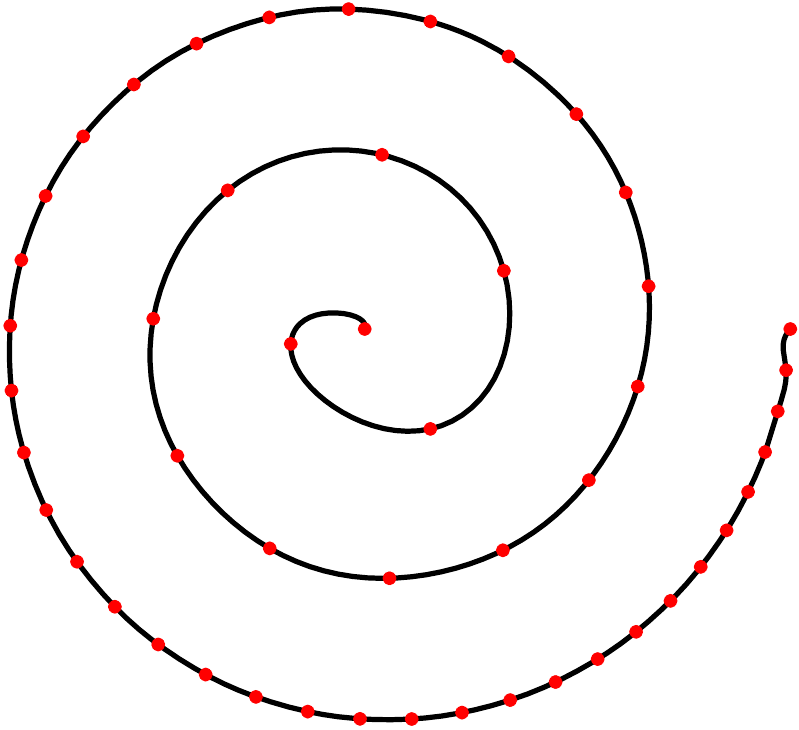}
    \caption{The curve after smoothing}
  \end{subfigure}
  \caption{The result of algorithm applied to \eqref{eq:openex1}. The red dots 
 mark the sample points.}
  \label{fig:open1}
\end{figure}

\begin{figure}[!h]
  \centering
  \begin{subfigure}[b]{0.45\linewidth}
    \includegraphics[width=\linewidth]{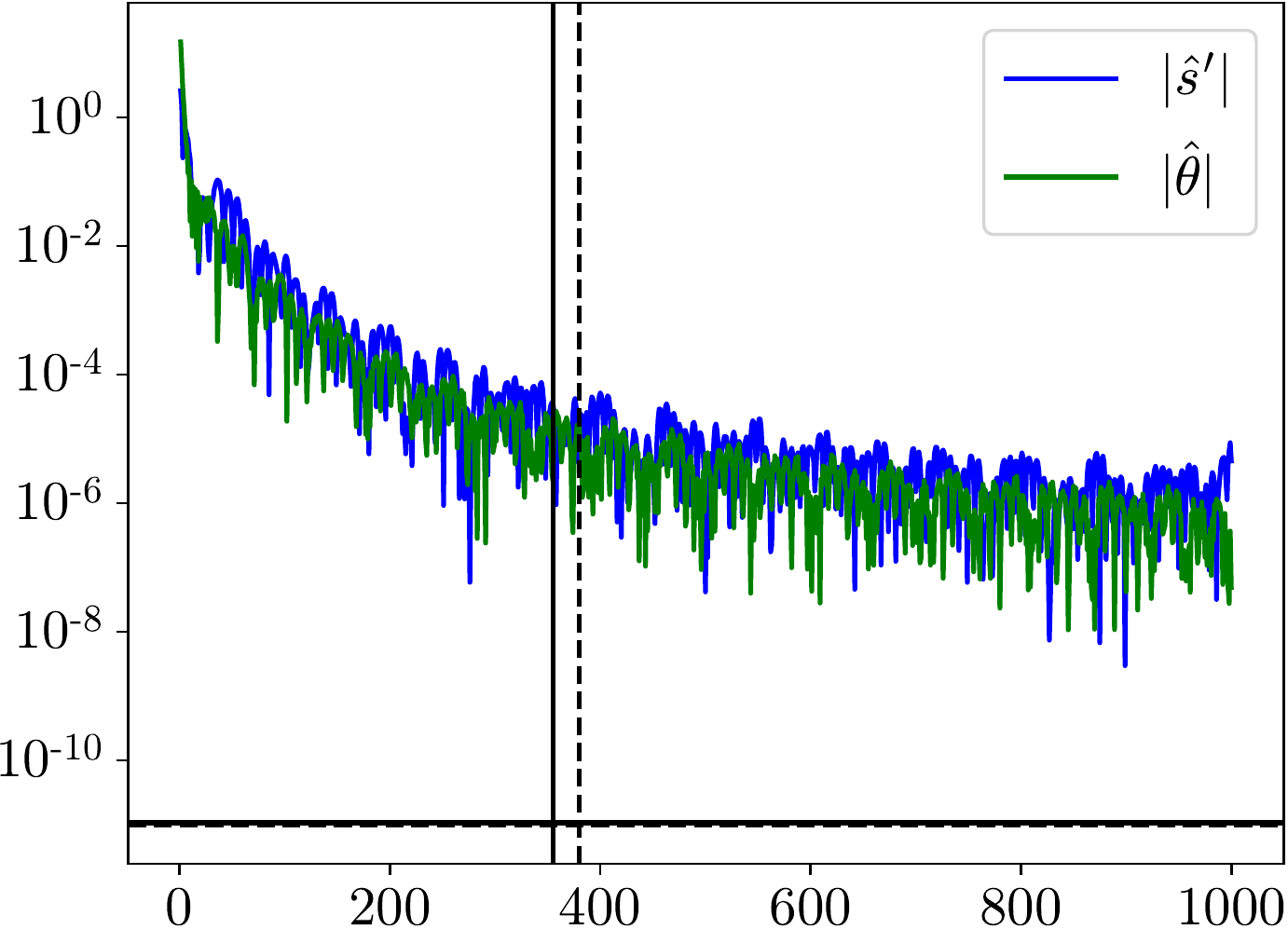}
    \caption{Before filtering}
  \end{subfigure}
  \begin{subfigure}[b]{0.45\linewidth}
    \includegraphics[width=\linewidth]{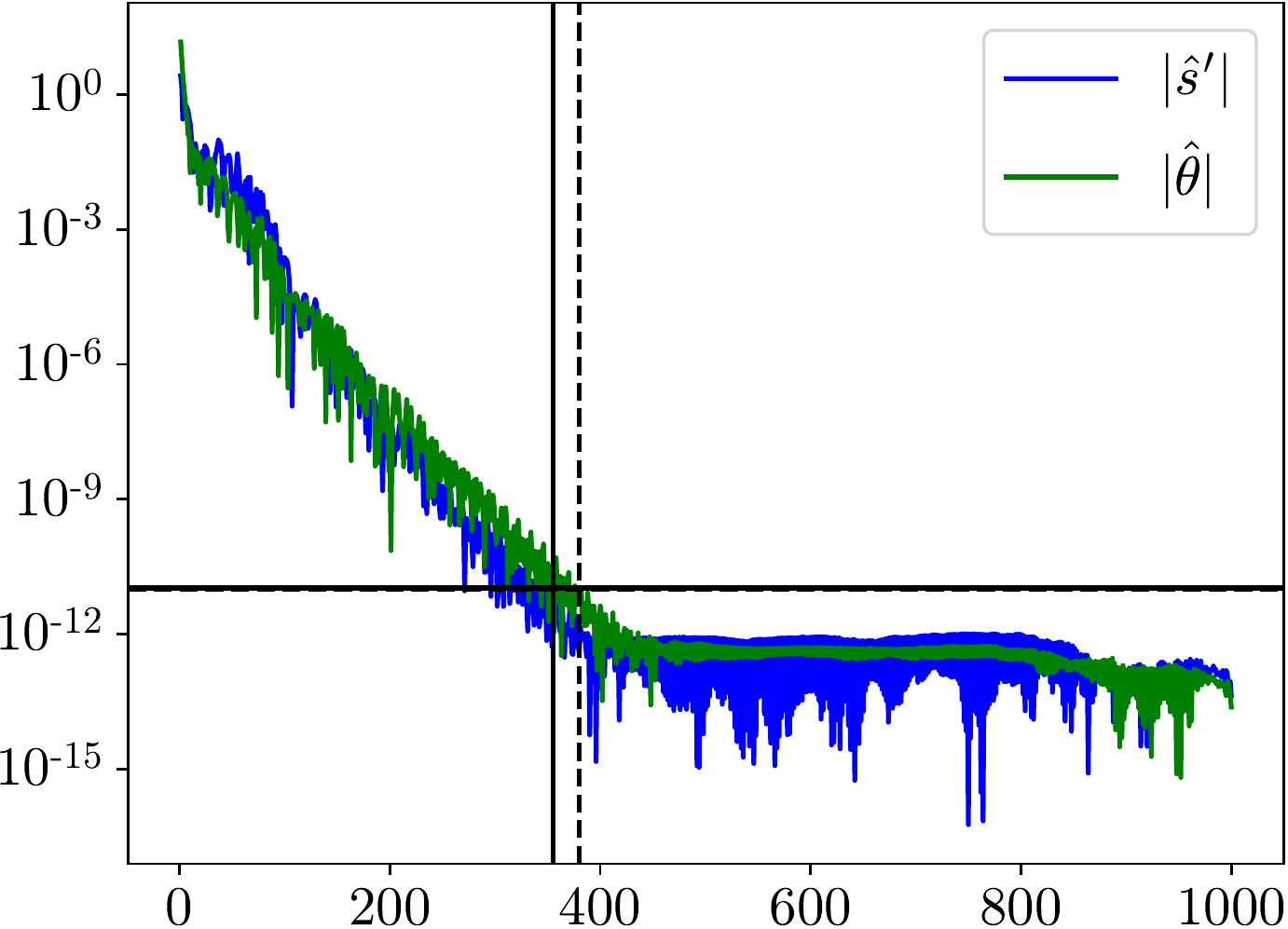}
    \caption{After filtering}
  \end{subfigure}
  \caption{Chebyshev coefficients of $s'(t)$ and $\theta(t)$ corresponding to 
  Figure~\ref{fig:open1}. The value of $\delta_{s'}$ is indicated by a horizontal solid line 
  and the value of $\delta_{\theta}$ is indicated by a horizontal dashed line.
  The $355$th coefficients of $s'(t)$ decays to $\delta_{s'}$, indicated by a 
  vertical solid line.
  The $380$th coefficients of $\theta(t)$ decays to $\delta_{\theta}$,
  indicated by a vertical dashed line.}
  \label{fig:open1c}
\end{figure}

\begin{figure}[!h]
  \centering
  \begin{subfigure}[b]{0.45\linewidth}
    \includegraphics[width=\linewidth]{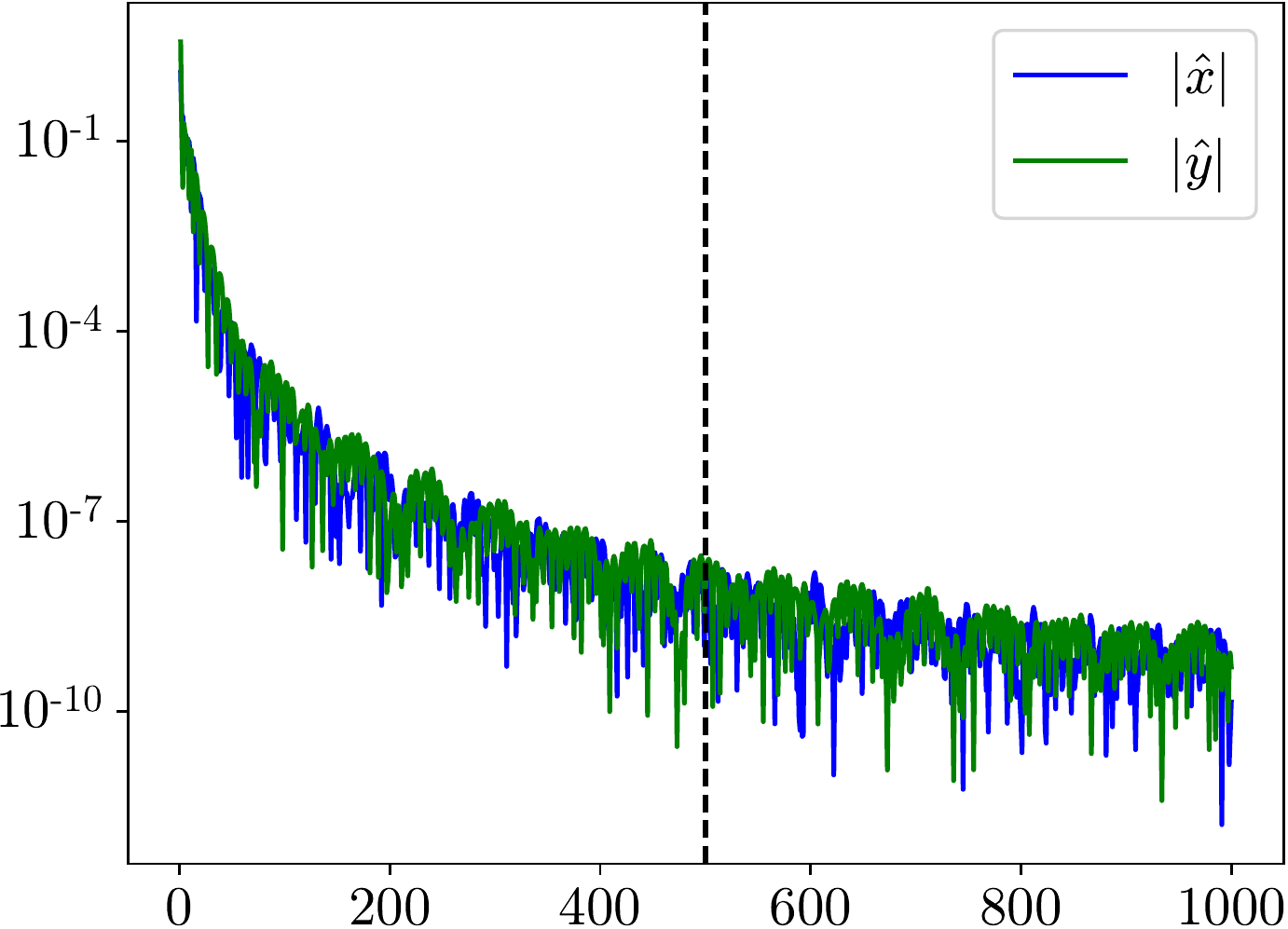}
    \caption{Chebyshev coefficients of the initial curve}
  \end{subfigure}
  \begin{subfigure}[b]{0.45\linewidth}
    \includegraphics[width=\linewidth]{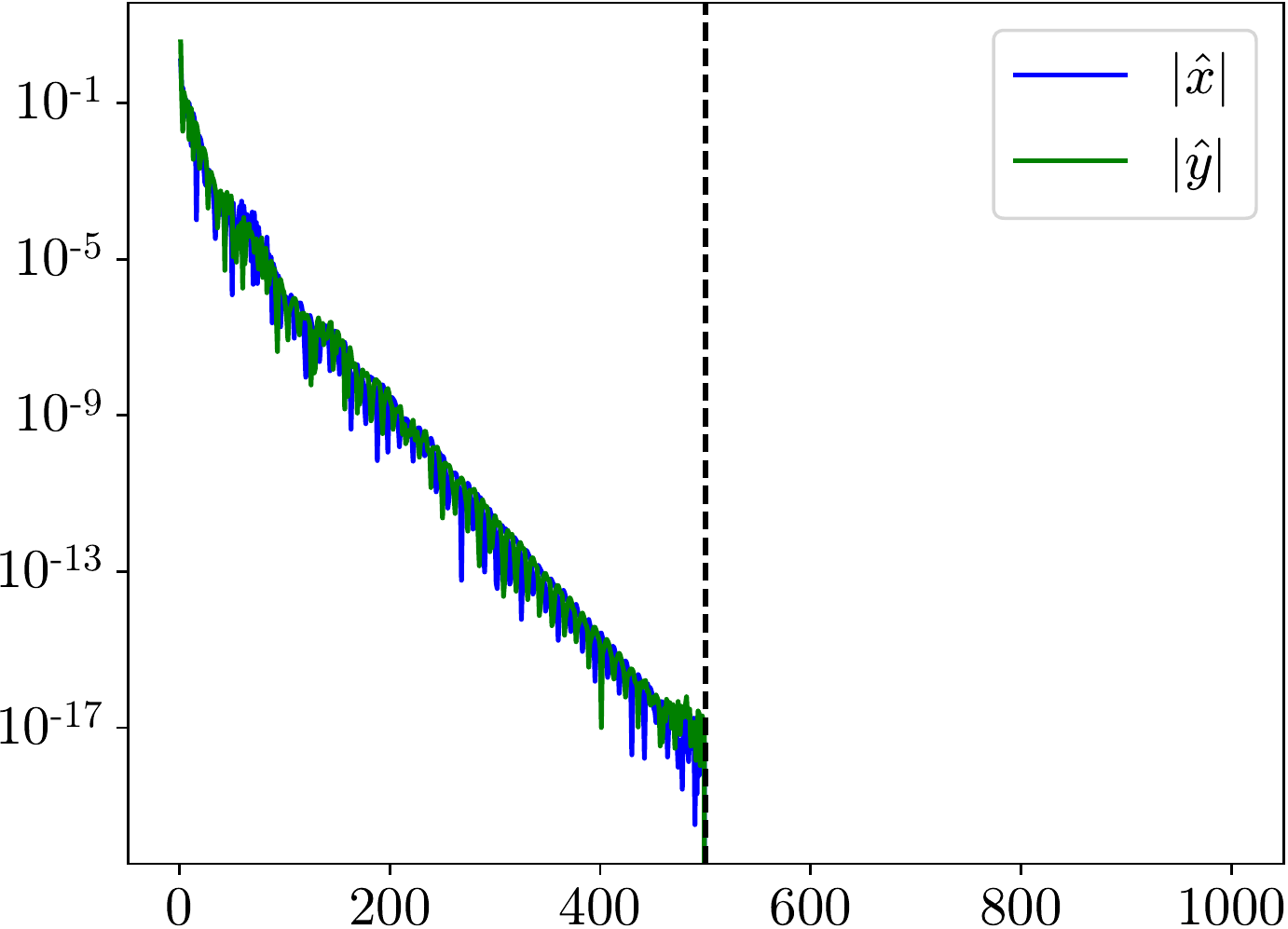}
    \caption{Chebyshev coefficients of the final curve}
  \end{subfigure}
  \caption{Chebyshev coefficients of $x(t)$ and $y(t)$ corresponding to
  Figure~\ref{fig:open1}.
  The value of $n_{\text{coefs}}$ is indicated by a vertical dashed line.} 
  \label{fig:open1c2}
\end{figure}

We sample some points from a spiral with the polar representation 
$(r(t)\cos{\varphi(t)}$, $r(t)\sin{\varphi(t)})$, where
  \begin{align}
  \label{eq:openex1} 
  \varphi(t)&=\frac{6 \pi}{\log{2}} \log{t}, \notag\\
  r(t)&=\varphi(t),  
  \end{align}
with $t \in [1,2]$, and construct 
the initial B\'ezier 
spline passing through the data points, as shown in Figure~\ref{fig:open1}(a).
The sample data points are scaled so that their width is $1$.
We set
$N=1000$, $n=50$, $x'_{\text{left}}=0.05$, $y'_{\text{left}}=0.05$,
$x'_{\text{right}}=0.05$,
$y'_{\text{right}}=0.05$,  
$n_{\text{iters}}=60$
$h_{\text{filter}}=\frac{1}{25}$,
$\epsilon=10^{-16}$, $n_{\text{coefs}}=500$,
$n_{\text{bands}}=8$. After 
$n_{\text{stop}}=16$ iterations, the algorithm terminates
and returns a curve represented
by only $500$ Chebyshev coefficients. 
We display the Chebyshev coefficients that are necessary to represent both the initial
and final curve in Figure~\ref{fig:open1c2}.
We can see that the shape of the final curve in 
Figure~\ref{fig:open1}(b) is smoother,
especially at the center of the spiral. 
Moreover, the resulting curve
curve passes through
the sample data points with an error of $E_{\text{samp}}=0.11548 \e{-13}$. 
The magnitudes
of the Chebyshev coefficients of $s'(t)$ and $\theta(t)$ before and after
filtering are displayed in Figure~\ref{fig:open1c}.

\begin{figure}[!h]
  \centering
  \begin{subfigure}[b]{0.25\linewidth}
    \includegraphics[width=\linewidth]{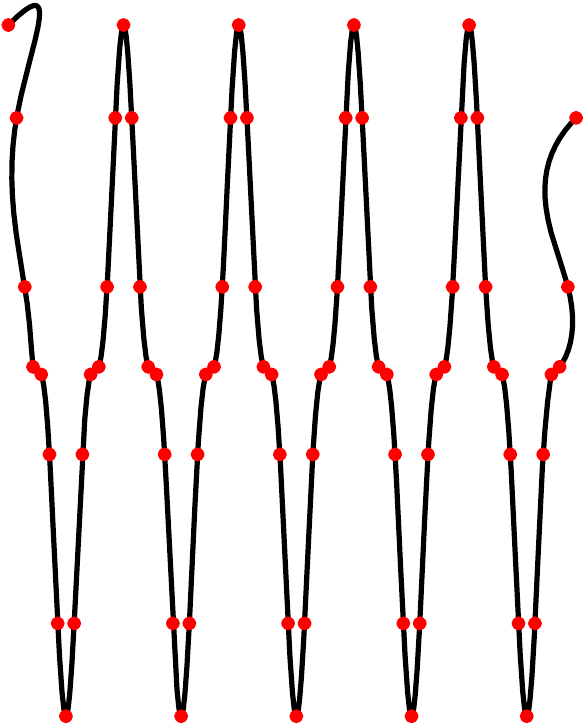}
    \caption{The curve before smoothing}
  \end{subfigure}\hspace{8mm}
  \begin{subfigure}[b]{0.25\linewidth}
    \includegraphics[width=\linewidth]{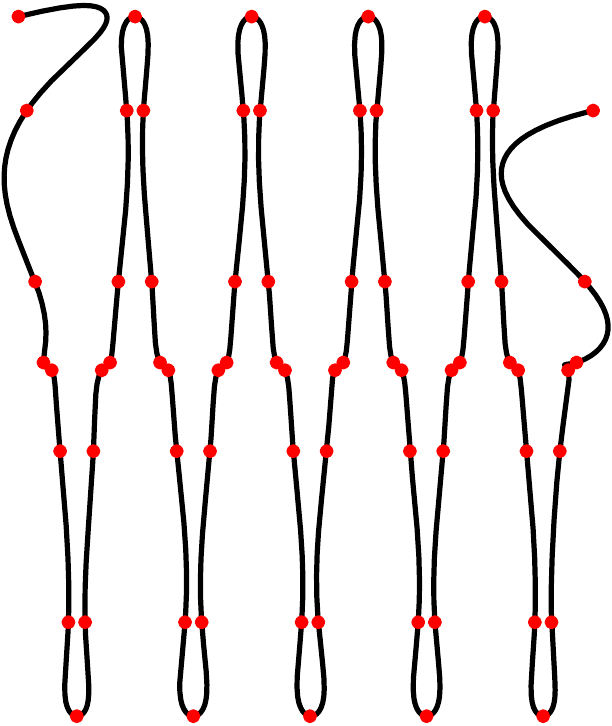}
    \caption{The curve after smoothing}
  \end{subfigure}
  \caption{The result of algorithm applied to \eqref{eq:openex2}. The red dots 
 mark the sample points.}
  \label{fig:open2}
\end{figure}

\begin{figure}[!h]
  \centering
  \begin{subfigure}[b]{0.04\linewidth}
   \captionsetup{width=3\linewidth}
    \includegraphics[width=\linewidth]{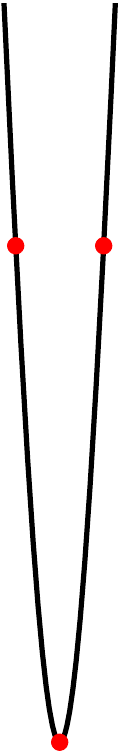}
    \caption{Before smoothing}
  \end{subfigure}\hspace{40mm} 
  \begin{subfigure}[b]{0.046\linewidth}
   \captionsetup{width=3\linewidth}
    \includegraphics[width=\linewidth]{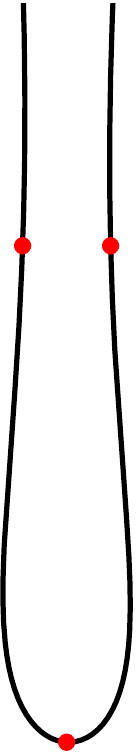}
    \caption{After smoothing}
  \end{subfigure}
  \caption{A detail of Figure~\ref{fig:open2}}
  \label{fig:open2t}
\end{figure}

\begin{figure}[!h]
  \centering
  \begin{subfigure}[b]{0.45\linewidth}
    \includegraphics[width=\linewidth]{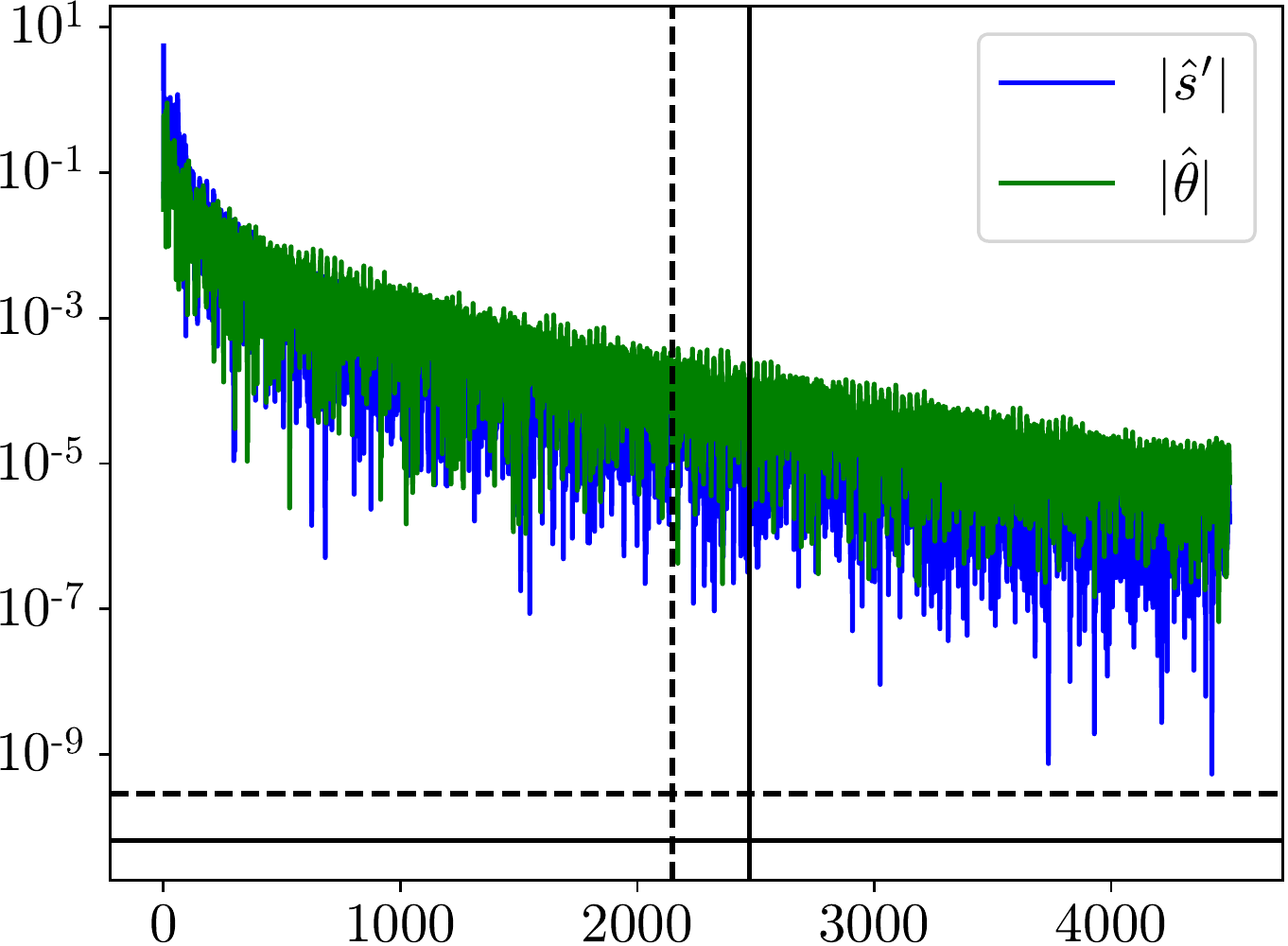}
    \caption{Before filtering}
  \end{subfigure}
  \begin{subfigure}[b]{0.45\linewidth}
    \includegraphics[width=\linewidth]{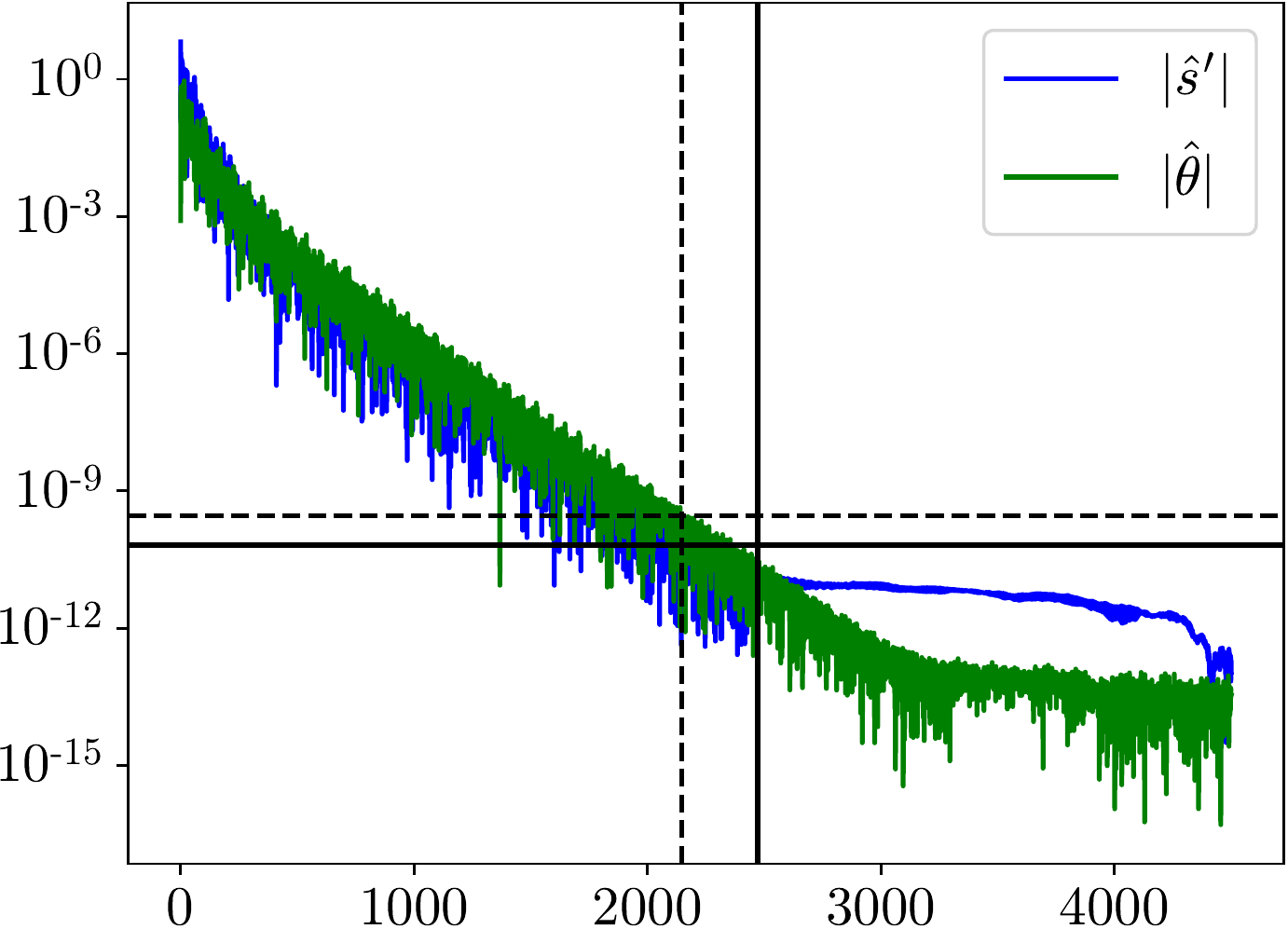}
    \caption{After filtering}
  \end{subfigure}
  \caption{Chebyshev coefficients of $s'(t)$ and $\theta(t)$ corresponding to 
  Figure~\ref{fig:open2}. The value of $\delta_{s'}$ is indicated by a horizontal solid line 
  and the value of $\delta_{\theta}$ is indicated by a horizontal dashed line. The 
  $2472$nd coefficients of $s'(t)$ decay to $\delta_{s'}$, indicated by a vertical
  solid line.
  The $2148$th coefficients of $\theta(t)$ decays to $\delta_{\theta}$,
  indicated by a vertical dashed line.} 
  \label{fig:open2c}
\end{figure}

\begin{figure}[!h]
  \centering
  \begin{subfigure}[b]{0.45\linewidth}
    \includegraphics[width=\linewidth]{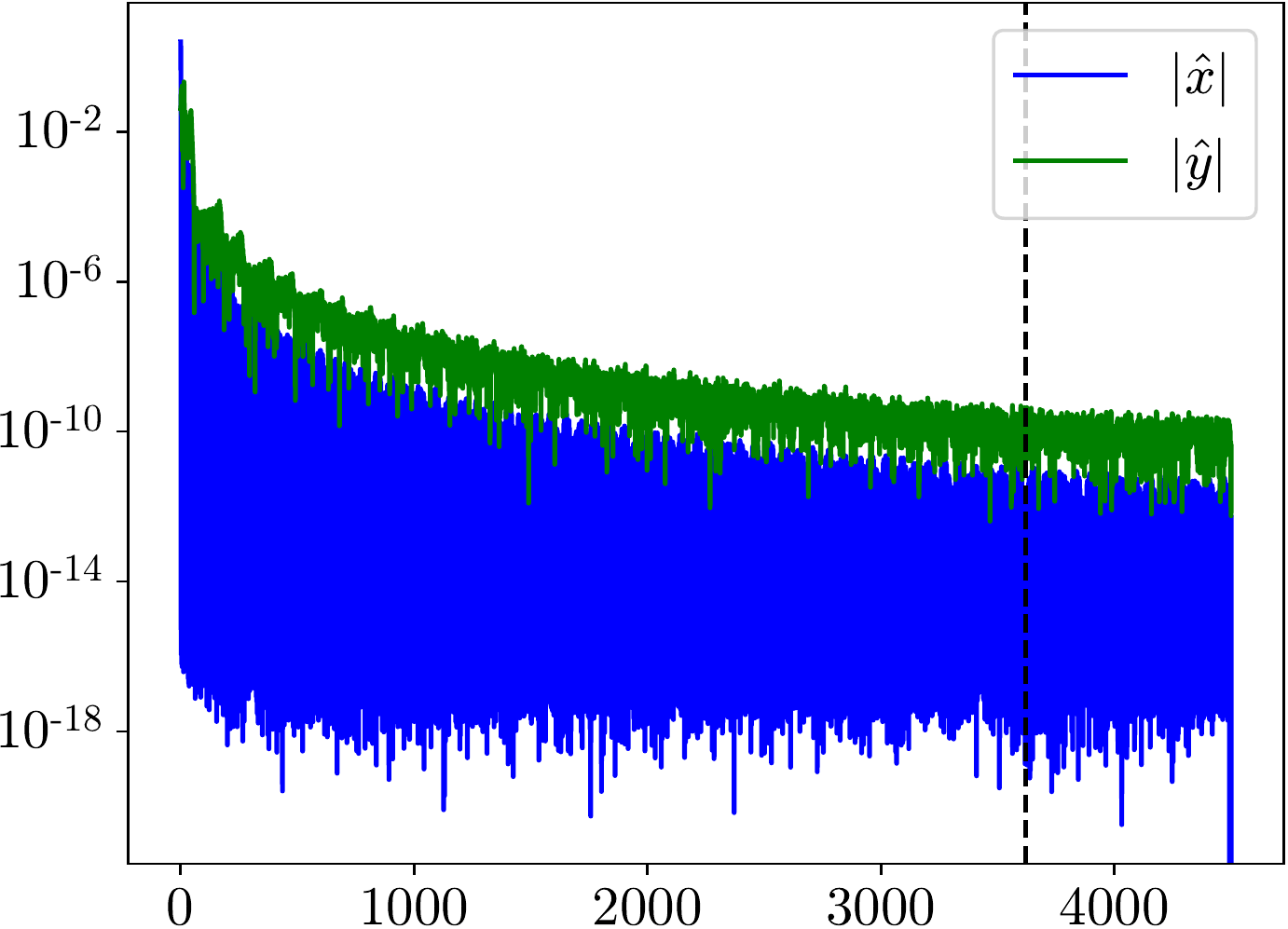}
    \caption{Chebyshev coefficients of the initial curve}
  \end{subfigure}
  \begin{subfigure}[b]{0.45\linewidth}
    \includegraphics[width=\linewidth]{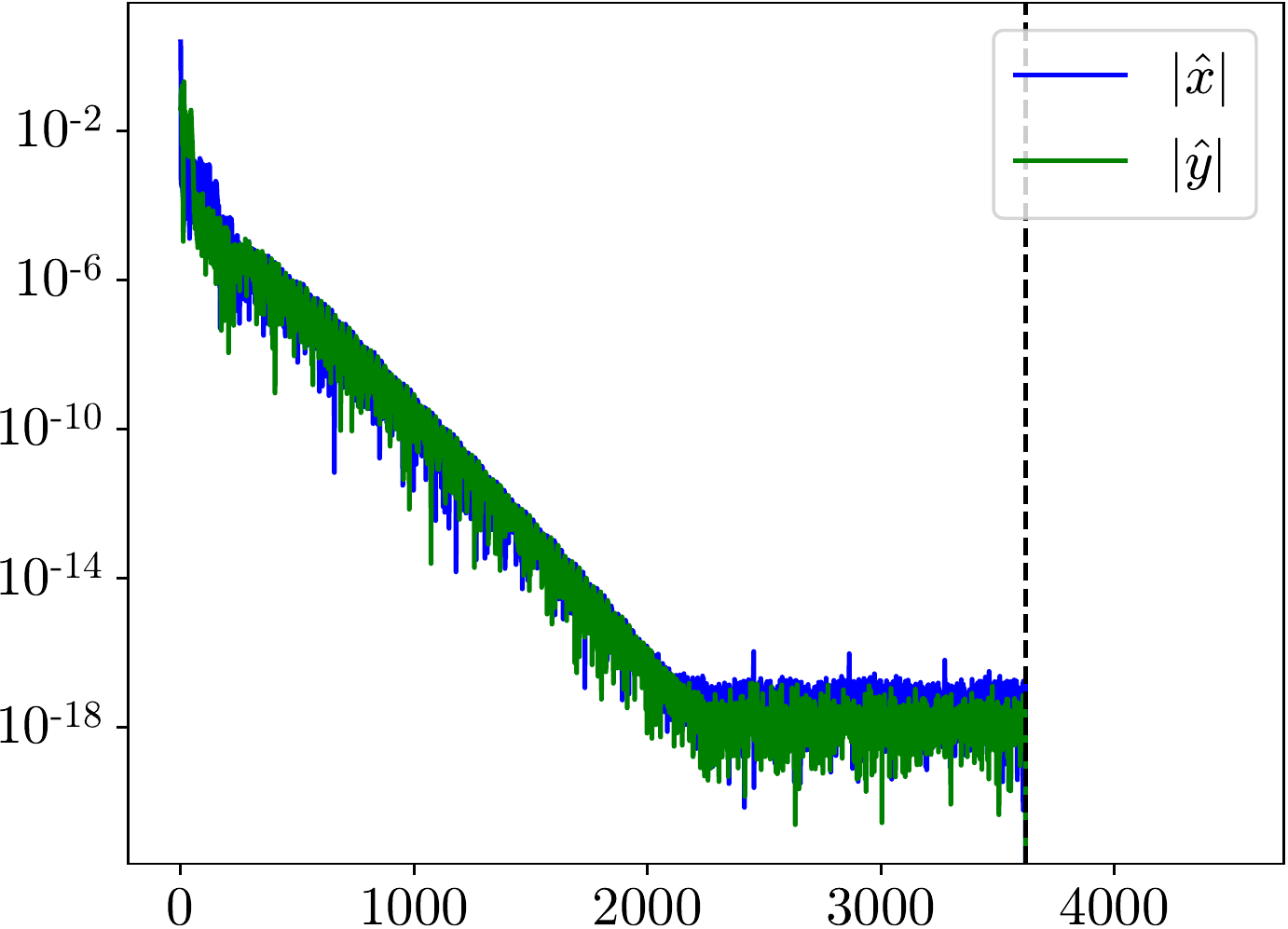}
    \caption{Chebyshev coefficients of the final curve}
  \end{subfigure}
  \caption{Chebyshev coefficients of $x(t)$ and $y(t)$ corresponding to
  Figure~\ref{fig:open2}.
  The value of $n_{\text{coefs}}$ is indicated by a vertical dashed line.} 
  \label{fig:open2c2}
\end{figure}

Another example depicted in Figure~\ref{fig:open2}(a)
is obtained by sampling from the curve
  \begin{align}
  \label{eq:openex2} 
  \gamma{(t)}=\{5t, 3 \cos{(10t \pi)^3}\}, \qquad t \in [0,1].
  \end{align}
The sample data points are scaled so that their height is $1$.
We run the algorithm by choosing  
$n=70$, $N=4500$, 
$x'_{\text{left}}=0.25$, $y'_{\text{left}}=0.25$,
$x'_{\text{right}}=0.25$,
$y'_{\text{right}}=0.25$,
$n_{\text{iters}}=70$, $h_{\text{filter}}=\frac{1}{45}$,
$\epsilon=10^{-16}$, $n_{\text{coefs}}=3620$,
$n_{\text{bands}}=6$. 
The curve before smoothing 
is observed to bend unnaturally when zooming in on some details, for example, those shown in Figure~\ref{fig:open2t}(a). 
Thus, a reasonably large number of Chebyshev
coefficients are required to represent $s'(t)$ and $\theta(t)$, as shown in Figure~\ref{fig:open2c}.
By looking at Figure~\ref{fig:open2}(b) and Figure~\ref{fig:open2t}(b),
the curve appears more like
a manually drawn smooth curve after $n_{\text{stop}}=60$ iterations.
The coefficients returned by the algorithm represent a curve passing through
the sample data points to within an error of $E_{\text{samp}}=0.16875 \e{-13}$. 
We display the magnitudes of the Chebyshev coefficients of both the initial 
and final curve 
in Figure~\ref{fig:open2c2}.
  
\begin{figure}[!h]
  \centering
  \begin{subfigure}[b]{0.18\linewidth}
    \includegraphics[width=\linewidth]{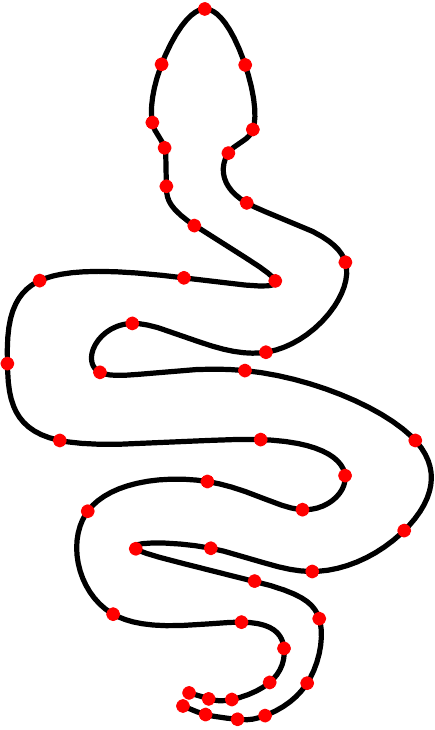}
    \caption{The curve before smoothing}
  \end{subfigure}\hspace{10mm}
  \begin{subfigure}[b]{0.18\linewidth}
    \includegraphics[width=\linewidth]{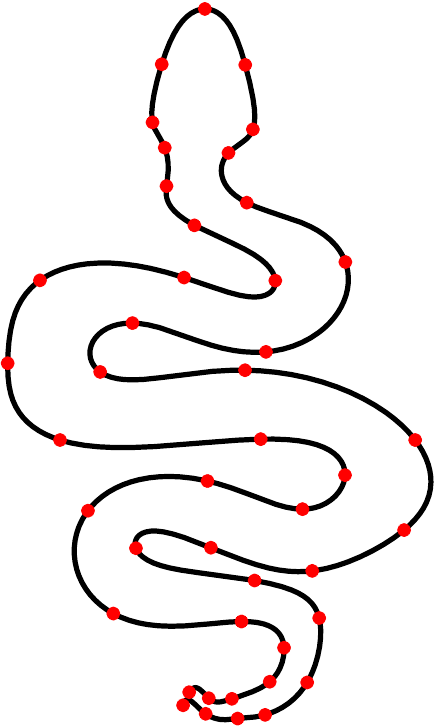}
    \caption{The curve after smoothing}
  \end{subfigure}
  \caption{A hand-drawn depiction of a snake shape. The red dots 
 mark the sample points.}
  \label{fig:open3}
\end{figure}

\begin{figure}[!h]
  \centering
  \begin{subfigure}[b]{0.45\linewidth}
    \includegraphics[width=\linewidth]{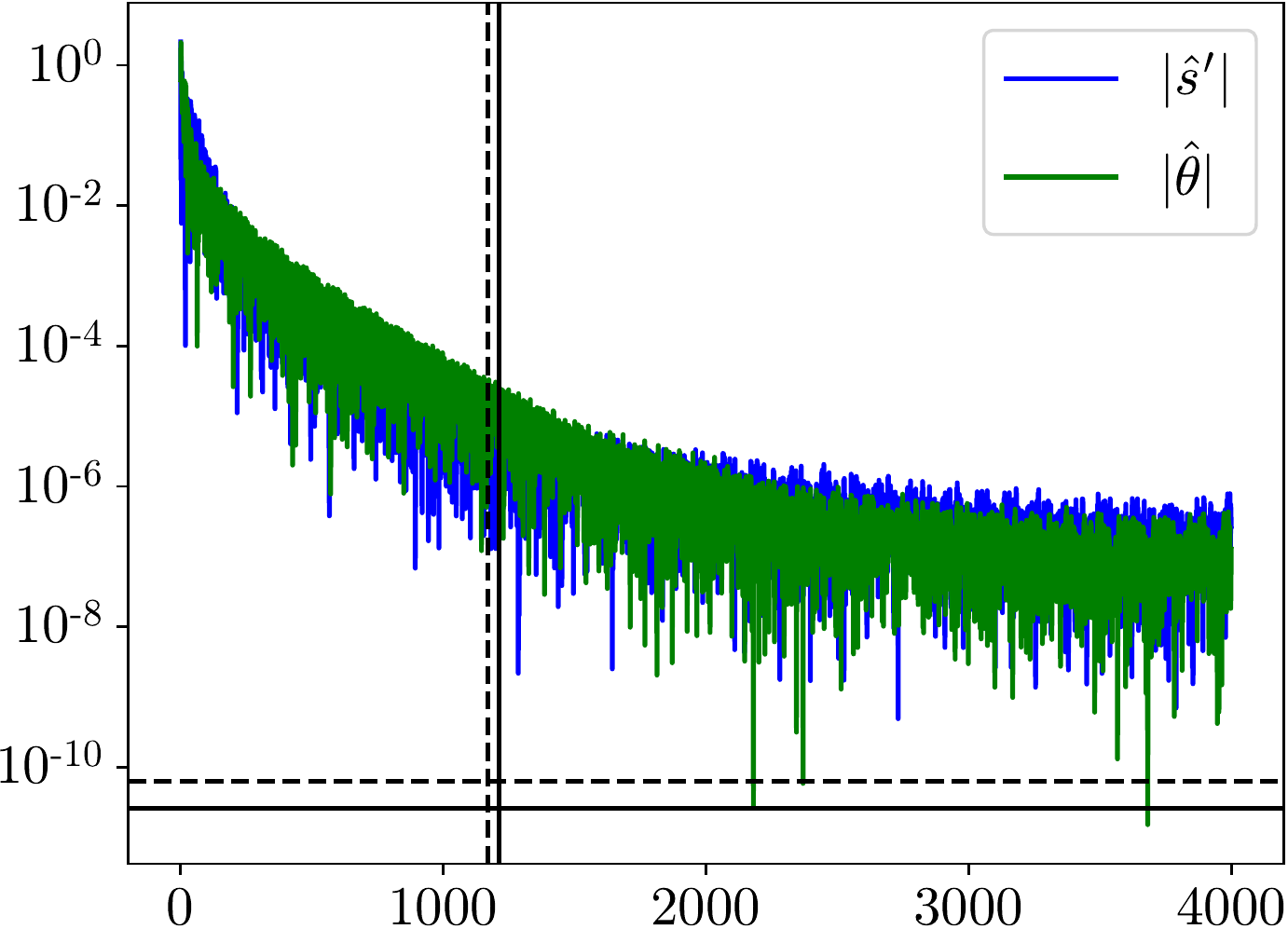}
    \caption{Before filtering}
  \end{subfigure}
  \begin{subfigure}[b]{0.45\linewidth}
    \includegraphics[width=\linewidth]{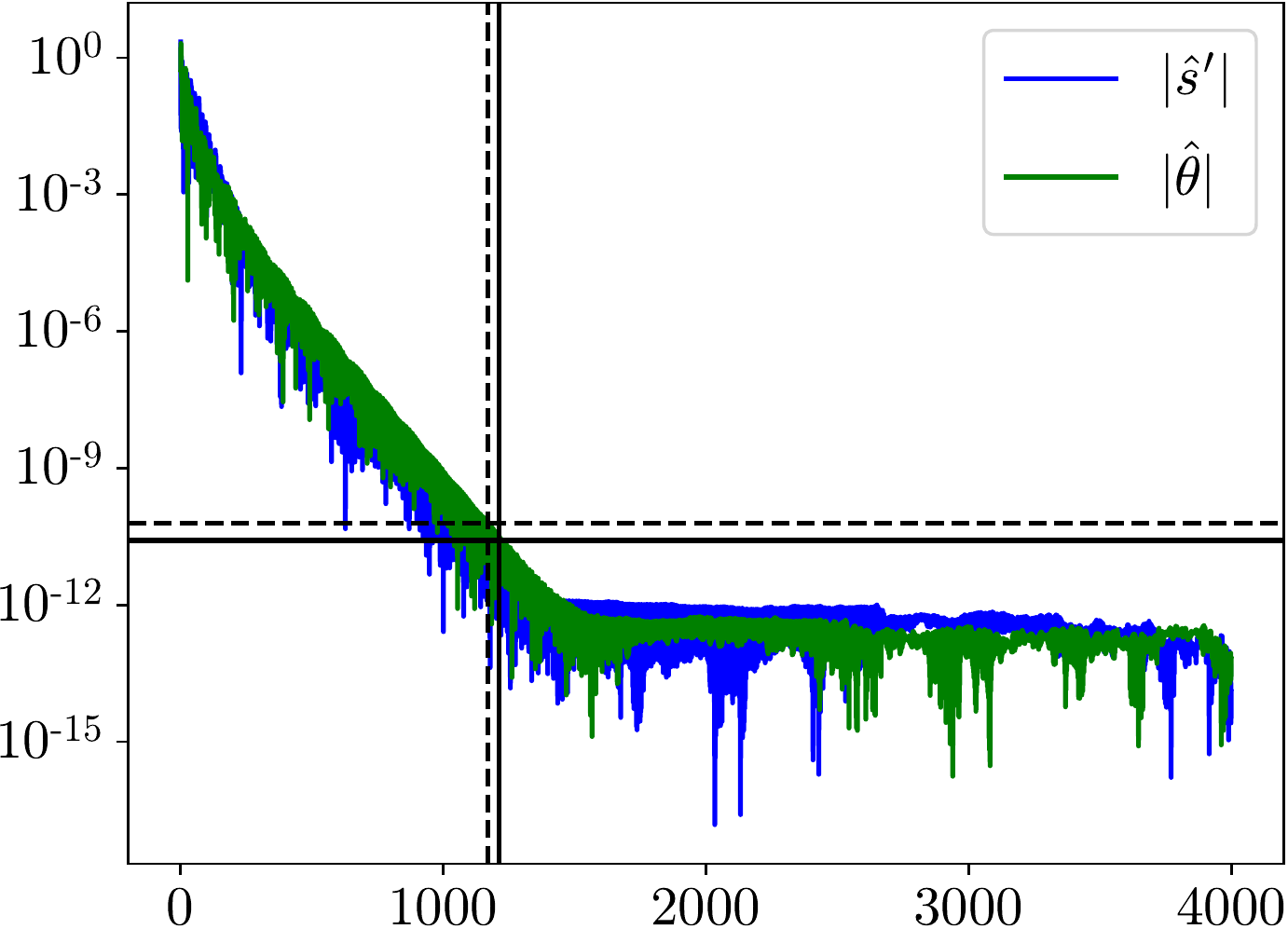}
    \caption{After filtering}
  \end{subfigure}
  \caption{Chebyshev coefficients of $s'(t)$ and $\theta(t)$ corresponding to 
  Figure~\ref{fig:open3}. The value of $\delta_{s'}$ is indicated by a horizontal solid line 
  and the value of $\delta_{\theta}$ is indicated by a horizontal dashed line.
  The $1214$th coefficients of $s'(t)$ decays to $\delta_{s'}$, indicated by a 
  vertical solid line.
  The $1171$st coefficients of $\theta(t)$ decays to $\delta_{\theta}$,
  indicated by a vertical dashed line.}
  \label{fig:open3c}
\end{figure}

\begin{figure}[!h]
  \centering
  \begin{subfigure}[b]{0.45\linewidth}
    \includegraphics[width=\linewidth]{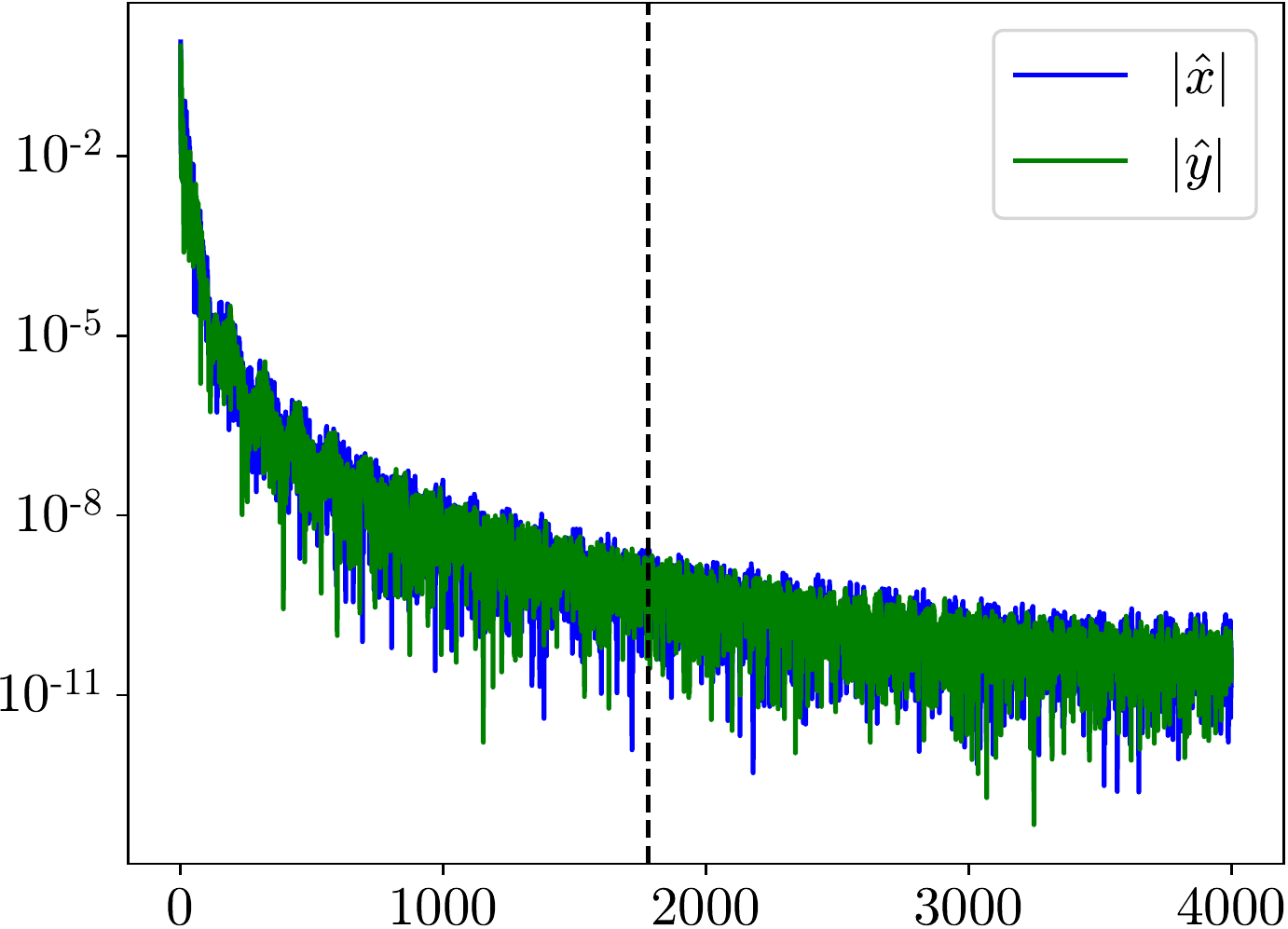}
    \caption{Chebyshev coefficients of the initial curve}
  \end{subfigure}
  \begin{subfigure}[b]{0.45\linewidth}
    \includegraphics[width=\linewidth]{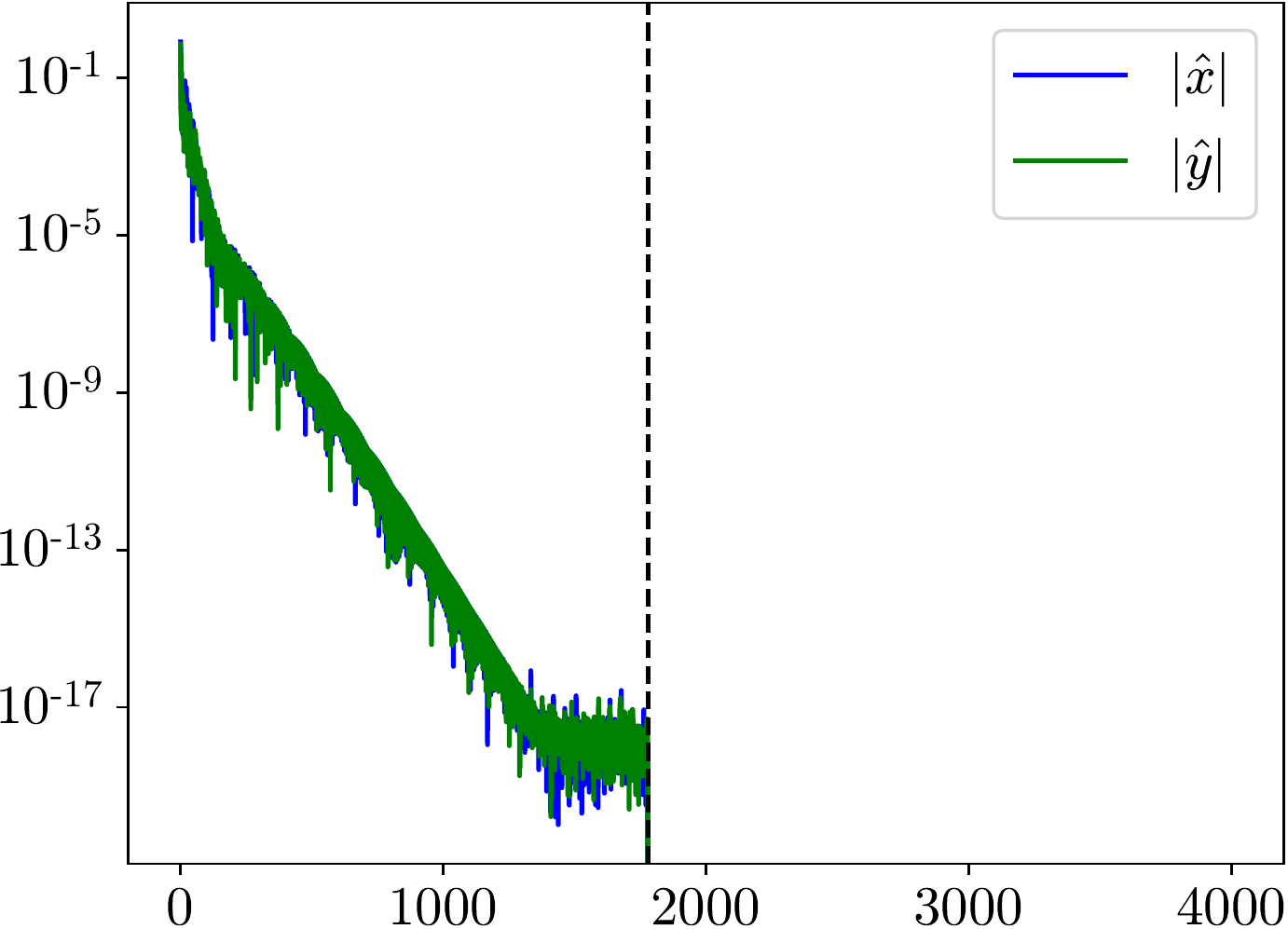}
    \caption{Chebyshev coefficients of the final curve}
  \end{subfigure}
  \caption{Chebyshev coefficients of $x(t)$ and $y(t)$ corresponding to
  Figure~\ref{fig:open3}.
  The value of $n_{\text{coefs}}$ is indicated by a vertical dashed line.} 
  \label{fig:open3c2}
\end{figure}

Figure~\ref{fig:open3}(a) shows a roughly sketched shape resembling 
a snake.
We scale the sample data points so that their height is $1$,
and run the algorithm by choosing $N=4000$, $n=44$,
$x'_{\text{left}}=0.05$, $y'_{\text{left}}=-0.02$,
$x'_{\text{right}}=-0.06$, $y'_{\text{right}}=0.02$,
$n_{\text{iters}}=80$, $h_{\text{filter}}=\frac{1}{50}$,
$\epsilon=10^{-16}$, $n_{\text{coefs}}=1780$,
$n_{\text{bands}}=6$.
The algorithm terminates at the $n_{\text{stop}}=71$st iteration, and the resulting curve 
passes through
the sample data points to within an error of $E_{\text{samp}}=0.35056 \e{-14}$.
We present the magnitudes of the coefficients of $s'(t)$ 
and $\theta(t)$ before and after filtering in Figure~\ref{fig:open3c}, 
and the magnitudes of the coefficients of $x(t)$ 
and $y(t)$ of both the initial and final curve in Figure~\ref{fig:open3c2}. 

\begin{figure}[!h]
  \centering
  \begin{subfigure}[b]{0.35\linewidth}
    \includegraphics[width=\linewidth]{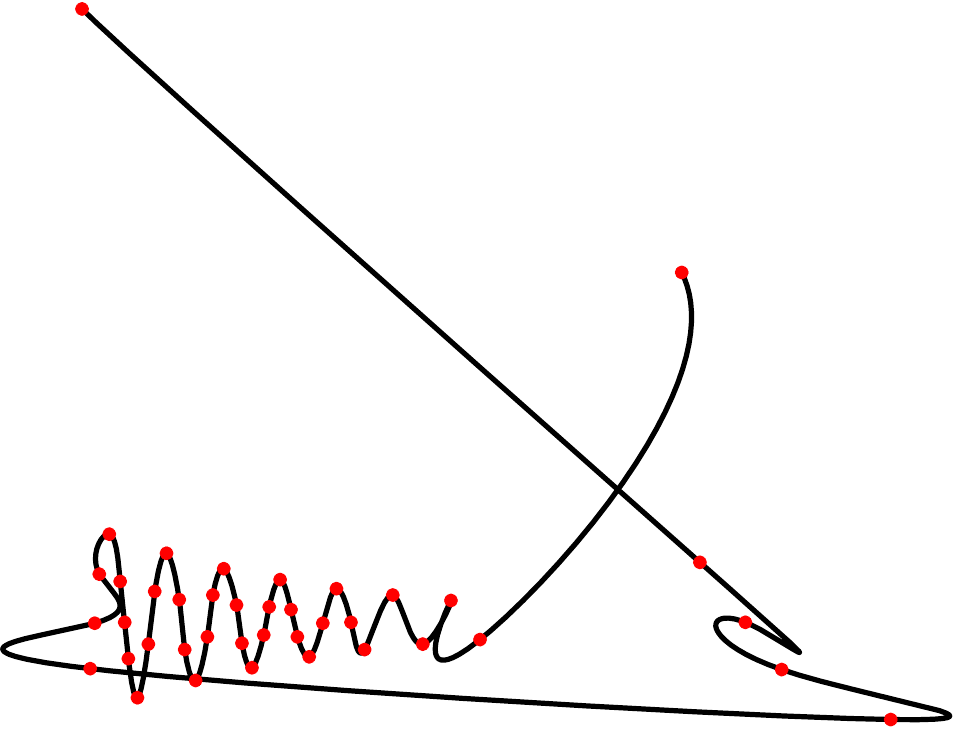}
    \caption{The curve before smoothing}
  \end{subfigure}\hspace{10mm}
  \begin{subfigure}[b]{0.35\linewidth}
    \includegraphics[width=\linewidth]{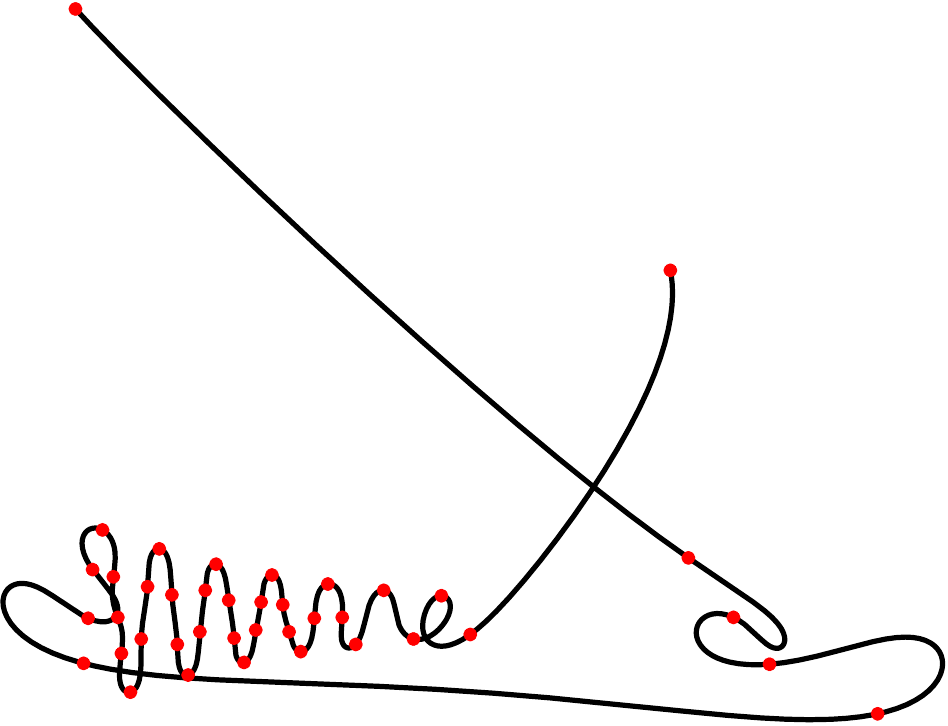}
    \caption{The curve after smoothing}
  \end{subfigure}
  \caption{A shape of oscillations that exhibit damping. The red dots 
 mark the sample points.}
  \label{fig:open4}
\end{figure}

\begin{figure}[!h]
  \centering
  \begin{subfigure}[b]{0.45\linewidth}
    \includegraphics[width=\linewidth]{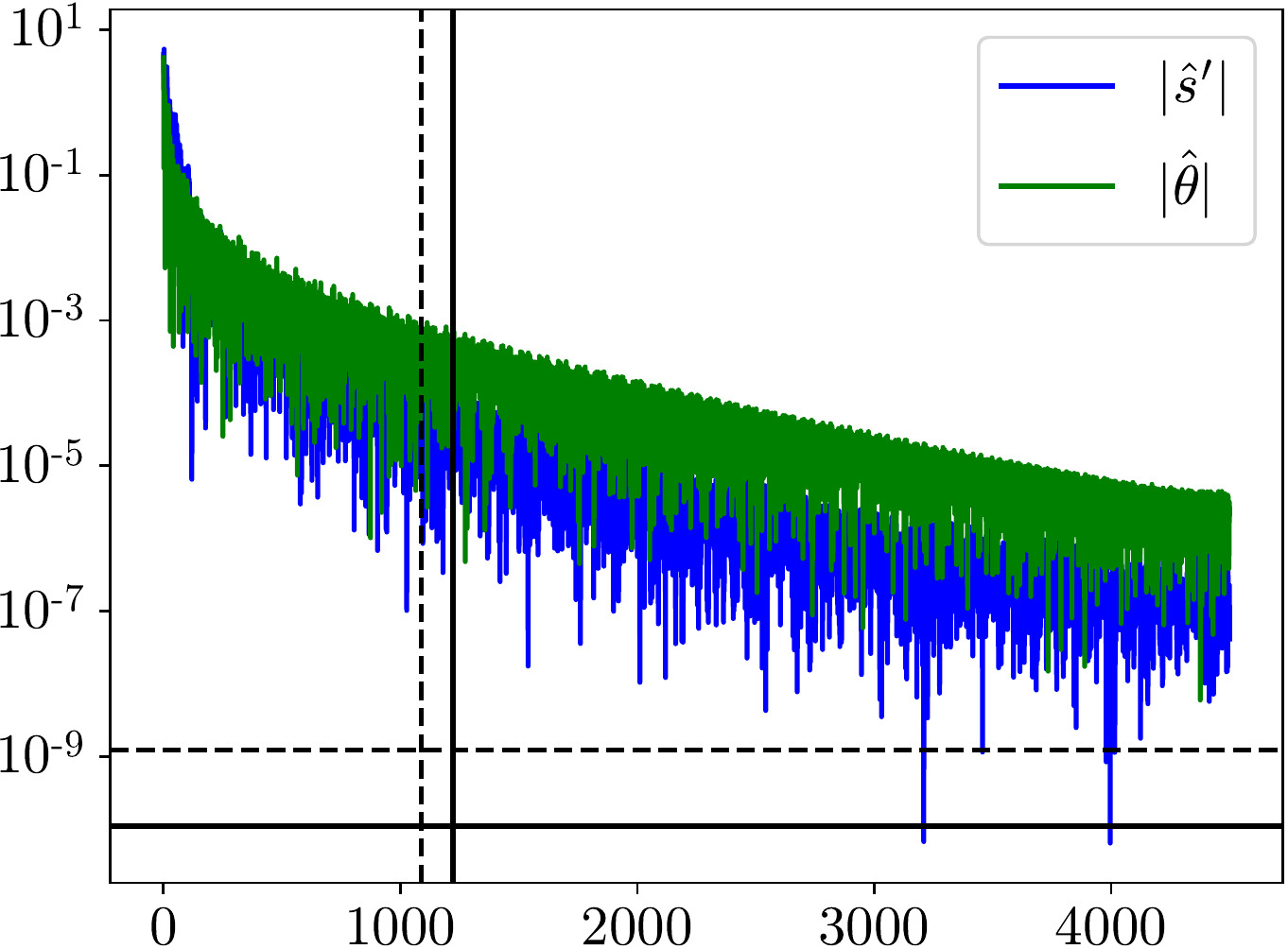}
    \caption{Before filtering}
  \end{subfigure}
  \begin{subfigure}[b]{0.45\linewidth}
    \includegraphics[width=\linewidth]{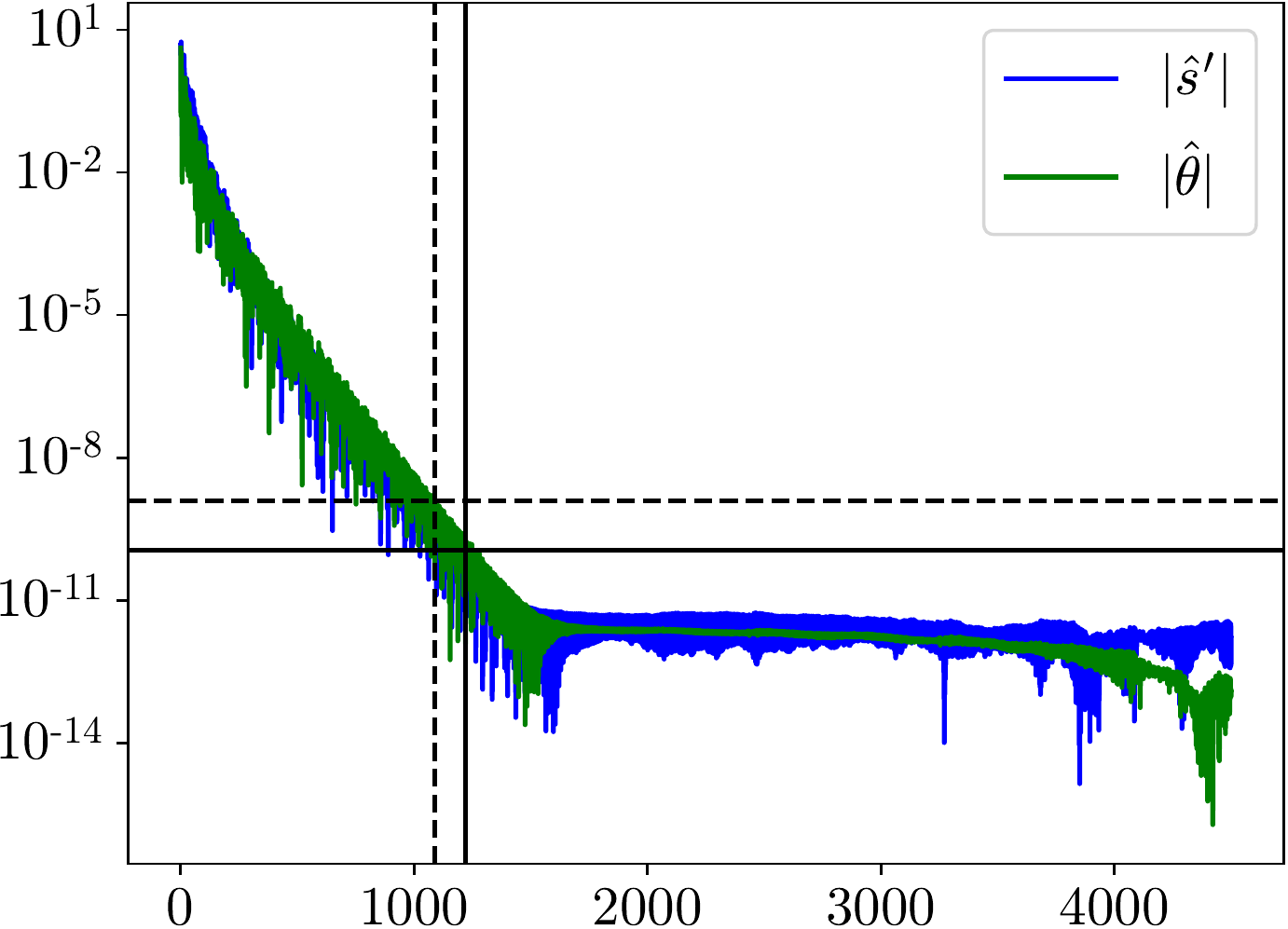}
    \caption{After filtering}
  \end{subfigure}
  \caption{Chebyshev coefficients of $s'(t)$ and $\theta(t)$ corresponding to 
  Figure~\ref{fig:open4}. The value of $\delta_{s'}$ is indicated by a horizontal solid line 
  and the value of $\delta_{\theta}$ is indicated by a horizontal dashed line.
  The $1222$nd coefficients of $s'(t)$ decays to $\delta_{s'}$, indicated by a 
  vertical solid line.
  The $1088$th coefficients of $\theta(t)$ decays to $\delta_{\theta}$,
  indicated by a vertical dashed line.}
  \label{fig:open4c}
\end{figure}

\begin{figure}[!h]
  \centering
  \begin{subfigure}[b]{0.45\linewidth}
    \includegraphics[width=\linewidth]{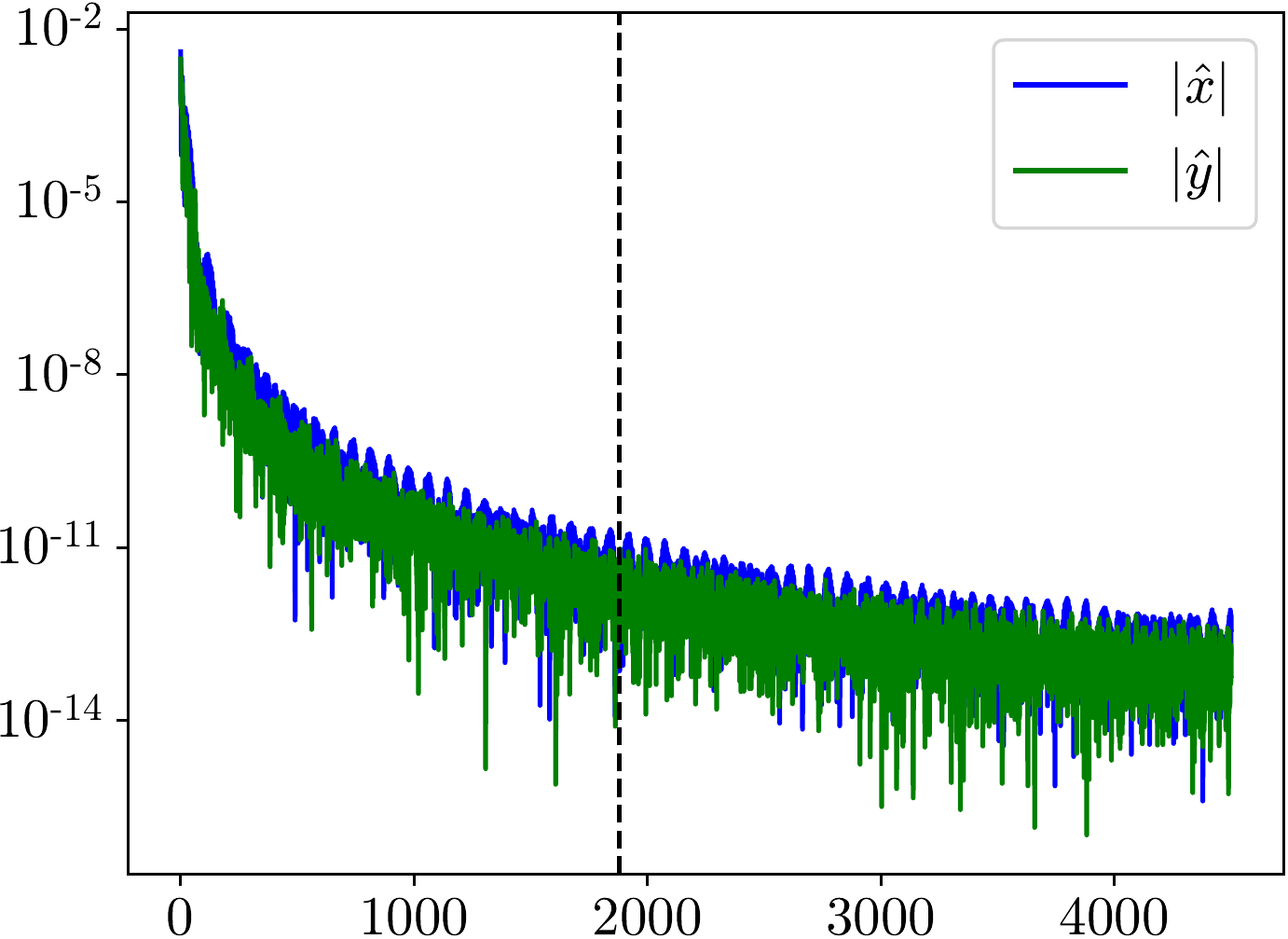}
    \caption{Chebyshev coefficients of the initial curve}
  \end{subfigure}
  \begin{subfigure}[b]{0.45\linewidth}
    \includegraphics[width=\linewidth]{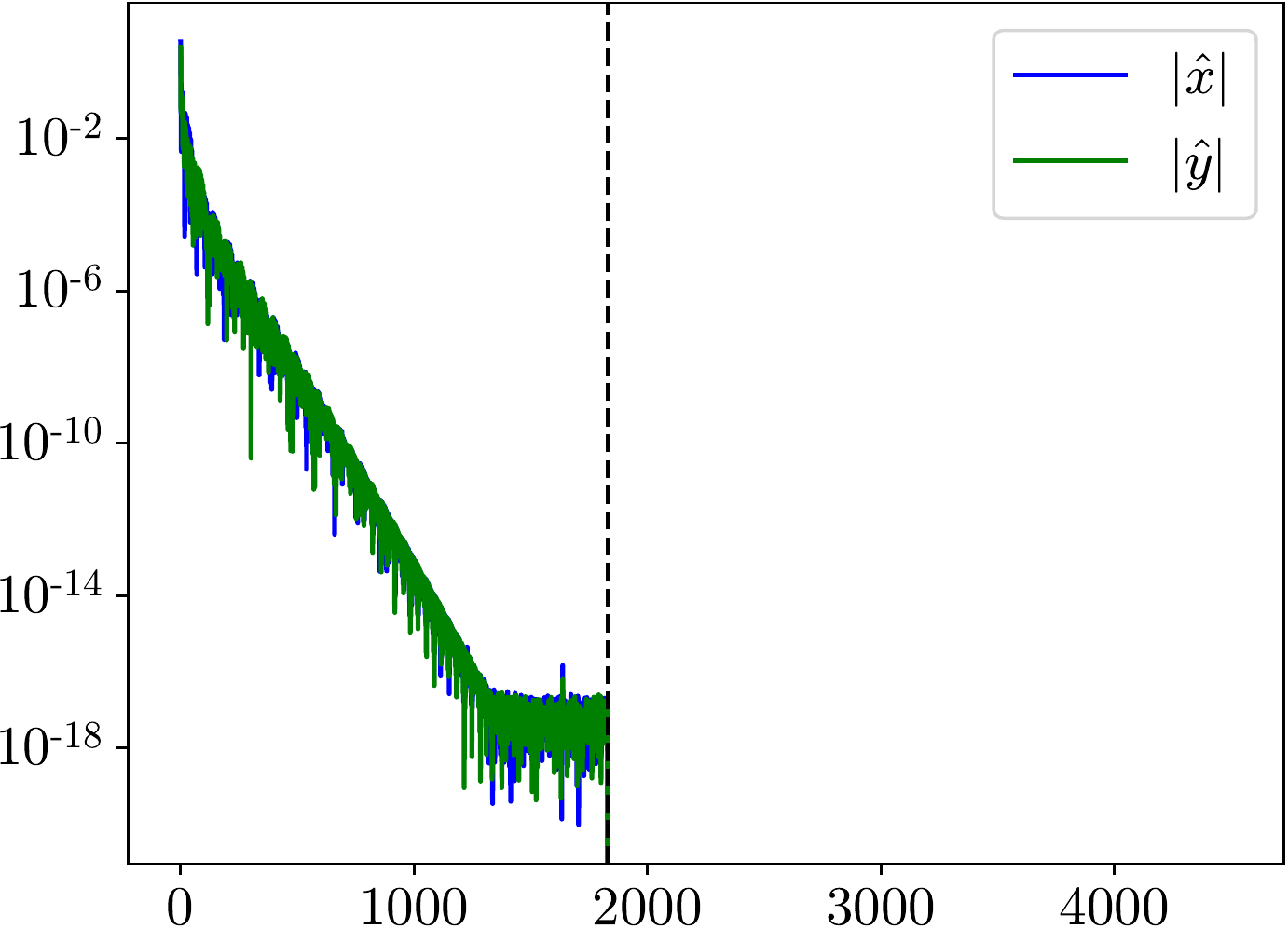}
    \caption{Chebyshev coefficients of the final curve}
  \end{subfigure}
  \caption{Chebyshev coefficients of $x(t)$ and $y(t)$ corresponding to
  Figure~\ref{fig:open4}.
  The value of $n_{\text{coefs}}$ is indicated by a vertical dashed line.} 
  \label{fig:open4c2}
\end{figure}

\begin{figure}[!h]
  \centering
  \begin{subfigure}[b]{0.35\linewidth}
    \includegraphics[width=\linewidth]{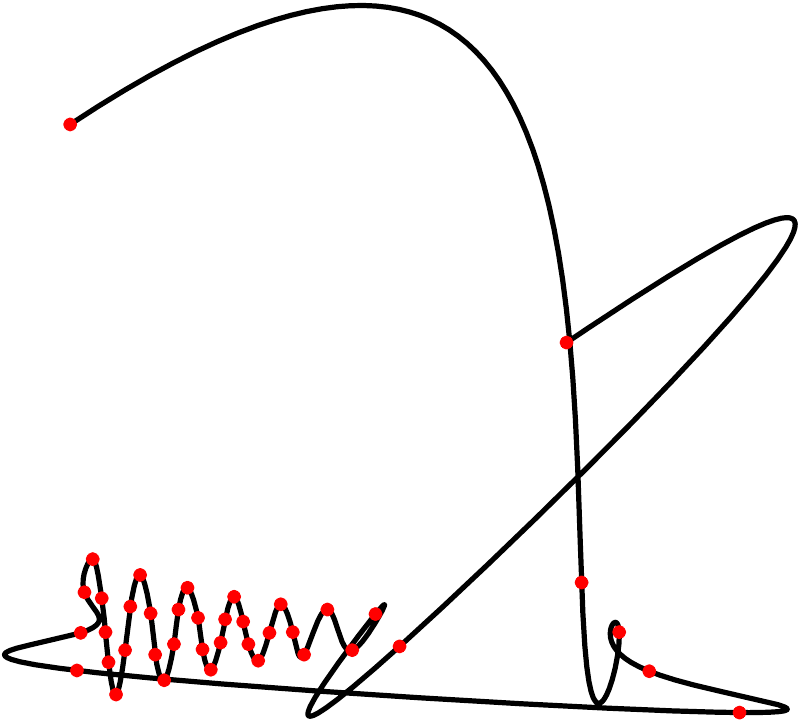}
    \caption{The interpolating curve. The red dots 
       mark the sample points.}
  \end{subfigure}\hspace{10mm}
  \begin{subfigure}[b]{0.45\linewidth}
    \includegraphics[width=\linewidth]{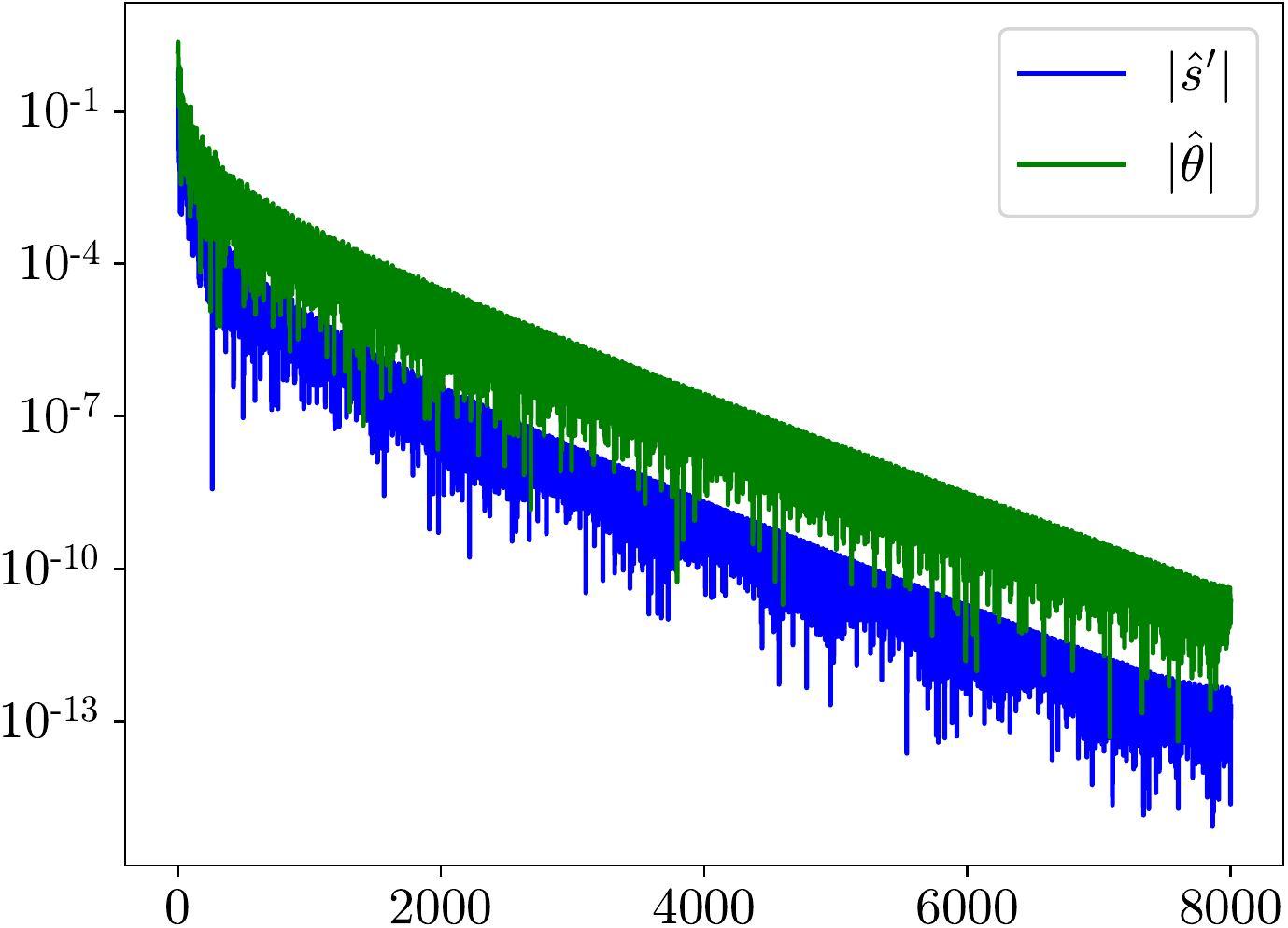}
    \caption{Chebyshev coefficients of $s'(t)$ and $\theta(t)$ corresponding to 
  Figure~\ref{fig:open4_2}(a).}
  \end{subfigure}
    \caption{The result of algorithm in \cite{zhang} applied to the 
    same data points in Figure~\ref{fig:open4}.}
  \label{fig:open4_2}
\end{figure}

We illustrate some damping oscillations, as displayed in Figure~\ref{fig:open4}(a).
The sample data points are scaled so that their width is $1$.
The initial curve has some sharp corners, and is distorted unnaturally. 
We set $N=4500$, $n=40$,
$x'_{\text{left}}=0.20$, $y'_{\text{left}}=-0.20$,
$x'_{\text{right}}=-0.20$,
$y'_{\text{right}}=0.40$,
$n_{\text{iters}}=70$, $h_{\text{filter}}=\frac{1}{40}$,
$\epsilon=10^{-16}$, $n_{\text{coefs}}=1830$,
$n_{\text{bands}}=8$. 
After $n_{\text{stop}}=69$ iterations, the algorithm terminates and
returns a curve 
passing through
the sample data points to within an error of $E_{\text{samp}}=0.22649 \e{-13}$.
The resulting curve in Figure~\ref{fig:open4}(b)
resembles a curve drawn by hand, with a completely smooth 
shape that naturally bends to pass through all the sample
data points to exhibit those damping
oscillations.
We present the magnitudes of the coefficients of $s'(t)$ 
and $\theta(t)$ before and after filtering in Figure~\ref{fig:open4c}, 
and the magnitudes of the coefficients of $x(t)$ 
and $y(t)$ of both the initial and final curve in Figure~\ref{fig:open4c2}. 
We apply the algorithm in \cite{zhang} to the same data points, by setting 
$a=0.2$, which produces the smoothest shape of the curve,
as displayed in Figure~\ref{fig:open4_2}(a).
Although 
the curve
in Figure~\ref{fig:open4_2}(a)
requires fewer coefficients to represent $x(t)$ and $y(t)$
compared to the curve in Figure \ref{fig:open4}(b),
it requires a much larger 
number of coefficients for $s'(t)$ and $\theta(t)$, as shown in 
\ref{fig:open4_2}(b).
This results in a curve with a high level of curvature.
With our algorithm,
any high curvature areas are effectively smoothed,
yielding a visually smoother curve
and requiring much fewer coefficients to represent $s'(t)$ and $\theta(t)$.

The runtimes per iteration for the open curve case are displayed in Table ~\ref{tab:table1}. Since we use the library \cite{fftw} for the 
implementation of the FFT, and the speed of the FFT routines in the library
depends in a complicated way on the input size,
we observed that the runtimes in Table~\ref{tab:table1} are not strictly proportional to the number of 
discretization points, $N$.






\begin{table}[h]
  \centering
\begin{tabular}{cccccc}
Case & $N=1025$ & $N=2049$ & $N=4097$ & $N=8193$  \\
\hline
\T Figure~\ref{fig:open1} &$0.14308\e{-02}$& $0.20470\e{-02}$ & $0.29810\e{-02}$& $0.51040\e{-02}$ \\
\end{tabular}

\caption{
  \label{tab:table1}
Average runtime per iteration, for an open curve, calculated by determining 
the total runtime for $100$ iterations  
and dividing by the number of iterations.
}

\end{table}

\subsection{Closed Curve Examples}
\begin{figure}[!h]
  \centering
  \begin{subfigure}[b]{0.28\linewidth}
    \includegraphics[width=\linewidth]{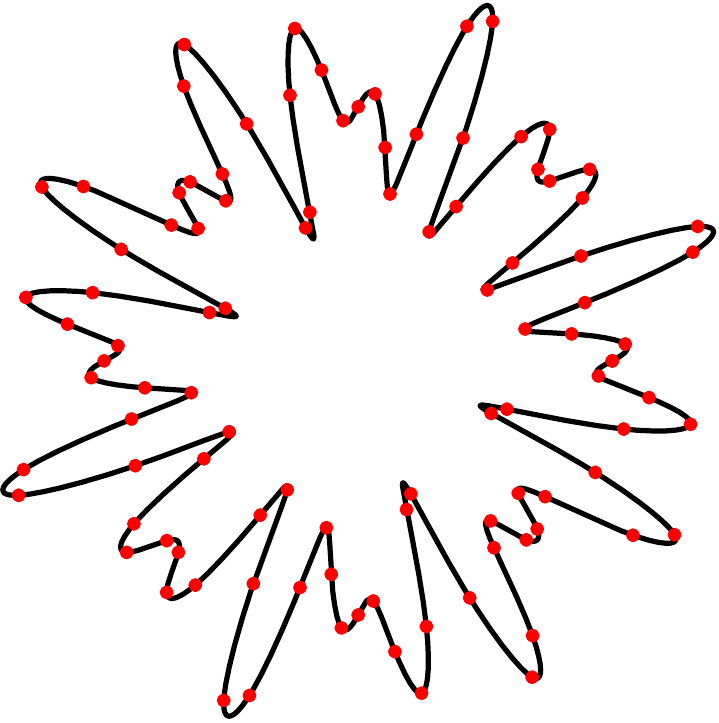}
    \caption{The curve before smoothing}
  \end{subfigure}\hspace{8mm}
  \begin{subfigure}[b]{0.28\linewidth}
    \includegraphics[width=\linewidth]{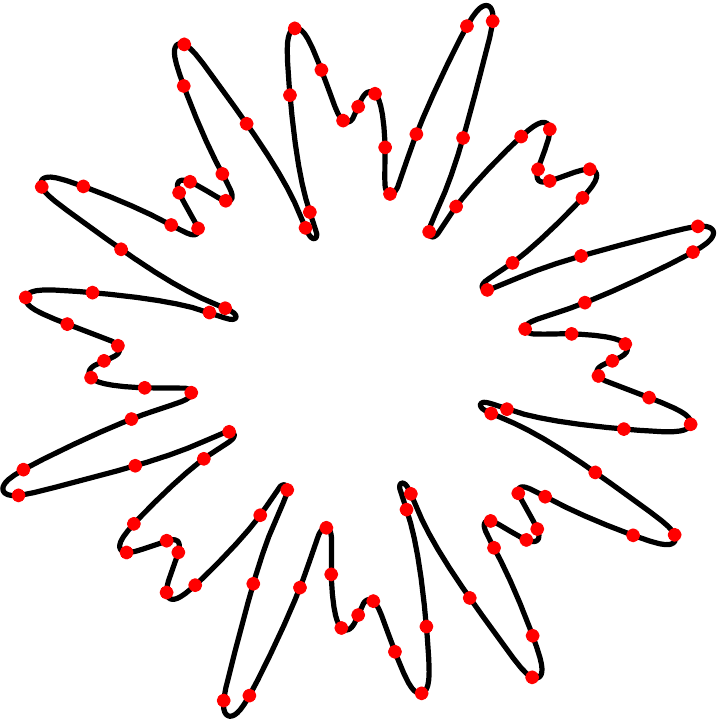}
    \caption{The curve after smoothing}
  \end{subfigure}
  \caption{The result of algorithm applied to \eqref{eq:closedex1} with $\alpha = 2$. The red dots 
 mark the sample points.}
  \label{fig:closed1}
\end{figure}

\begin{figure}[!h]
  \centering
  \begin{subfigure}[b]{0.45\linewidth}
    \includegraphics[width=\linewidth]{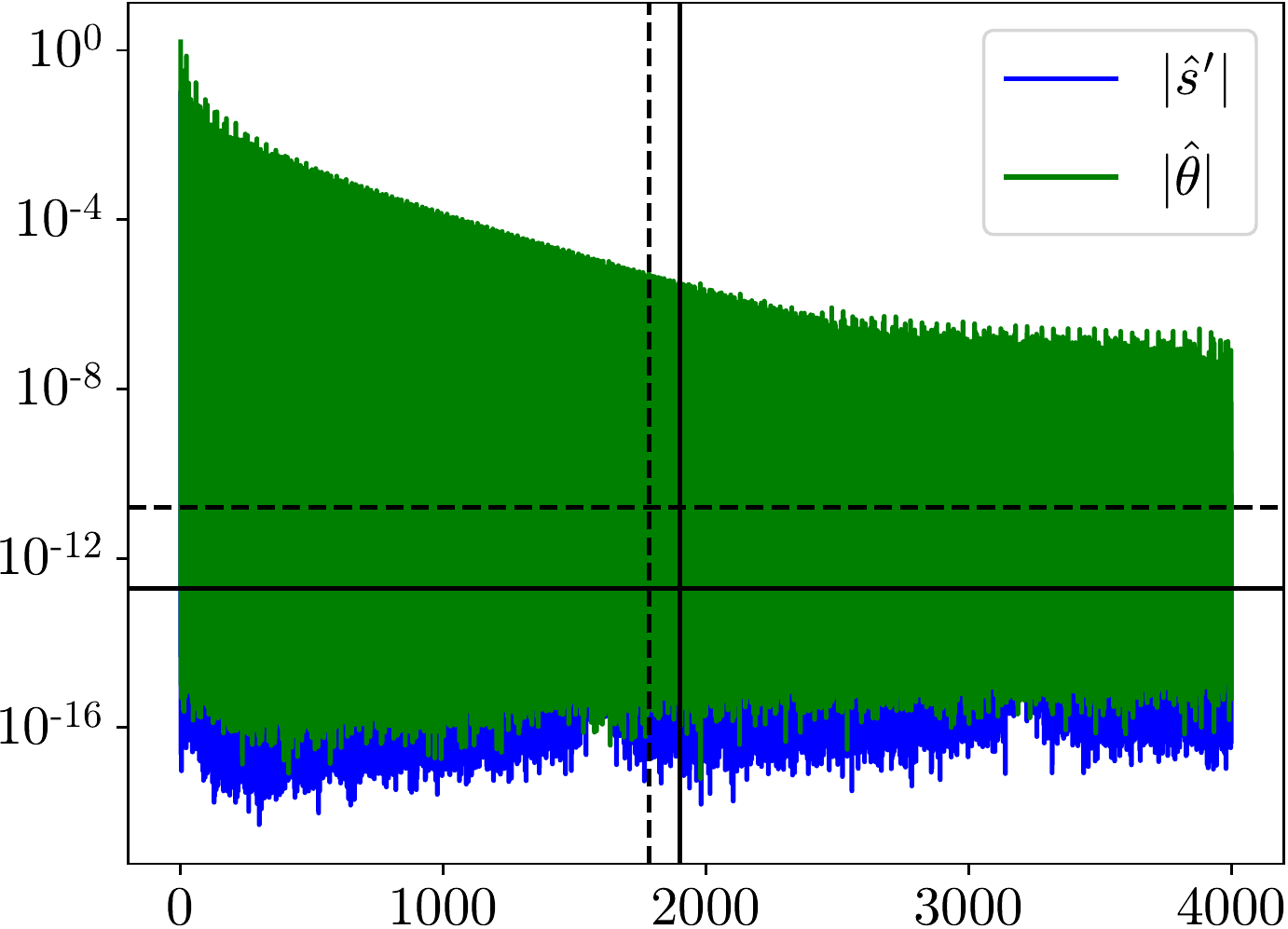}
    \caption{Before filtering}
  \end{subfigure}
  \begin{subfigure}[b]{0.45\linewidth}
    \includegraphics[width=\linewidth]{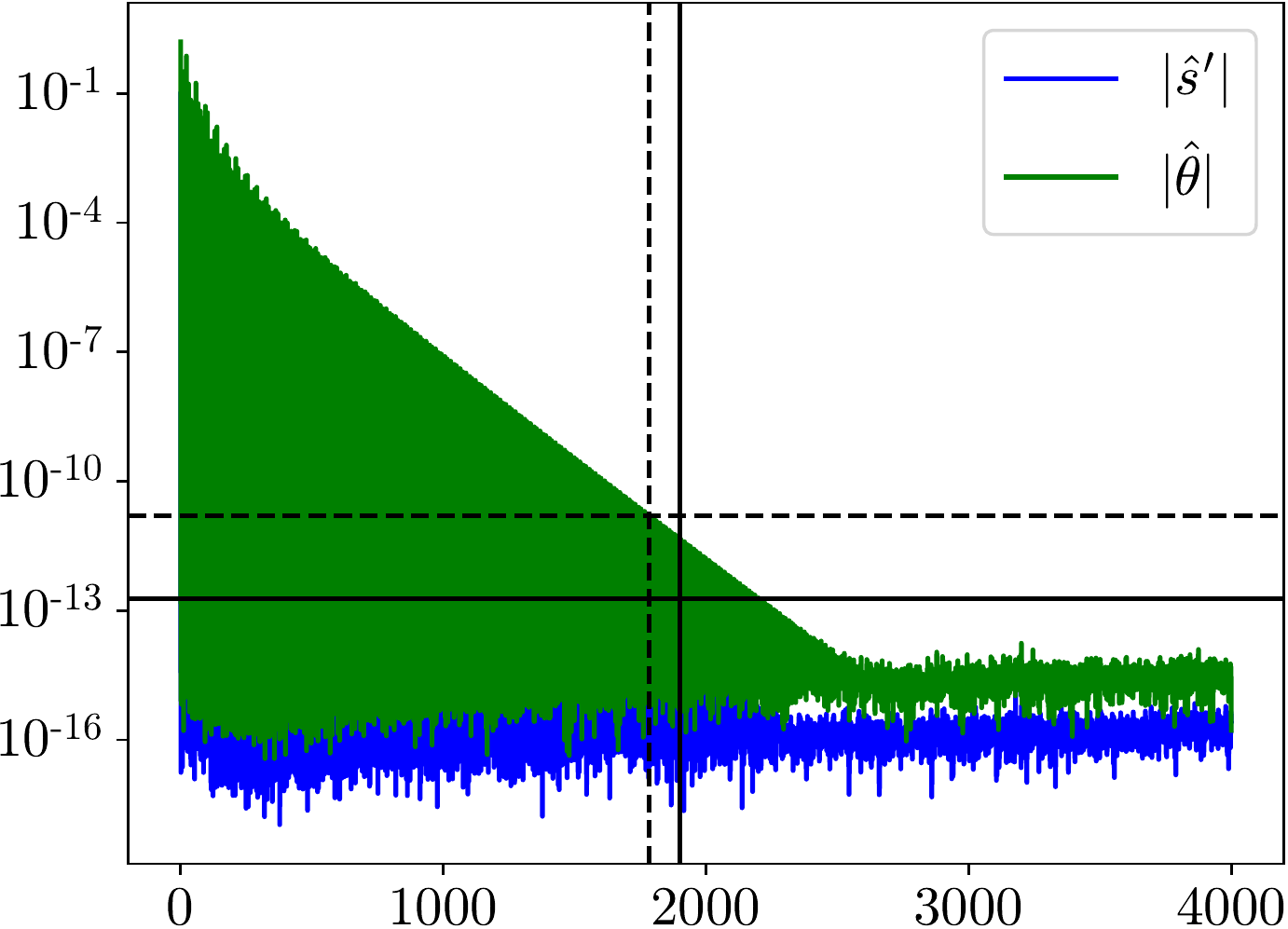}
    \caption{After filtering}
  \end{subfigure}
  \caption{Fourier coefficients of $s'(t)$ and $\theta(t)$ corresponding to 
  Figure~\ref{fig:closed1}.
  The value of $\delta_{s'}$ is indicated by a horizontal solid line 
  and the value of $\delta_{\theta}$ is indicated by a horizontal dashed line.
  The $1901$st coefficients of $s'(t)$ decays to $\delta_{s'}$, indicated by
  a vertical solid line.
  The $1785$th coefficients of $\theta(t)$ decays to $\delta_{\theta}$,
  indicated by a vertical dashed line.}
  \label{fig:closed1cs}
\end{figure}

\begin{figure}[!h]
  \centering
  \begin{subfigure}[b]{0.45\linewidth}
    \includegraphics[width=\linewidth]{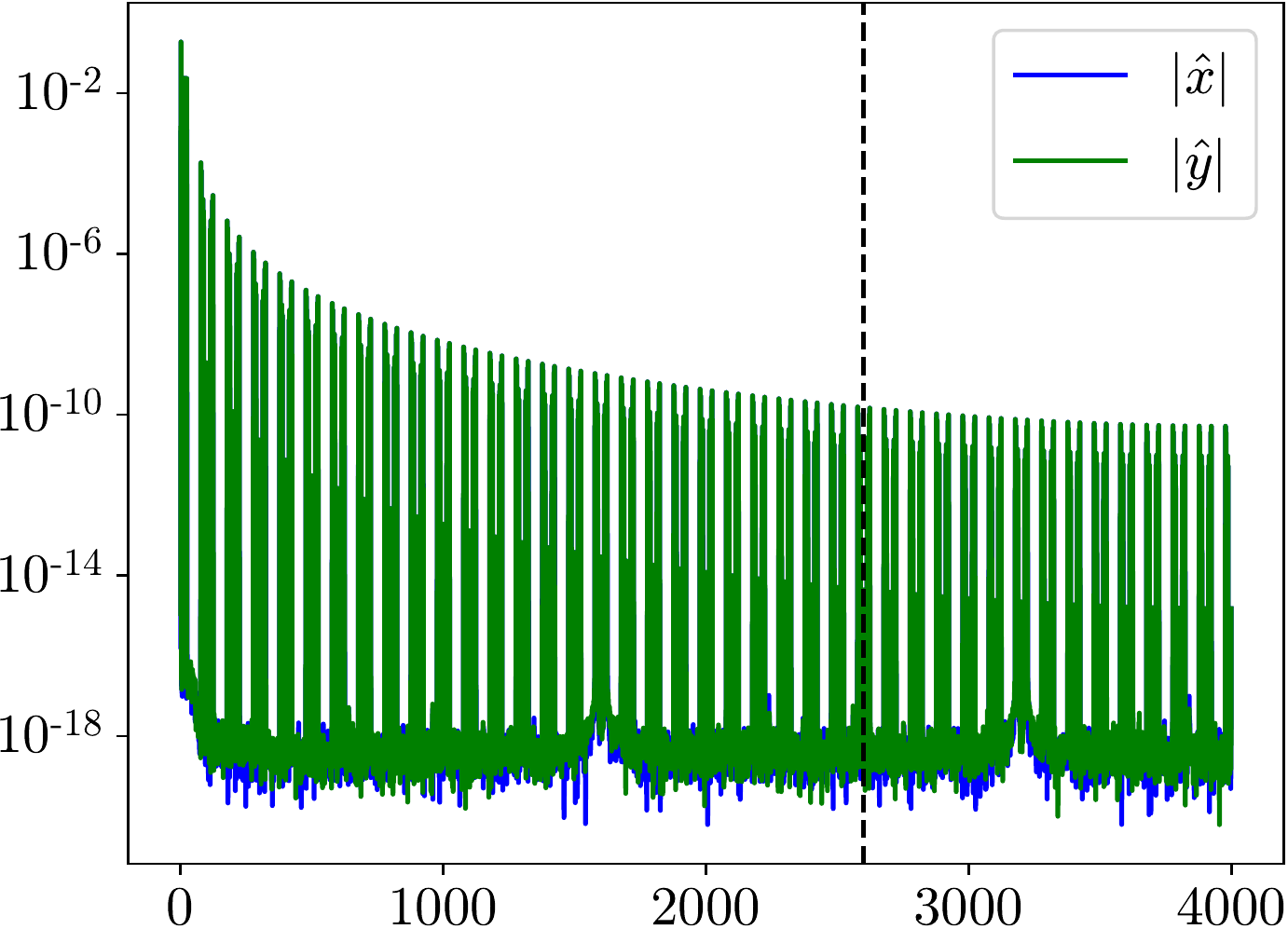}
    \caption{Fourier coefficients of the initial curve}
  \end{subfigure}
  \begin{subfigure}[b]{0.45\linewidth}
    \includegraphics[width=\linewidth]{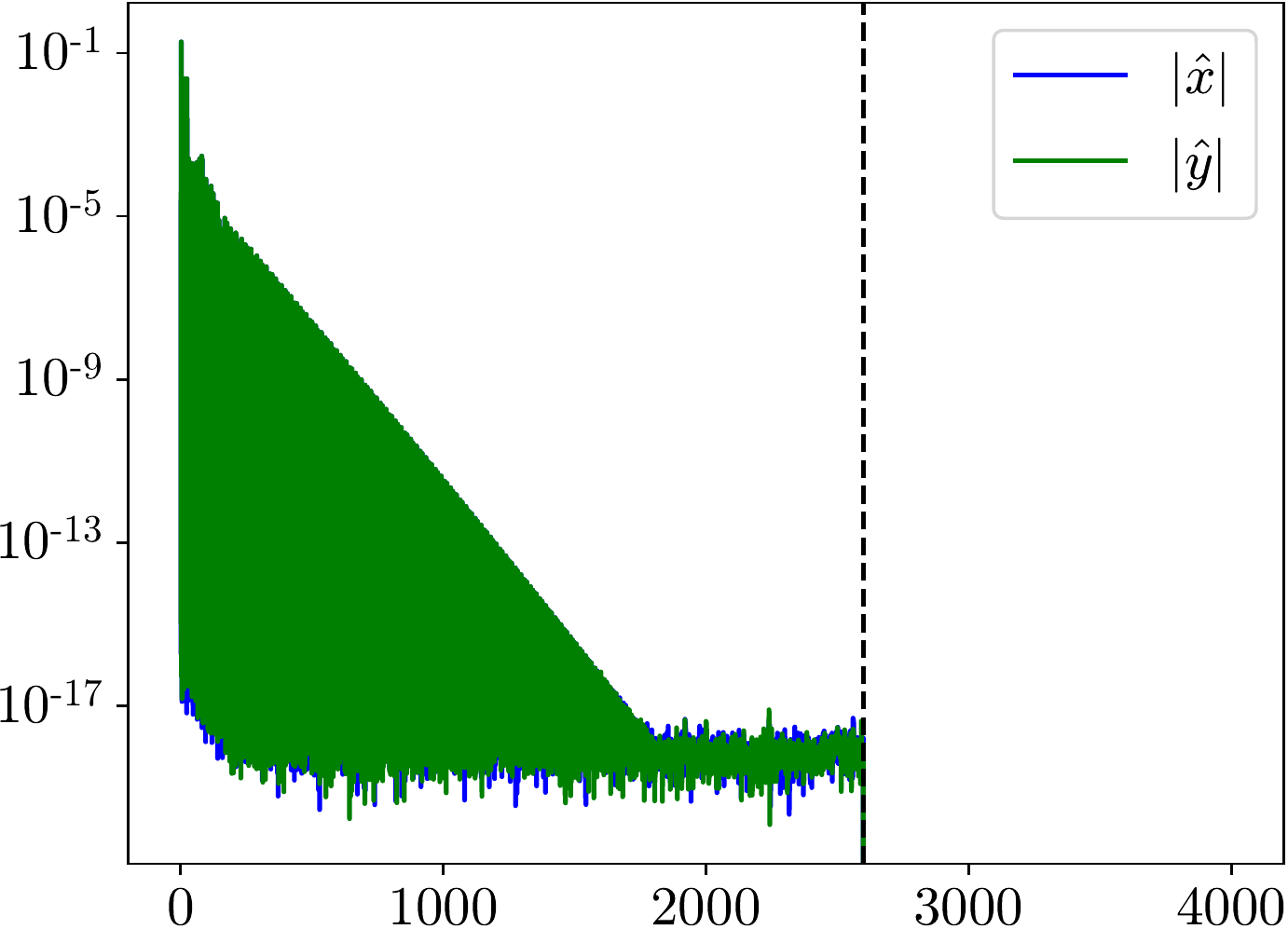}
    \caption{Fourier coefficients of the final curve}
  \end{subfigure}
  \caption{Fourier coefficients of $x(t)$ and $y(t)$ corresponding to
  Figure~\ref{fig:closed1}.
  The value of $n_{\text{coefs}}$ is indicated by a vertical dashed line.} 
  \label{fig:closed1cs2}
\end{figure}

The first closed curve example is obtained by sampling from the polar representation 
$(r(t)\cos{\varphi(t)}$, $r(t)\sin{\varphi(t)})$, where
\clearpage
  \begin{align}
  \label{eq:closedex1} 
  \varphi(t)&=2 \pi t, \notag\\
  r(t)&=(1+\frac{1}{\alpha} \cos{(18\varphi(t))} \sin{(4\varphi(t))},
  \end{align}
with $t \in [0,1]$, and $\alpha$ is a tuning parameter.
The sample data points are scaled so that both their width and height are $1$.
We sample the curve \eqref{eq:closedex1}
with $\alpha=2$, $N=8000$
and $n=100$ to obtain the initial curve in Figure~\ref{fig:closed1}(a). Applying 
the algorithm with 
$n_{\text{iters}}=70$, $h_{\text{filter}}=\frac{1}{35}$,
$\epsilon=10^{-16}$, $n_{\text{coefs}}=5200$,
$n_{\text{bands}}=12$, we obtained the filtered 
coefficients of $\theta(t)$ and $s'(t)$, as displayed in Figure~\ref{fig:closed1cs}.
We find that, after $n_{\text{stop}}=67$ iterations, $5200$ coefficients of $x(t)$ and
$y(t)$ are necessary to represent the smooth curve 
displayed in Figure~\ref{fig:closed1}(b),
to within an error of $E_{\text{samp}}=0.22453 \e{-14}$.
The magnitudes of the coefficients of both the initial 
and final curve are displayed 
in Figure~\ref{fig:closed1cs2}.

\begin{remark}
  \label{rem:dft}
Note that, since the DFT, $X_{-\frac{N}{2}},\ldots,X_{\frac{N}{2}-1}$,
of a real sequence, $x_0$, \ldots, $x_{N-1}$, satisfies the relation: 
  \begin{align}
  \label{eq:rmkdft}
X_k=\overline{X}_{-k}, \qquad k=-\frac{N}{2},\ldots,\frac{N}{2}-1,
  \end{align}
we only show the magnitudes of the coefficients,
for $k=0$, \ldots, $\frac{N}{2}-1$.  
\end{remark}

\begin{figure}[!h]
  \centering
  \begin{subfigure}[b]{0.28\linewidth}
    \includegraphics[width=\linewidth]{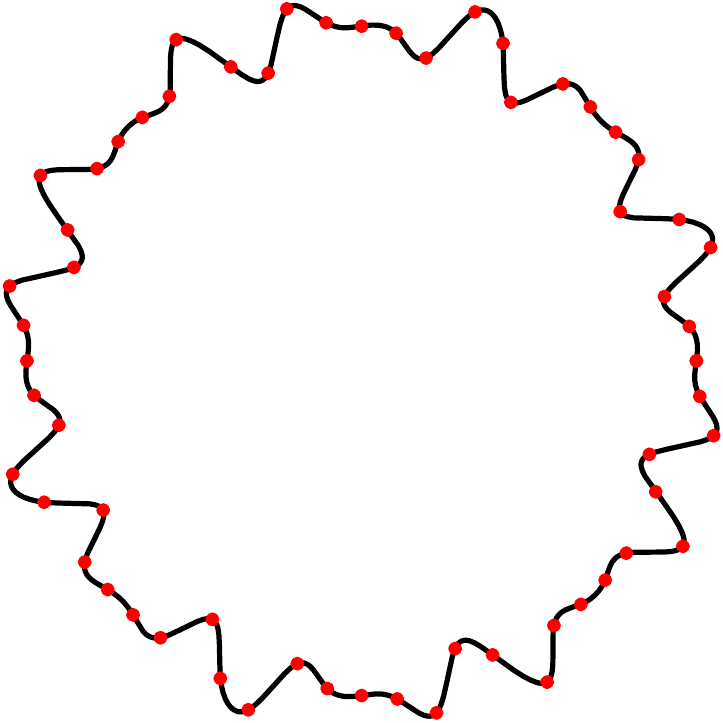}
    \caption{The curve before smoothing}
  \end{subfigure}\hspace{8mm}
  \begin{subfigure}[b]{0.28\linewidth}
    \includegraphics[width=\linewidth]{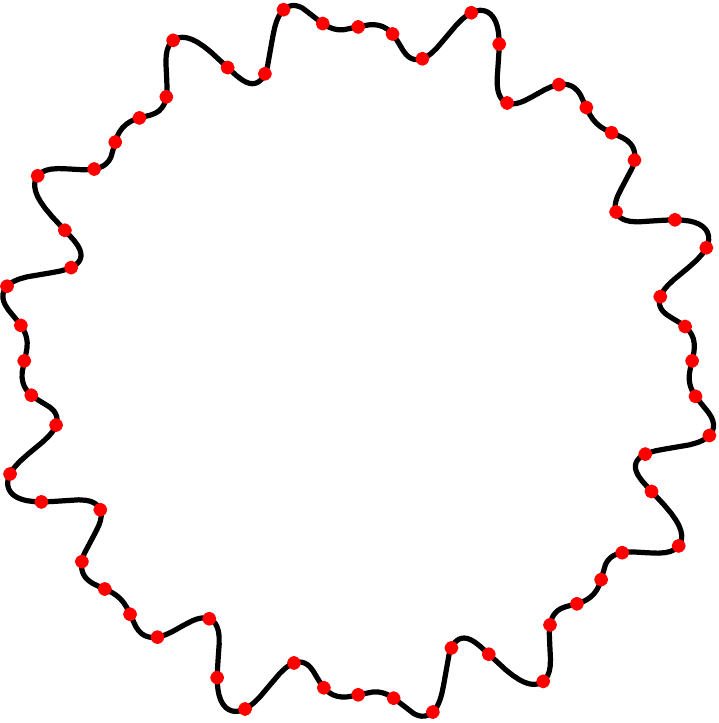}
    \caption{The curve after smoothing}
  \end{subfigure}
  \caption{The result of algorithm applied to \eqref{eq:closedex1} with $\alpha = 8$. The red dots 
 mark the sample points.}
  \label{fig:closed2}
\end{figure}

\begin{figure}[!h]
  \centering
  \begin{subfigure}[b]{0.45\linewidth}
    \includegraphics[width=\linewidth]{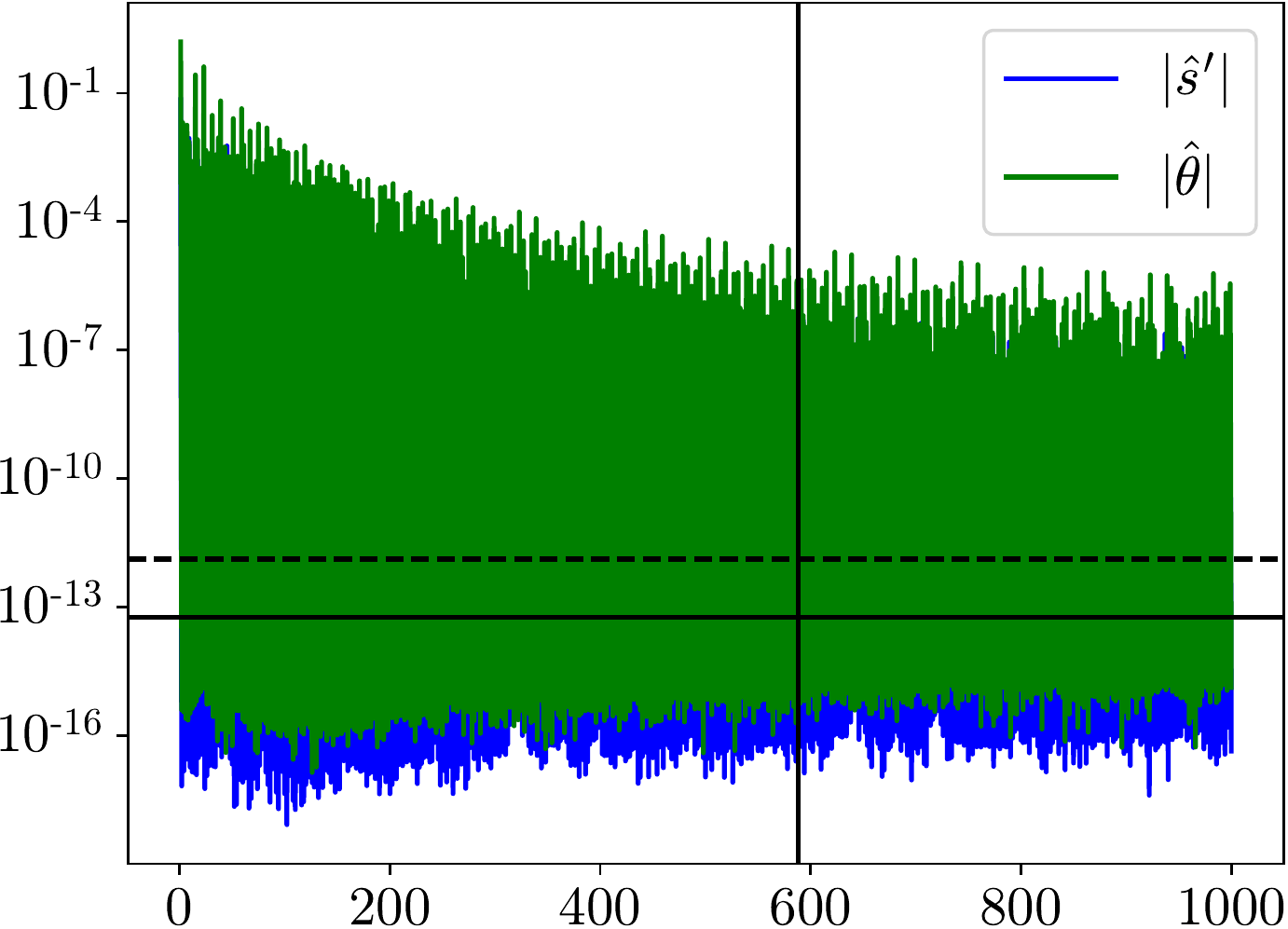}
    \caption{Before filtering}
  \end{subfigure}
  \begin{subfigure}[b]{0.45\linewidth}
    \includegraphics[width=\linewidth]{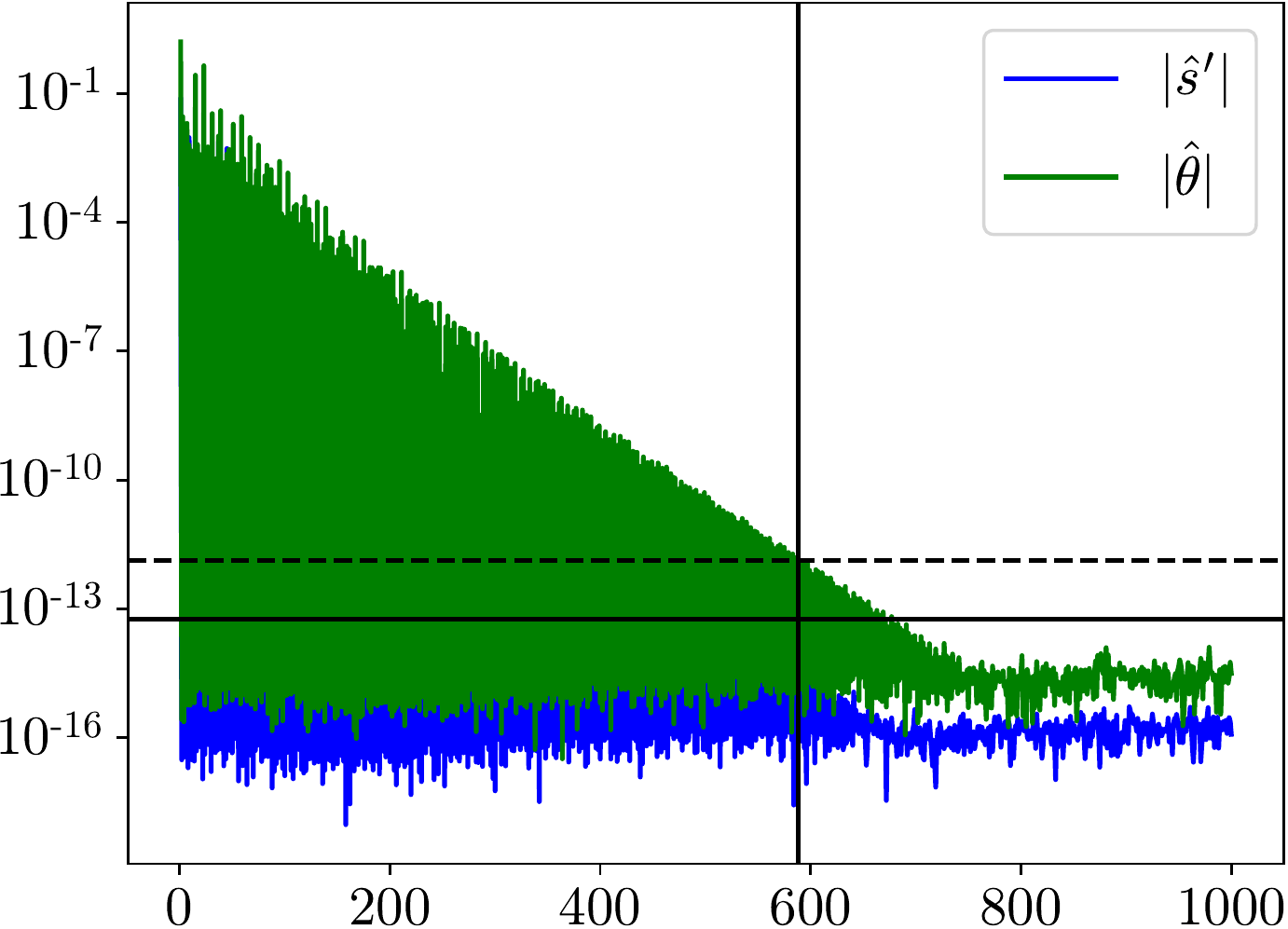}
    \caption{After filtering}
  \end{subfigure}
  \caption{Fourier coefficients of $s'(t)$ and $\theta(t)$ corresponding to
  Figure~\ref{fig:closed2}.
  The value of $\delta_{s'}$ is indicated by a horizontal solid line 
  and the value of $\delta_{\theta}$ is indicated by a horizontal dashed line.
  The $588$th coefficients of $s'(t)$ decays to $\delta_{s'}$, indicated by
  a vertical solid line.
  The $588$th coefficients of $\theta(t)$ decays to $\delta_{\theta}$,
  indicated by a vertical dashed line.}
  \label{fig:closed2cs}
\end{figure}

\begin{figure}[!h]
  \centering
  \begin{subfigure}[b]{0.45\linewidth}
    \includegraphics[width=\linewidth]{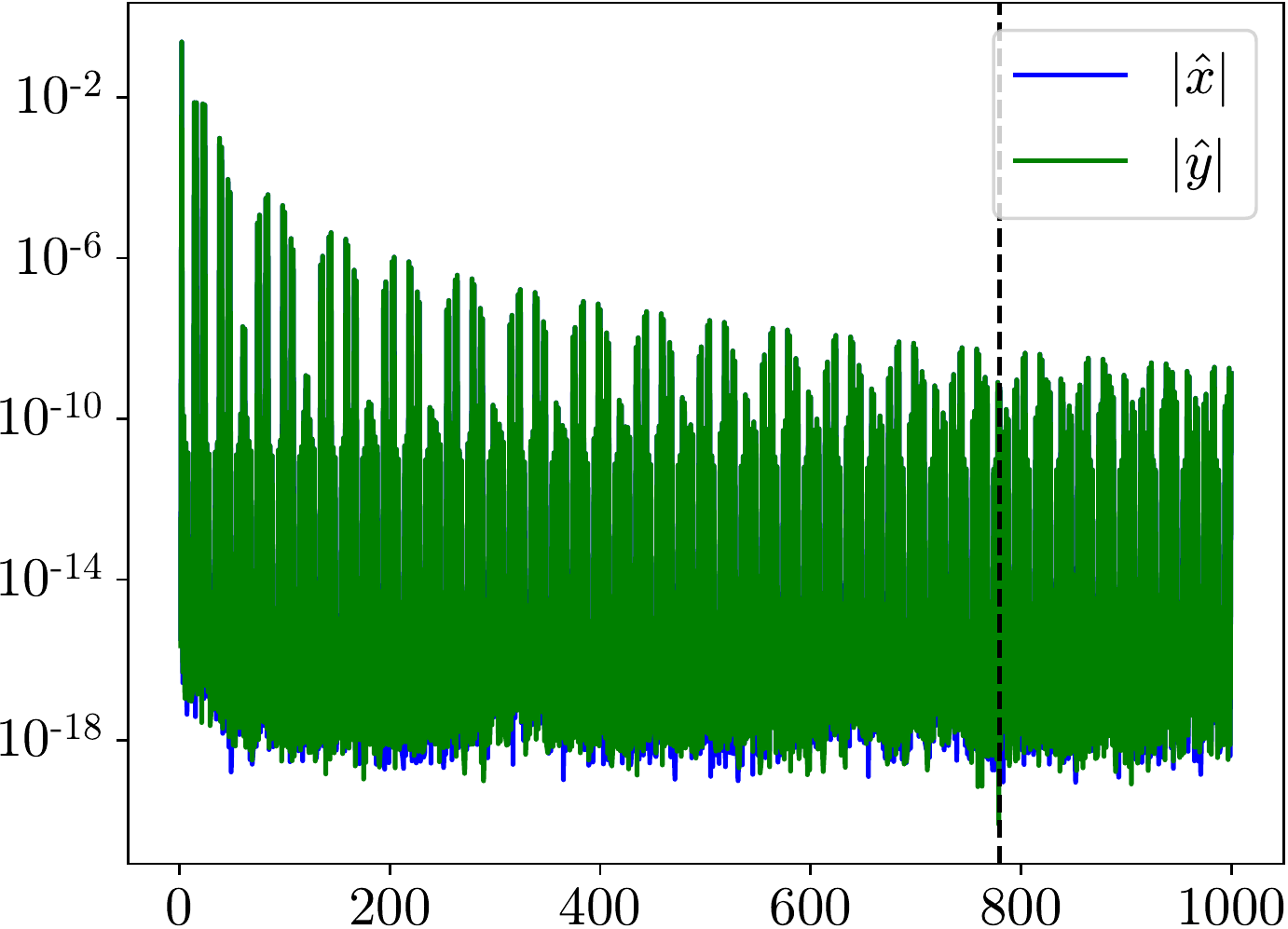}
    \caption{Fourier coefficients of the initial curve}
  \end{subfigure}
  \begin{subfigure}[b]{0.45\linewidth}
    \includegraphics[width=\linewidth]{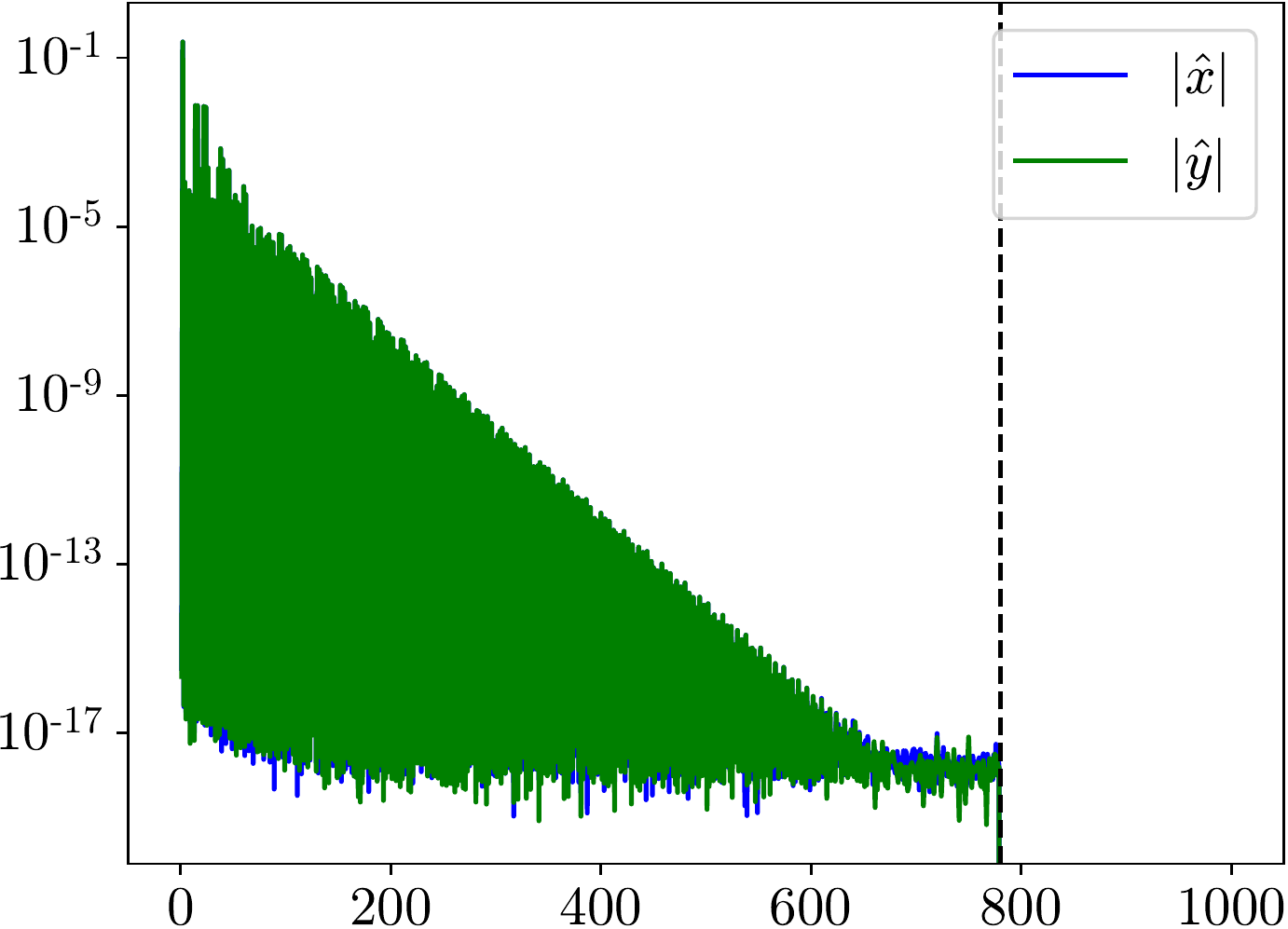}
    \caption{Fourier coefficients of the final curve}
  \end{subfigure}
  \caption{Fourier coefficients of $x(t)$ and $y(t)$ corresponding to
  Figure~\ref{fig:closed2}.
  The value of $n_{\text{coefs}}$ is indicated by a vertical dashed line.} 
  \label{fig:closed2cs2}
\end{figure}

Another example is shown in Figure~\ref{fig:closed2}, by sampling the curve \eqref{eq:closedex1} 
with $\alpha=8$, $N=2000$ and $n=60$.
The sample data points are scaled so that both their width and height are $1$.
It is obvious that the curve has fewer wobbles than the previous curve.
In this case, we set 
$n_{\text{iters}}=60$, $h_{\text{filter}}=\frac{1}{35}$,
$\epsilon=10^{-16}$, $n_{\text{coefs}}=1560$,
$n_{\text{bands}}=8$. The algorithm terminates after 
$n_{\text{stop}}=34$ iterations, and the error between the final curve
and the sample data points is $E_{\text{samp}}=0.11008 \e{-14}$.
The magnitudes of the coefficients of $s'(t)$ 
and $\theta(t)$ before and after filtering are displayed in Figure~\ref{fig:closed2cs}, and
the magnitudes of the coefficients of both the initial 
and final curve are displayed 
in Figure~\ref{fig:closed2cs2}.
While
the shapes of the curves appear similarly before and after filtering,
the coefficients change dramatically.

\begin{figure}[!h]
  \centering
  \begin{subfigure}[b]{0.35\linewidth}
    \includegraphics[width=\linewidth]{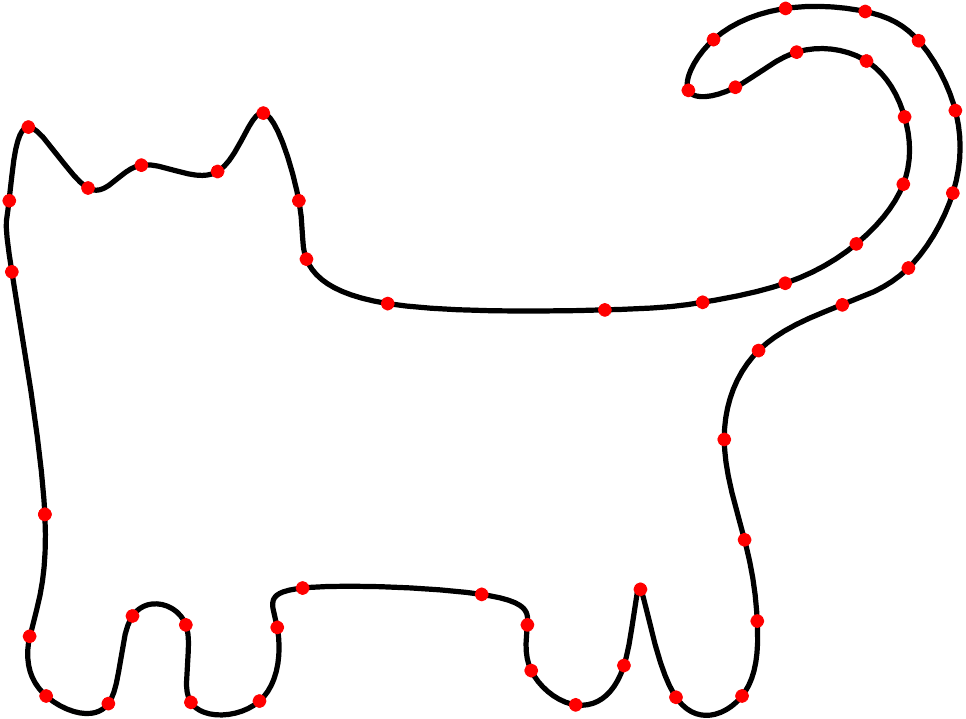}
    \caption{The curve before smoothing}
  \end{subfigure}\hspace{8mm}
  \begin{subfigure}[b]{0.35\linewidth}
    \includegraphics[width=\linewidth]{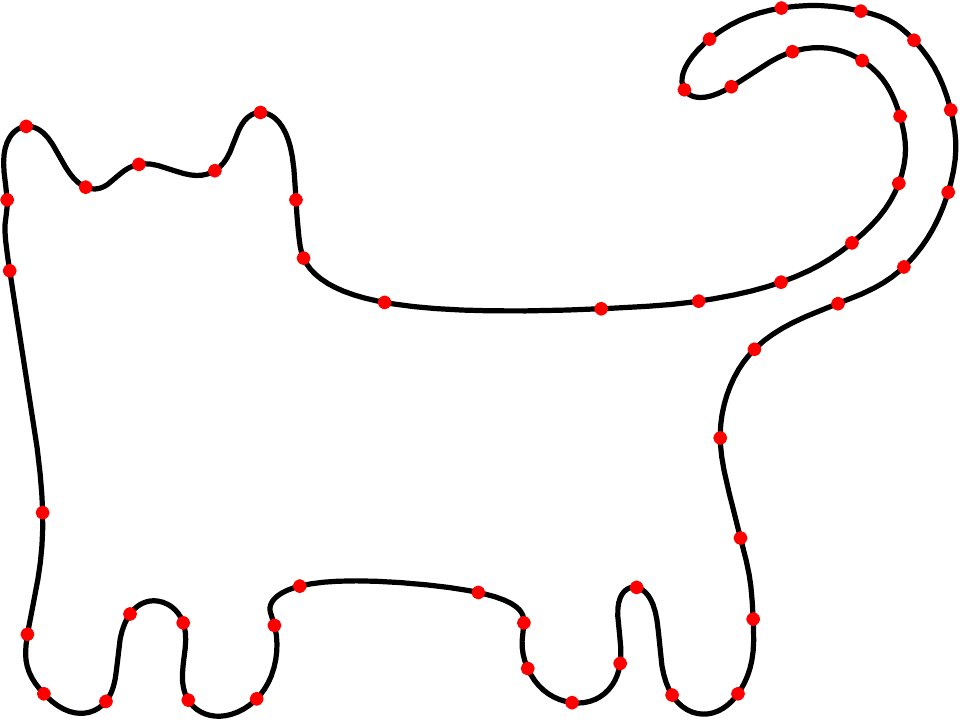}
    \caption{The curve after smoothing}
  \end{subfigure}
  \caption{A hand-drawn depiction of a cat shape. The red dots 
 mark the sample points.}
  \label{fig:closed4}
\end{figure}

\begin{figure}[!h]
  \centering
  \begin{subfigure}[b]{0.45\linewidth}
    \includegraphics[width=\linewidth]{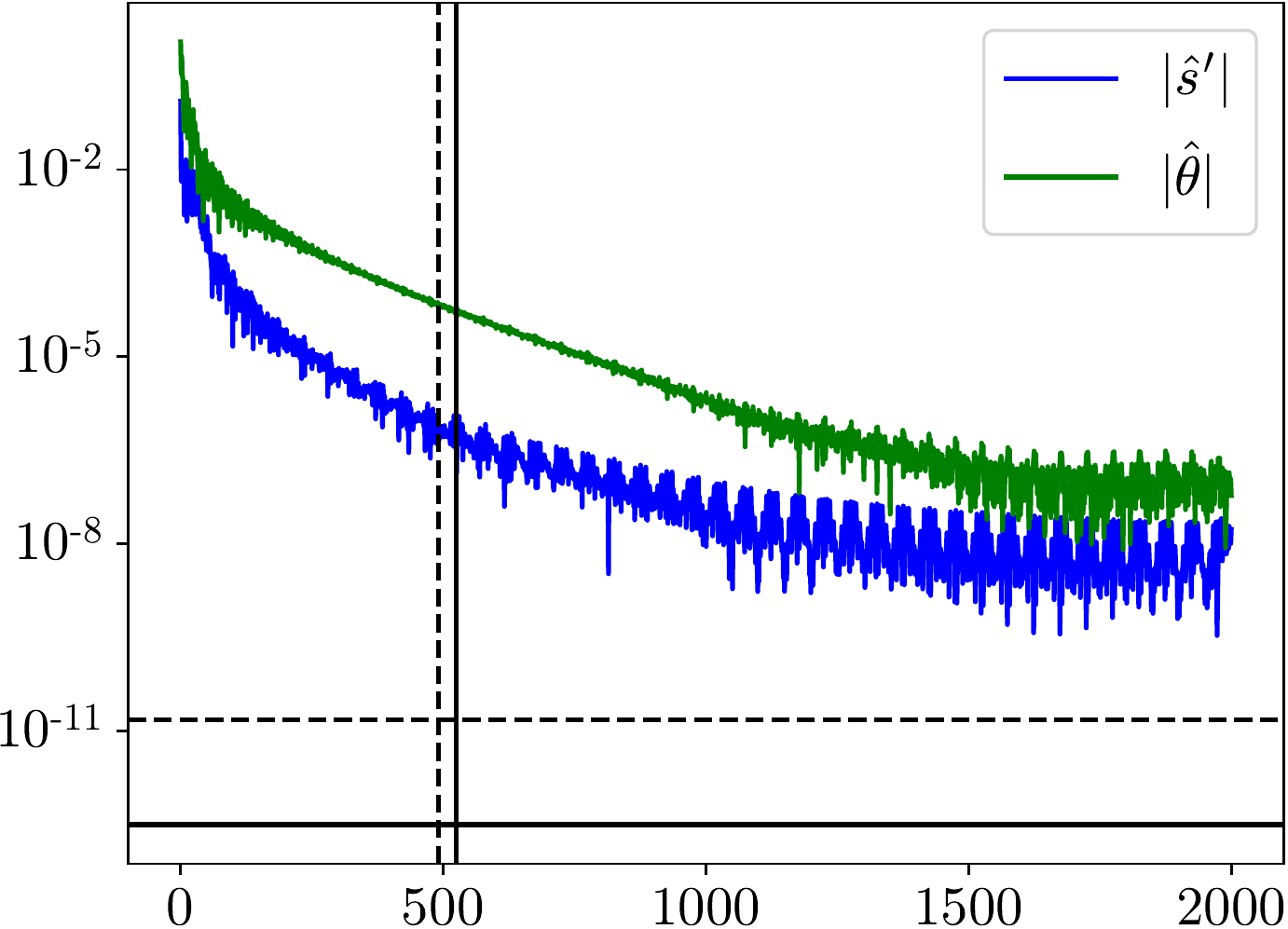}
    \caption{Before filtering}
  \end{subfigure}
  \begin{subfigure}[b]{0.45\linewidth}
    \includegraphics[width=\linewidth]{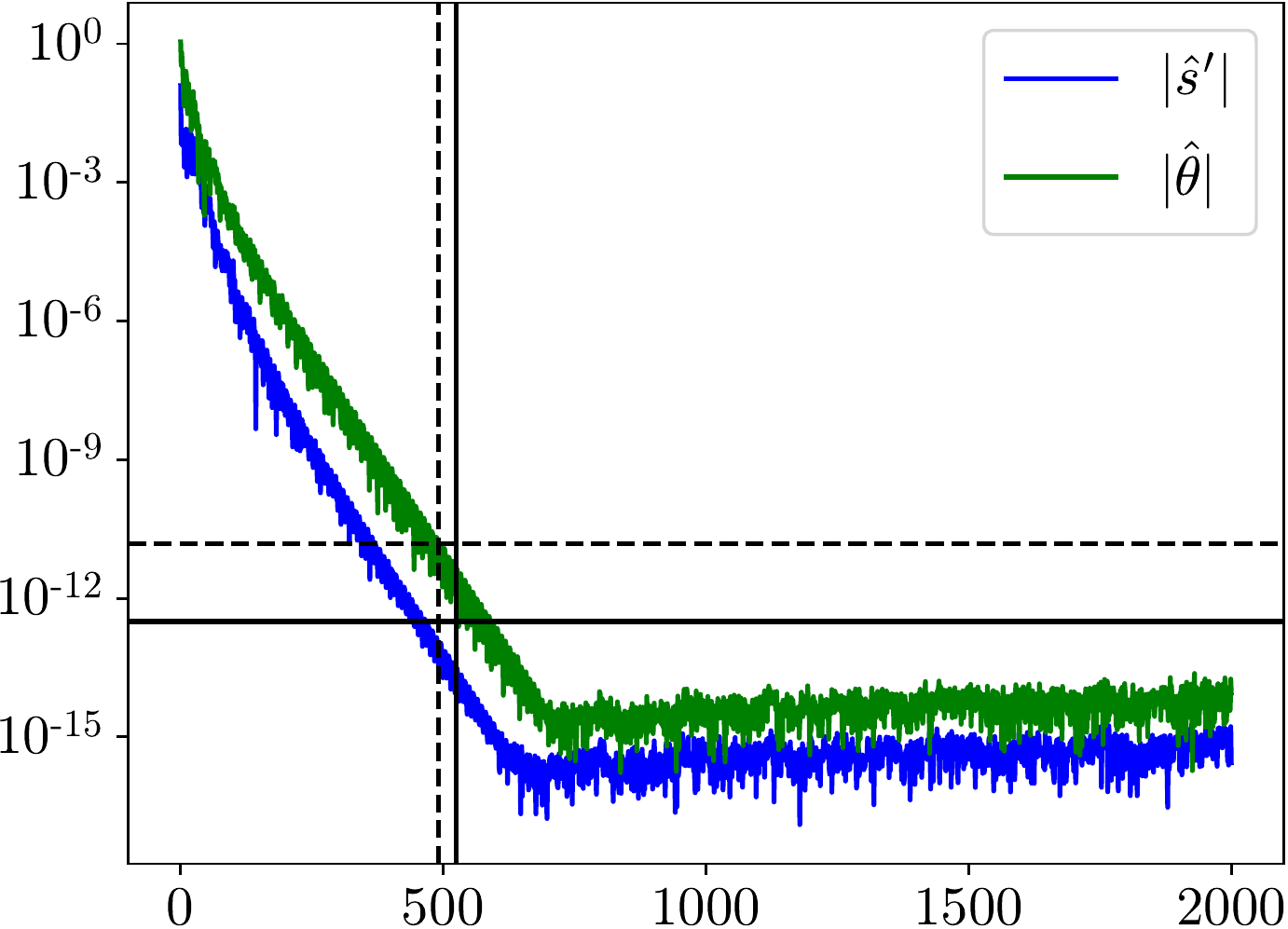}
    \caption{After filtering}
  \end{subfigure}
  \caption{Fourier coefficients of $s'(t)$ and $\theta(t)$ corresponding to
  Figure~\ref{fig:closed4}.
  The value of $\delta_{s'}$ is indicated by a horizontal solid line 
  and the value of $\delta_{\theta}$ is indicated by a horizontal dashed line.
  The $525$th coefficients of $s'(t)$ decays to $\delta_{s'}$, indicated by
  a vertical solid line.
  The $492$nd coefficients of $\theta(t)$ decays to $\delta_{\theta}$,
  indicated by a vertical dashed line.}
  \label{fig:closed4cs}
\end{figure}

\begin{figure}[!h]
  \centering
  \begin{subfigure}[b]{0.45\linewidth}
    \includegraphics[width=\linewidth]{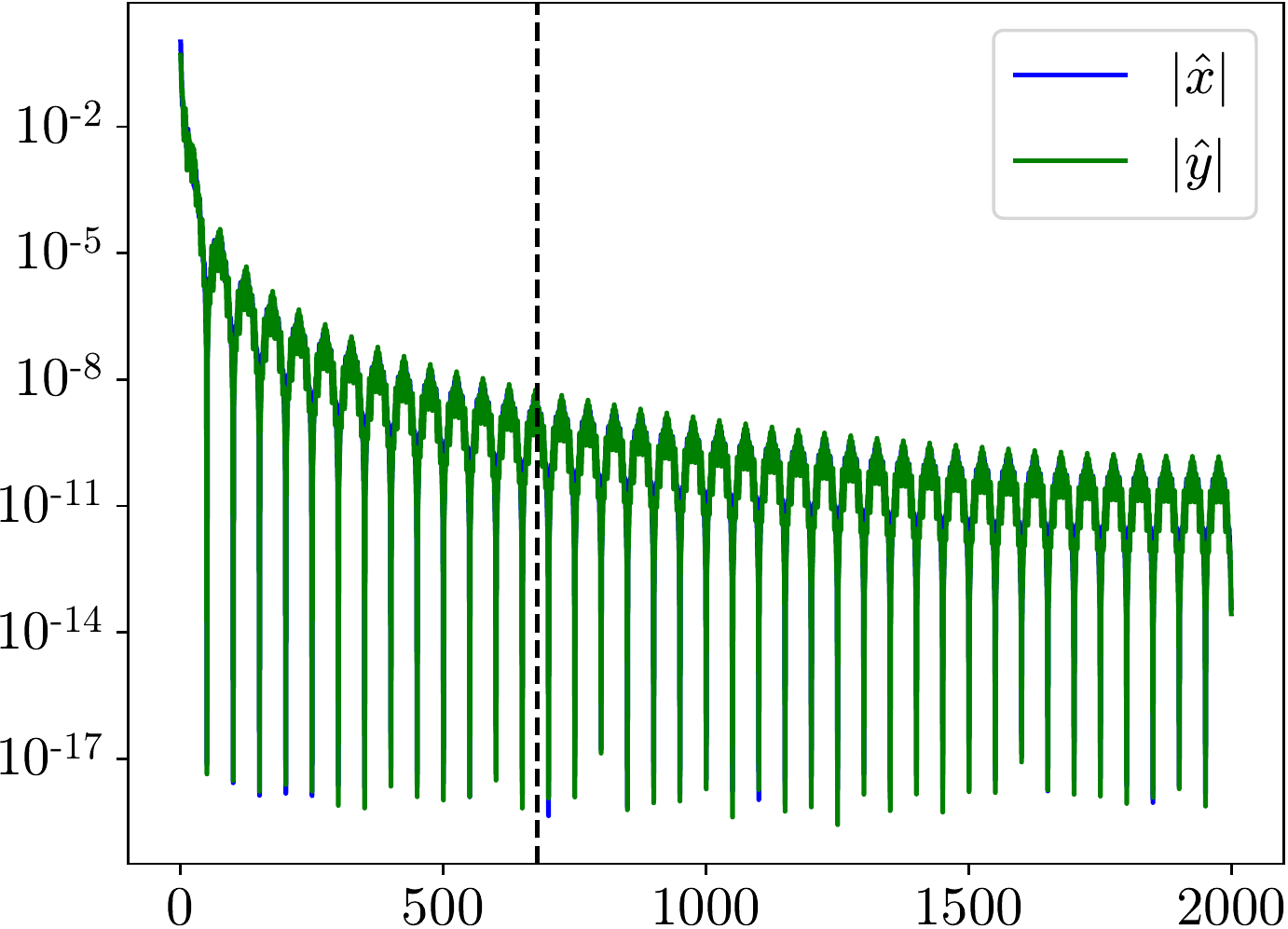}
    \caption{Fourier coefficients of the initial curve}
  \end{subfigure}
  \begin{subfigure}[b]{0.45\linewidth}
    \includegraphics[width=\linewidth]{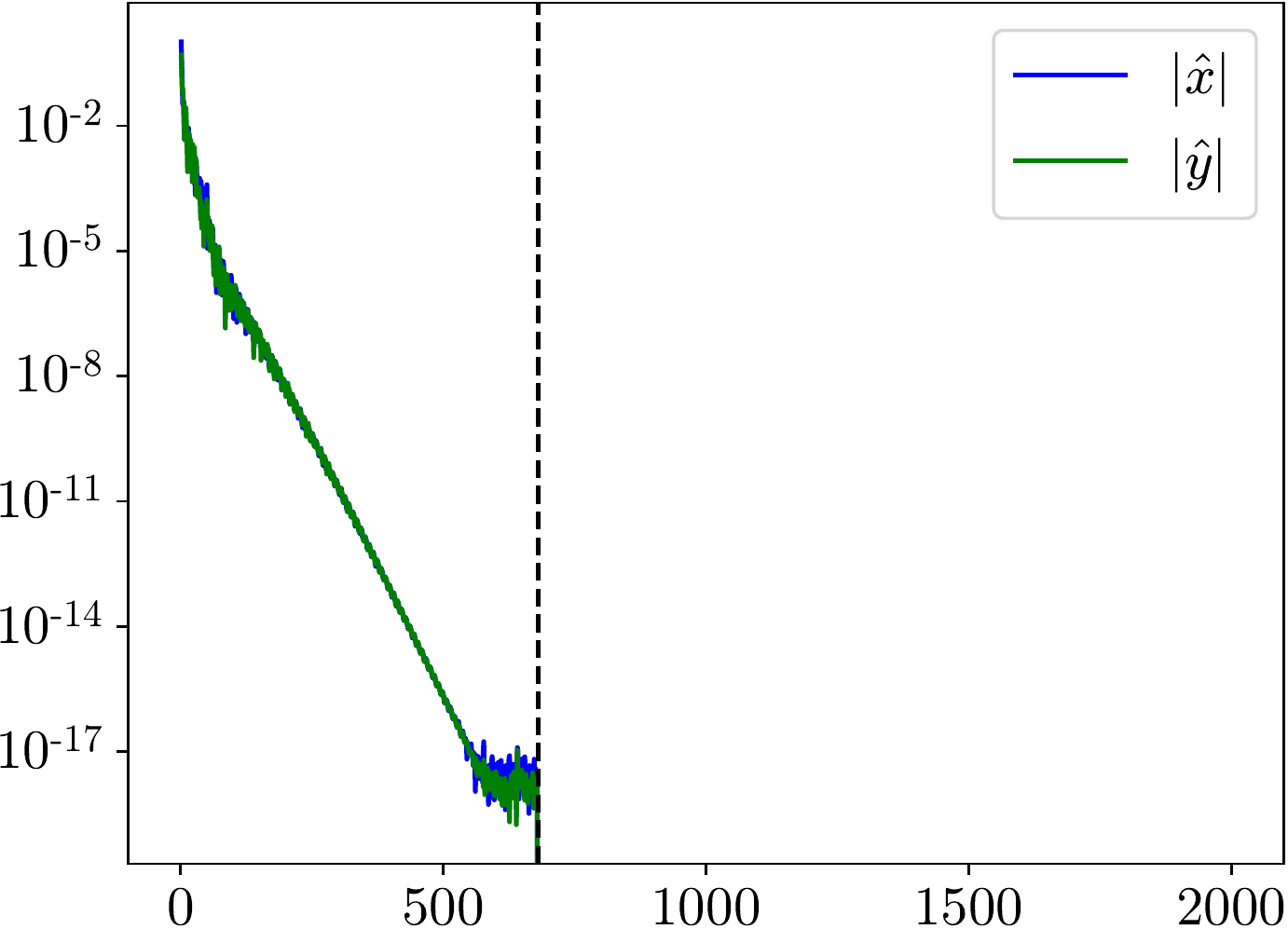}
    \caption{Fourier coefficients of the final curve}
  \end{subfigure}
  \caption{Fourier coefficients of $x(t)$ and $y(t)$ corresponding to
  Figure~\ref{fig:closed4}.
  The value of $n_{\text{coefs}}$ is indicated by a vertical dashed line.} 
  \label{fig:closed4cs2}
\end{figure}

\begin{figure}[!h]
  \centering
  \begin{subfigure}[b]{0.35\linewidth}
    \includegraphics[width=\linewidth]{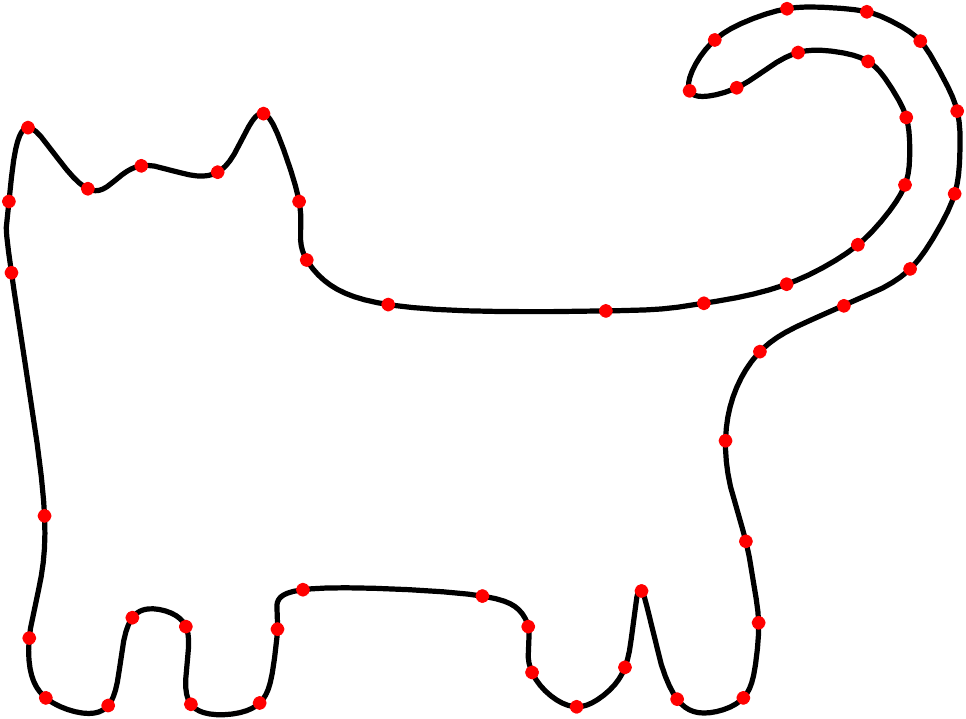}
    \caption{The interpolating curve. The red dots 
       mark the sample points.}
  \end{subfigure}\hspace{10mm}
  \begin{subfigure}[b]{0.45\linewidth}
    \includegraphics[width=\linewidth]{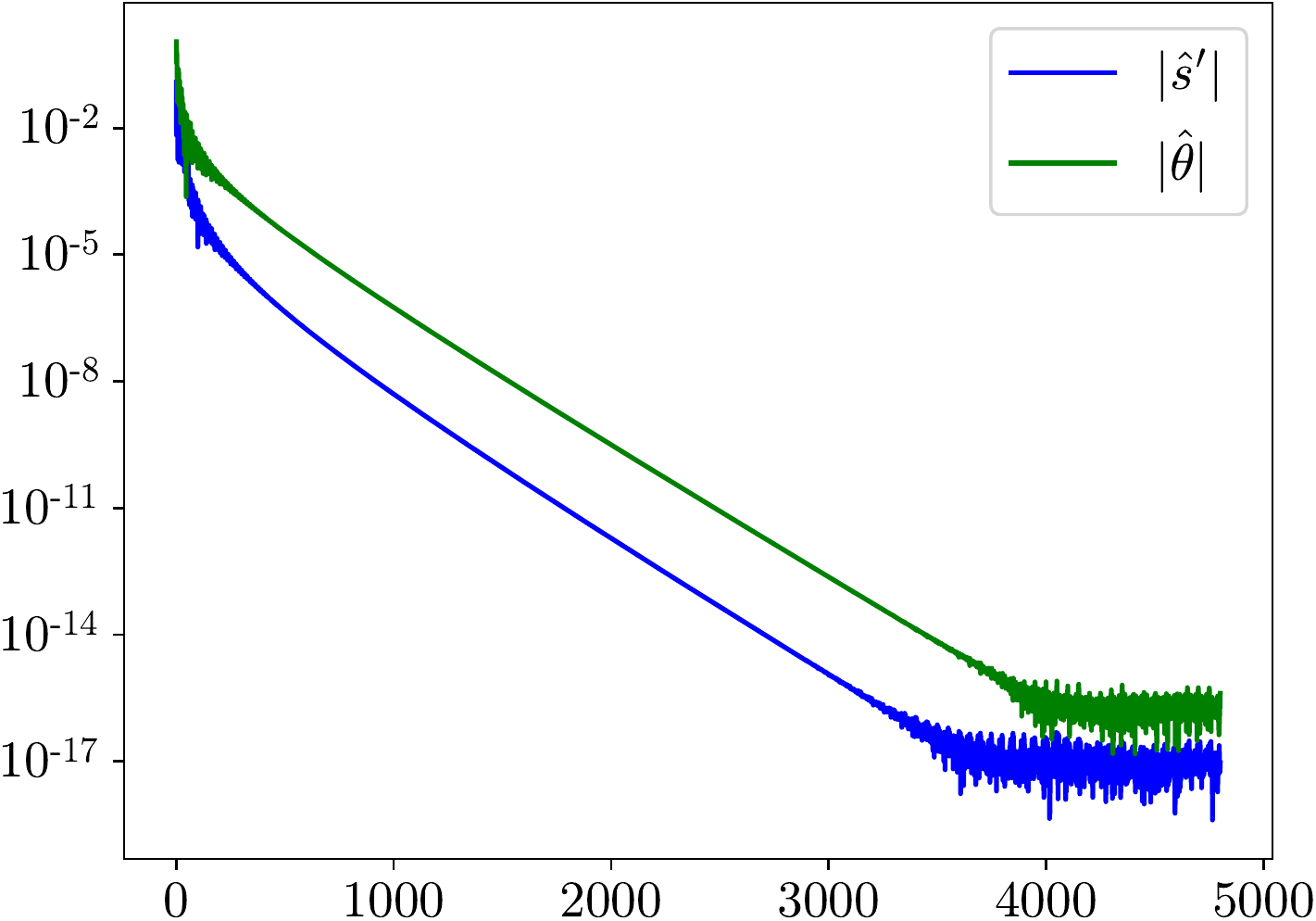}
    \caption{Fourier coefficients of $s'(t)$ and $\theta(t)$ corresponding to 
  Figure~\ref{fig:closed4_2}(a).}
  \end{subfigure}
    \caption{The result of algorithm in \cite{zhang} applied to the 
    same data points in Figure~\ref{fig:closed4}.}
  \label{fig:closed4_2}
\end{figure}

The example shown in Figure~\ref{fig:closed4} has the shape
of a cat, corresponding to $n=50$ sample points.
We scale the sample data points so that their height is $1$,
and discretize the curve with $N=4000$, and set
$n_{\text{iters}}=100$, $h_{\text{filter}}=\frac{1}{45}$,
$\epsilon=10^{-15}$, $n_{\text{coefs}}=1360$,
$n_{\text{bands}}=4$. After 
$n_{\text{stop}}=98$ iterations,
$1360$ coefficients
are sufficient to represent 
the curve, and
the error between the final curve
and the sample data points is $E_{\text{samp}}=0.64403 \e{-14}$.
The magnitudes of the coefficients of $s'(t)$ 
and $\theta(t)$ before and after filtering are displayed in Figure~\ref{fig:closed4cs}, and
the magnitudes of the coefficients of both the initial 
and final curve are displayed 
in Figure~\ref{fig:closed4cs2}.
We observe that 
the sharp edges on the curve in Figure~\ref{fig:closed4}(a) becomes
more rounded, and the
resulting curve more closely resembles the shape of a cat.
Meanwhile,
Figure~\ref{fig:closed4_2}(a) is the result of the algorithm in \cite{zhang} applied 
to the same data points, with $a=0.4$. The number of coefficients that are
necessary to represent $s'(t)$ and $\theta(t)$ are displayed in 
Figure \ref{fig:closed4_2}(b).
As we can observe from
Figure~\ref{fig:closed4_2}(a), the curve contains 
sharp corners, and its visual smoothness is similar to the initial curve displayed
in Figure~\ref{fig:closed4}(a), prior to applying our algorithm.
In contrast, our algorithm eliminates high curvature areas, resulting in a much smoother
appearance. 

\begin{figure}[!h]
  \centering
  \begin{subfigure}[b]{0.38\linewidth}
    \includegraphics[width=\linewidth]{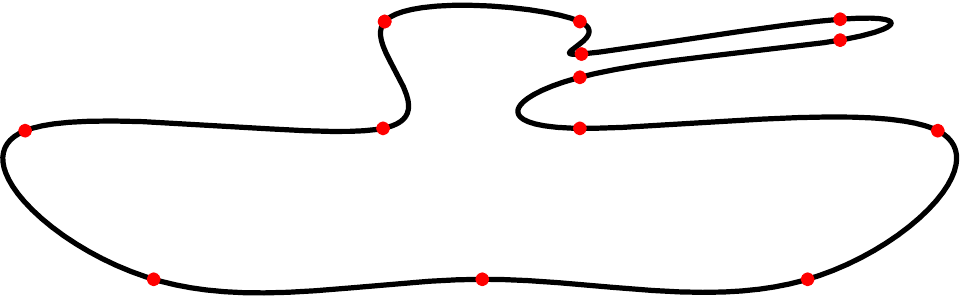}
    \caption{The curve before smoothing}
  \end{subfigure}\hspace{8mm}
  \begin{subfigure}[b]{0.38\linewidth}
    \includegraphics[width=\linewidth]{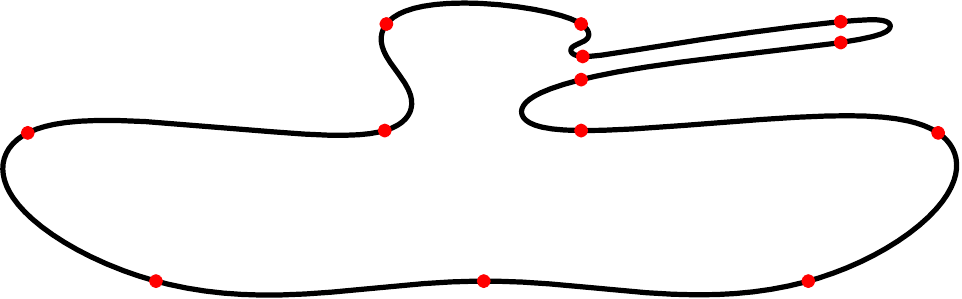}
    \caption{The curve after smoothing}
  \end{subfigure}
  \caption{The result of algorithm applied to Figure $4.5$ in~\cite{bandlimited}. The red dots 
 mark the sample points.}
  \label{fig:closed3}
\end{figure}

Inspired by Figure 4.5 and Figure 4.3 in~\cite{bandlimited},
we apply our algorithm 
to the same sample data points to show that our algorithm produces a smoother
curve and represents the curve with fewer coefficients.
We start with the example in Figure 4.5.
The sample data points are scaled so that their width is $1$.
With the parameters
$N=1600$, $n=13$, $n_{\text{iters}}=80$, $h_{\text{filter}}=\frac{1}{35}$,
$\epsilon=10^{-16}$, $n_{\text{coefs}}=840$, and
$n_{\text{bands}}=4$,
our algorithm produces a curve represented by only $840$ coefficients,
while the algorithm of~\cite{bandlimited}
produces a curve represented by $2 \cdot 25,000 = 50,000$
coefficients.
The algorithm terminates at the $n_{\text{stop}}=58$th iteration,
and the resulting curve passes through 
the sample data points within an error of $E_{\text{samp}}=0.15713 \e{-13}$.
Notice that the corners in Figure $4.5$ are eliminated and the final curve 
in Figure~\ref{fig:closed3}(b) looks 
much smoother.   
We display the magnitudes of the coefficients of $s'(t)$ 
and $\theta(t)$ before and after filtering in Figure~\ref{fig:closed3cs}, and
the magnitudes of the coefficients of both the initial 
and final curve 
in Figure~\ref{fig:closed3cs2}.

Since the initial curve is constructed by using
smooth splines to connect the sample data points,
and our algorithm filters the curve further during the filtering process,
this causes the shape of the curve to deviate from that of
Figure 4.5 in~\cite{bandlimited}.
In order to preserve the shape of the curve in Figure 4.5,
we increase the number of sample data points.
The sample data points are scaled so that their width is $1$.
Applying the algorithm to the curve displayed in Figure~\ref{fig:closed6}(a), with  
$N=4000$, $n=59$, $n_{\text{iters}}=80$, $h_{\text{filter}}=\frac{1}{45}$,
$\epsilon=10^{-16}$, $n_{\text{coefs}}=1700$,
$n_{\text{bands}}=4$,
we obtain $E_{\text{samp}}=0.83564 \e{-13}$ at the 
$n_{\text{stop}}=72$nd iteration. 
The magnitudes of the coefficients of $s'(t)$ 
and $\theta(t)$ before and after filtering are displayed in Figure~\ref{fig:closed6cs}, and
the magnitudes of the coefficients of both the initial 
and final curve are displayed in 
in Figure~\ref{fig:closed6cs2}.
It is noteworthy that the shape of the curve has been preserved, while
achieving a smoother curve than in \cite{bandlimited}.
Moreover, the number of coefficients required to represent the curve increases only
moderately when compared to those shown in Figure~\ref{fig:closed3cs}.

For Figure $4.3$ in~\cite{bandlimited}, 
we scale the sample data points so that their height is $1$,
and apply the algorithm to the reproduced curve 
in Figure~\ref{fig:closed5}(a), 
with
$N=2000$, $n=41$,
$n_{\text{iters}}=70$, $h_{\text{filter}}=\frac{1}{40}$,
$\epsilon=10^{-16}$, $n_{\text{coefs}}=680$,
$n_{\text{bands}}=4$.
Although the difference can not be distinguished visually, 
after $n_{\text{stop}}=55$ iterations, $680$ coefficients
are necessary to represent 
the curve, to within an error of $E_{\text{samp}}=0.15102\e{-13}$.
The magnitudes of the coefficients of $s'(t)$ and $\theta(t)$ 
are displayed in 
in Figure~\ref{fig:closed5cs}, and
the magnitudes of the coefficients of both the initial and final curve  
are displayed in 
in Figure~\ref{fig:closed5cs2}. Notice that 
the algorithm of~\cite{bandlimited} requires approximately $2 \cdot 7000 = 14,000$
coefficients to fit a curve passing through the same sample data points.
 
The runtimes per iteration for the first two closed curve cases are displayed in Table ~\ref{tab:table2}.
We observe that, as in the open curve case, the runtimes are not strictly proportional to $N$.




\begin{figure}[!h]
  \centering
  \includegraphics[scale=0.45]{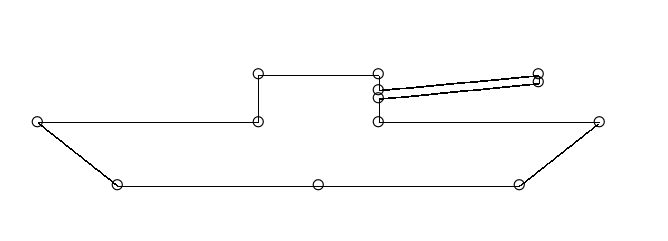}
  \caption{Figure $4.5$ in \cite{bandlimited}}
  \label{fig:closed3ex}
\end{figure}
\begin{figure}[!h]
  \centering
  \begin{subfigure}[b]{0.45\linewidth}
    \includegraphics[width=\linewidth]{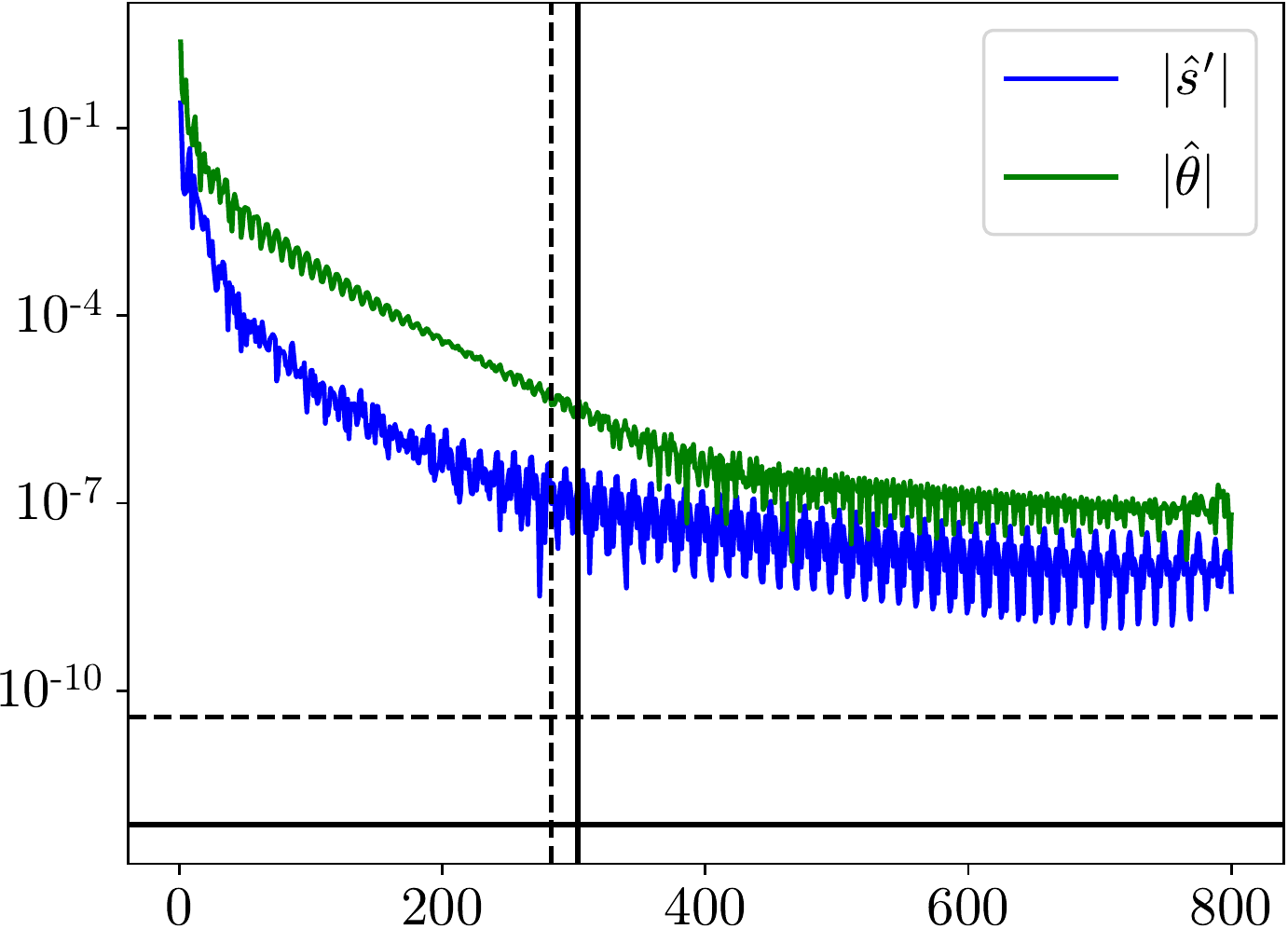}
    \caption{Before filtering}
  \end{subfigure}
  \begin{subfigure}[b]{0.45\linewidth}
    \includegraphics[width=\linewidth]{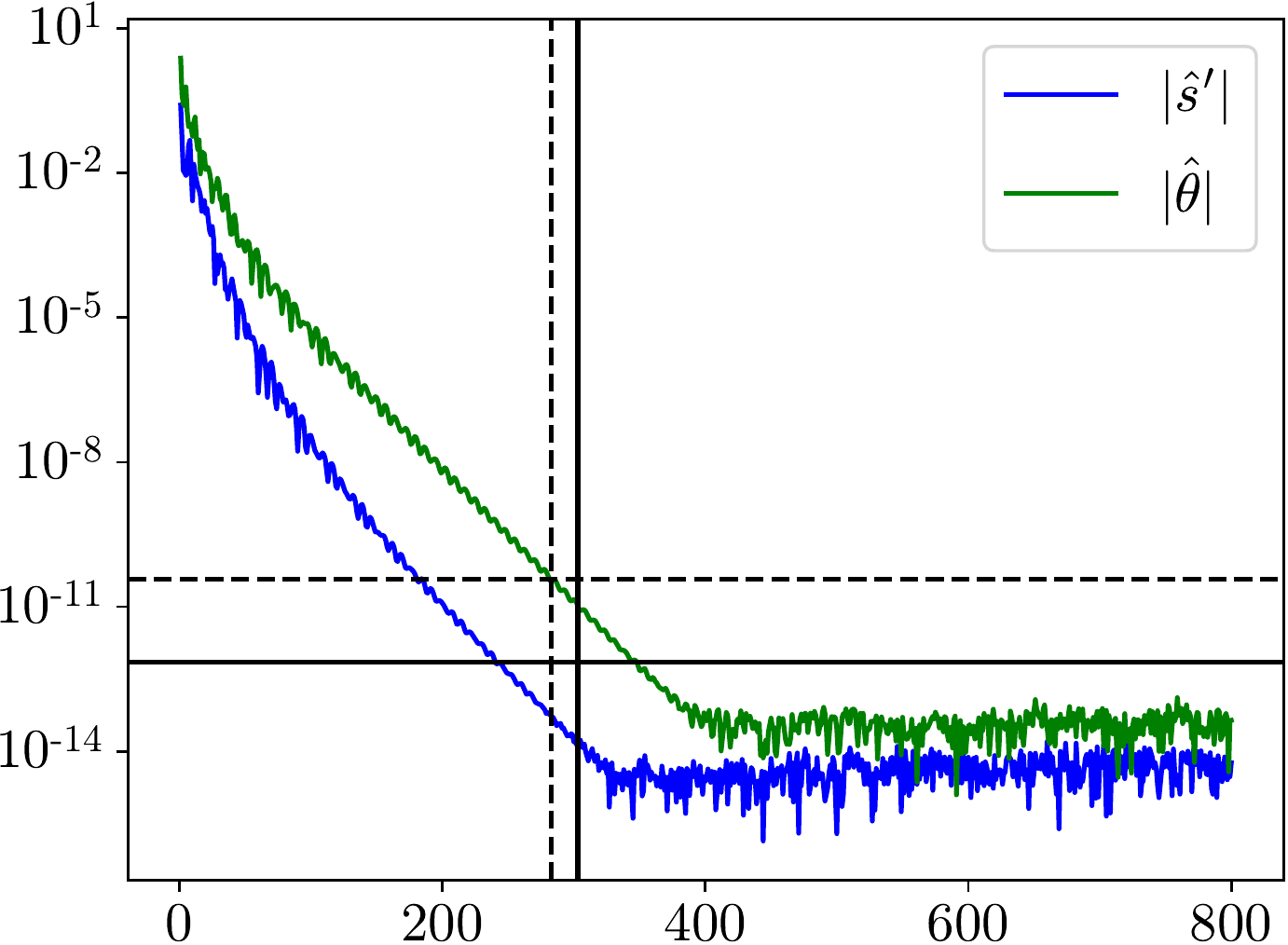}
    \caption{After filtering}
  \end{subfigure}
  \caption{Fourier coefficients of $s'(t)$ and $\theta(t)$ corresponding to
  Figure~\ref{fig:closed3}.
  The value of $\delta_{s'}$ is indicated by a horizontal solid line 
  and the value of $\delta_{\theta}$ is indicated by a horizontal dashed line.
  The $303$rd coefficients of $s'(t)$ decays to the $\delta_{s'}$, indicated
  by a vertical solid line. 
  The $283$rd coefficients of $\theta(t)$ decays to the $\delta_{\theta}$,
  indicated by a vertical dashed line.}
  \label{fig:closed3cs}
\end{figure}
\begin{figure}[!h]
  \centering
  \begin{subfigure}[b]{0.45\linewidth}
    \includegraphics[width=\linewidth]{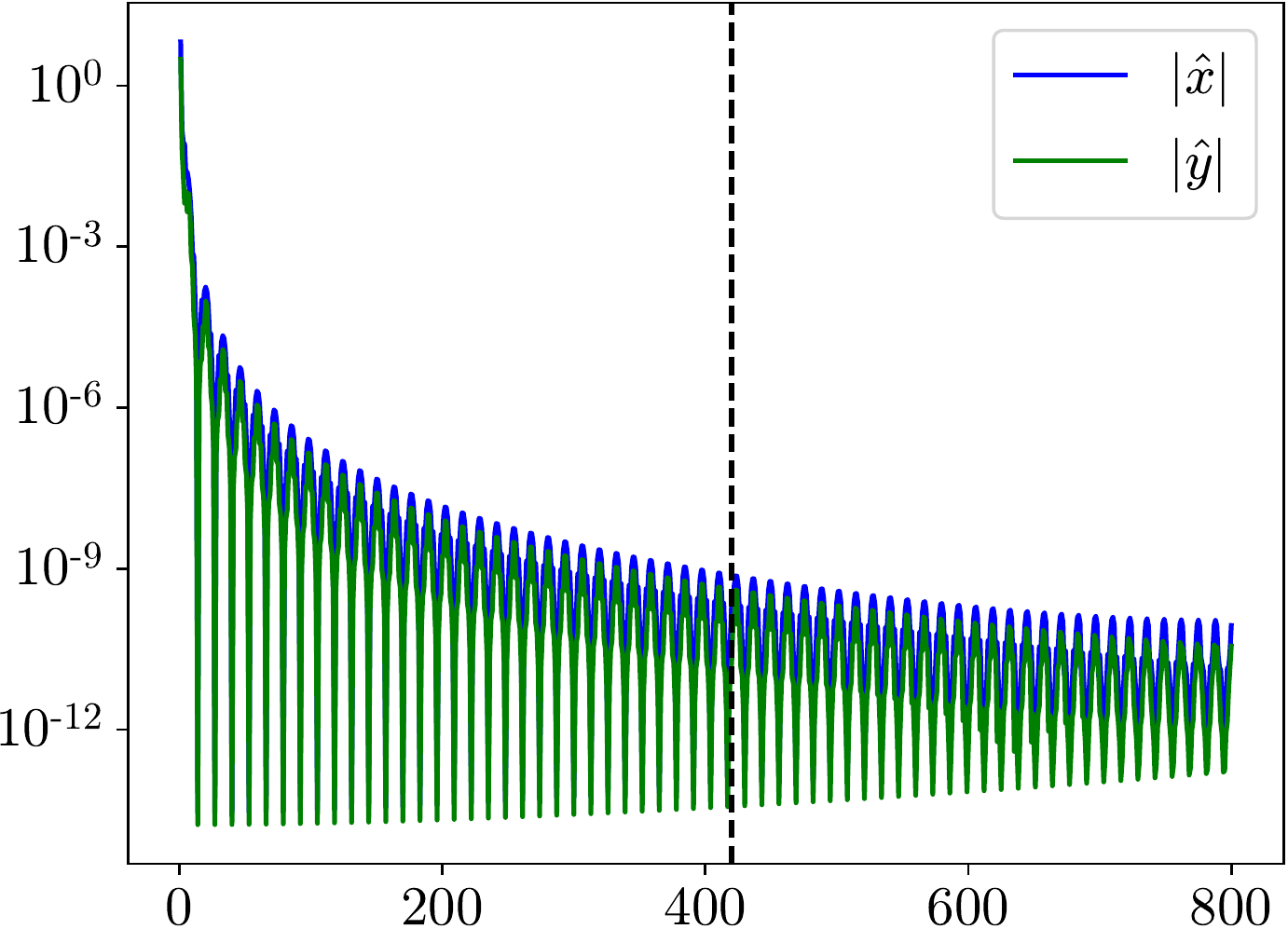}
    \caption{Fourier coefficients of the initial curve}
  \end{subfigure}
  \begin{subfigure}[b]{0.45\linewidth}
    \includegraphics[width=\linewidth]{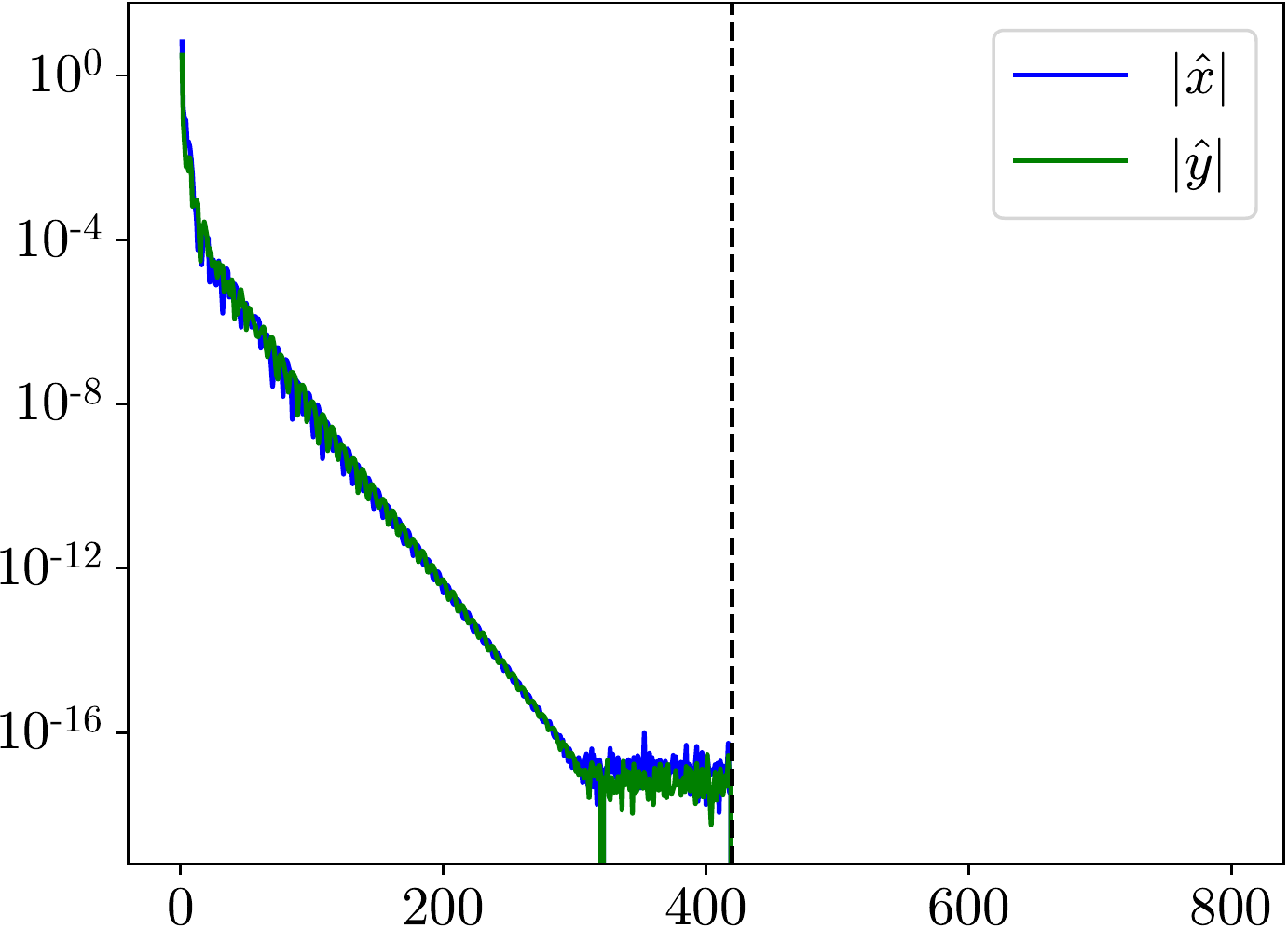}
    \caption{Fourier coefficients of the final curve}
  \end{subfigure}
  \caption{Fourier coefficients of $x(t)$ and $y(t)$ corresponding to
  Figure~\ref{fig:closed3}.
  The value of $n_{\text{coefs}}$ is indicated by a vertical dashed line.} 
  \label{fig:closed3cs2}
\end{figure}
\begin{figure}[!h]
  \centering
  \begin{subfigure}[b]{0.38\linewidth}
    \includegraphics[width=\linewidth]{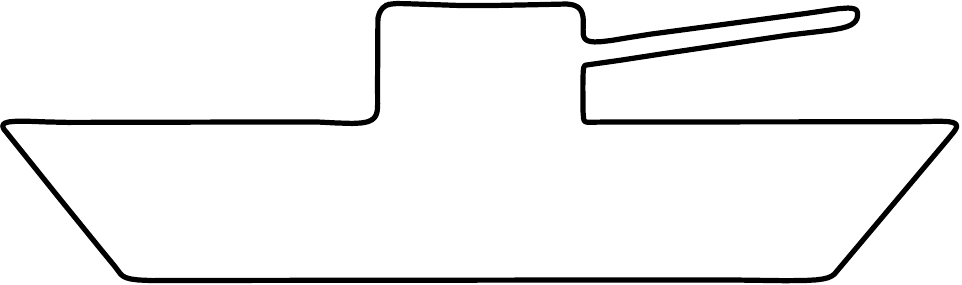}
    \caption{The curve before smoothing}
  \end{subfigure}\hspace{8mm}
  \begin{subfigure}[b]{0.38\linewidth}
    \includegraphics[width=\linewidth]{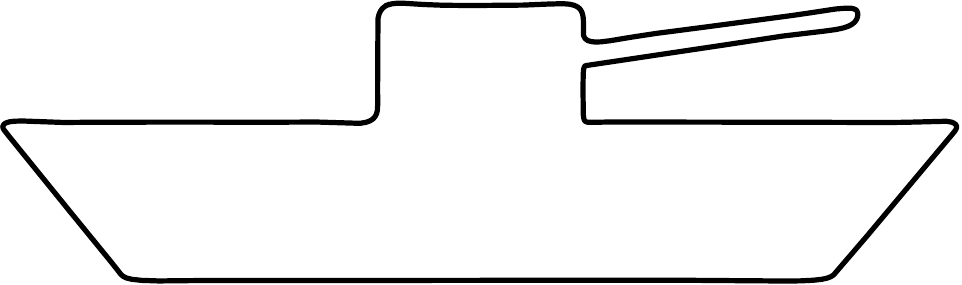}
    \caption{The curve after smoothing}
  \end{subfigure}
  \caption{The result of algorithm applied to Figure $4.5$ in~\cite{bandlimited},
   with more sample data points.
   Due to the large quantity and non-uniform distribution of 
   the sample data points, we choose not to display them in the plot.}
  \label{fig:closed6}
\end{figure}

\begin{figure}[!h]
  \centering
  \begin{subfigure}[b]{0.45\linewidth}
    \includegraphics[width=\linewidth]{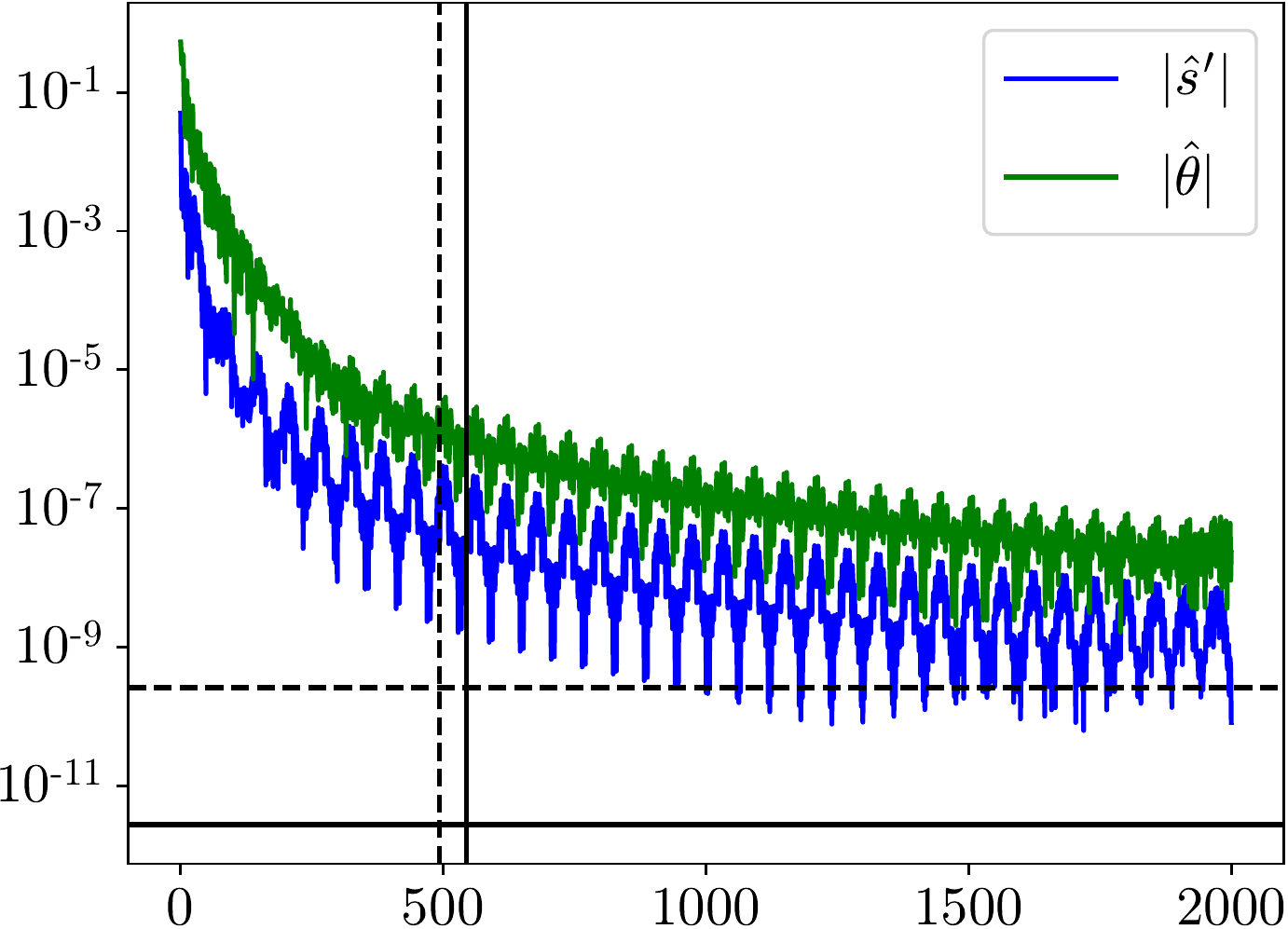}
    \caption{Before filtering}
  \end{subfigure}
  \begin{subfigure}[b]{0.45\linewidth}
    \includegraphics[width=\linewidth]{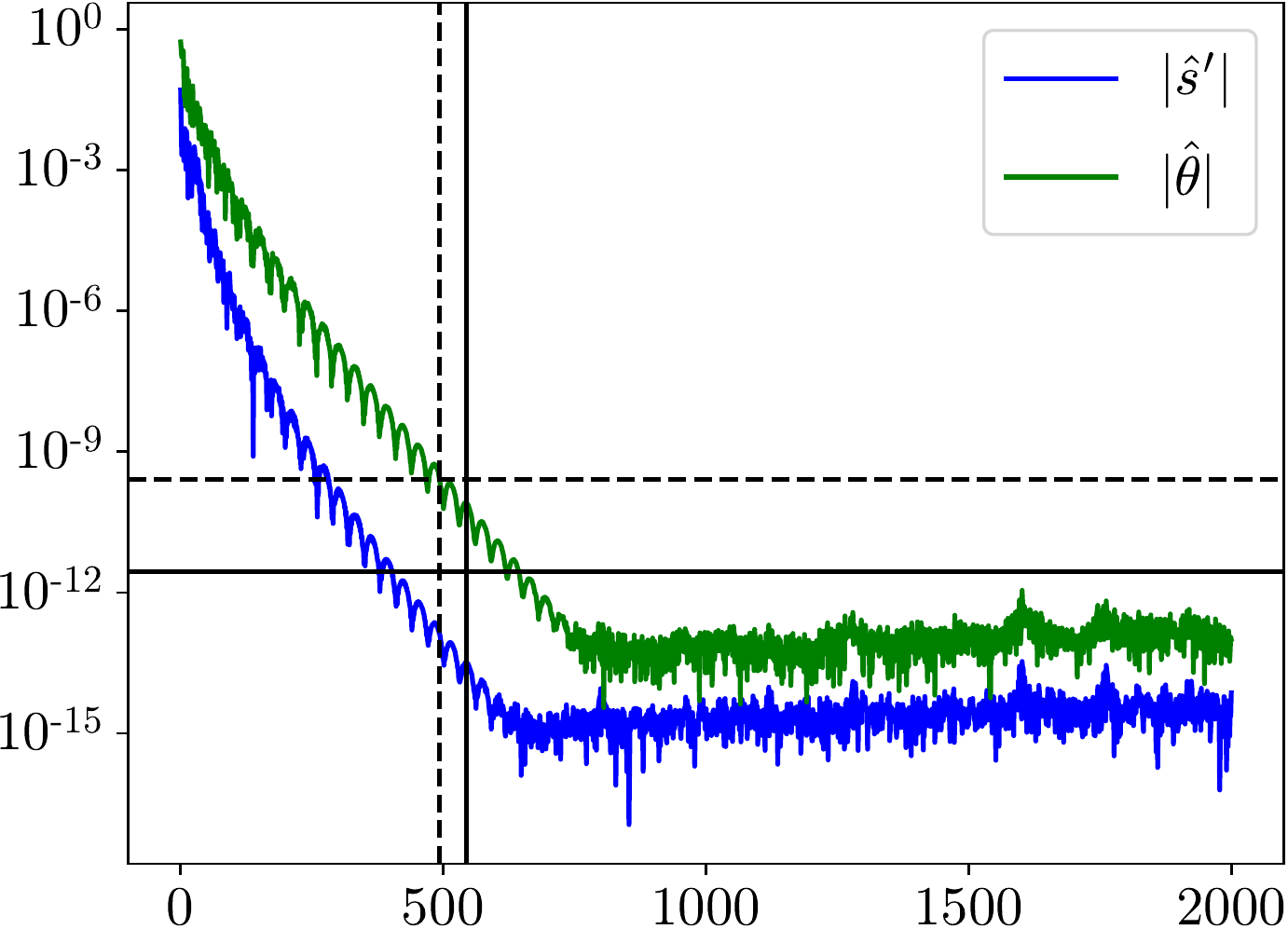}
    \caption{After filtering}
  \end{subfigure}
  \caption{Fourier coefficients of $s'(t)$ and $\theta(t)$ corresponding to
  Figure~\ref{fig:closed6}.
  The value of $\delta_{s'}$ is indicated by a horizontal solid line 
  and the value of $\delta_{\theta}$ is indicated by a horizontal dashed line.
  The $545$th coefficients of $s'(t)$ decays to the $\delta_{s'}$, indicated
  by a vertical solid line. 
  The $494$th coefficients of $\theta(t)$ decays to the $\delta_{\theta}$,
  indicated by a vertical dashed line.}
  \label{fig:closed6cs}
\end{figure}
\begin{figure}[!h]
  \centering
  \begin{subfigure}[b]{0.45\linewidth}
    \includegraphics[width=\linewidth]{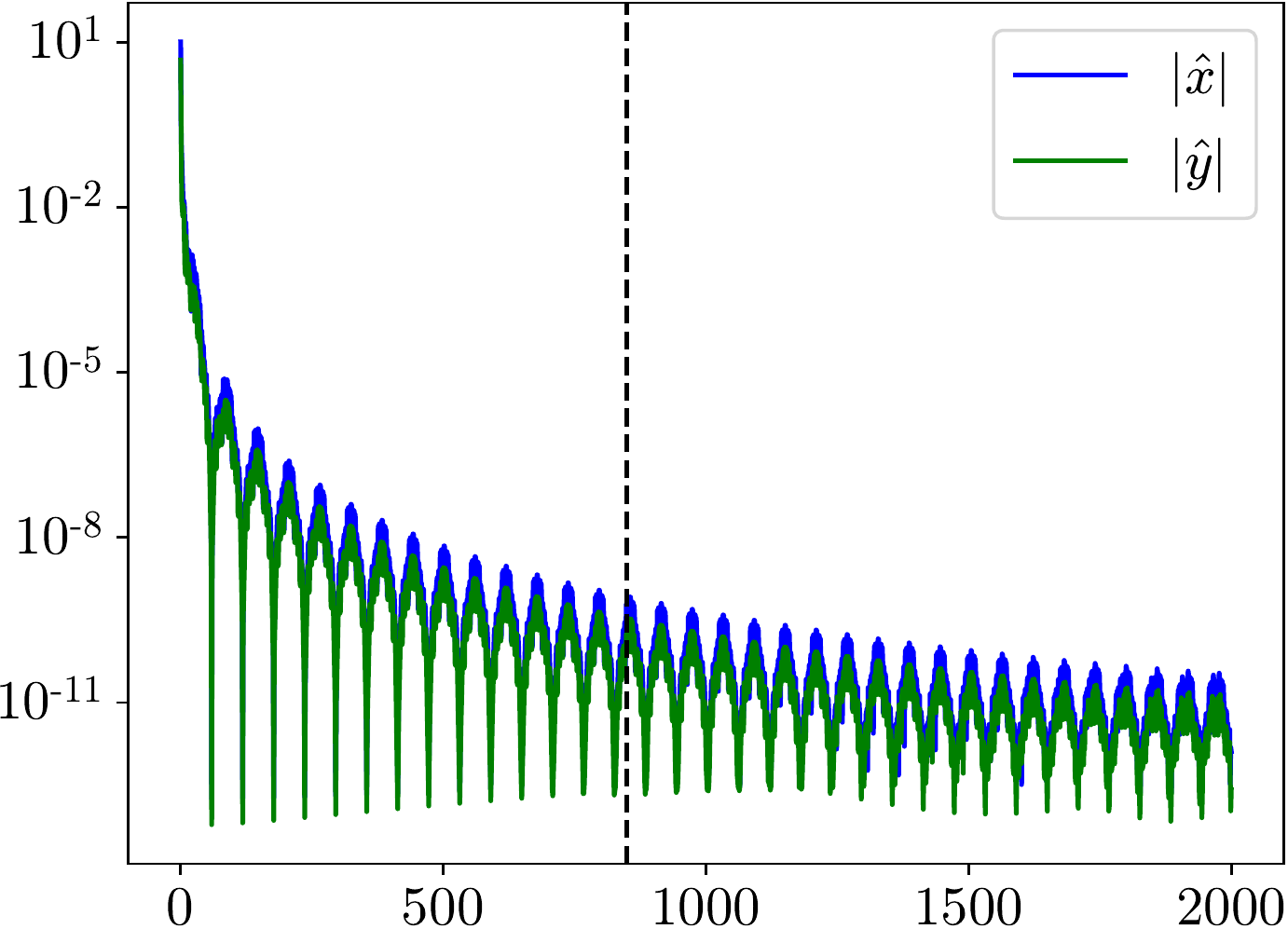}
    \caption{Fourier coefficients of the initial curve}
  \end{subfigure}
  \begin{subfigure}[b]{0.45\linewidth}
    \includegraphics[width=\linewidth]{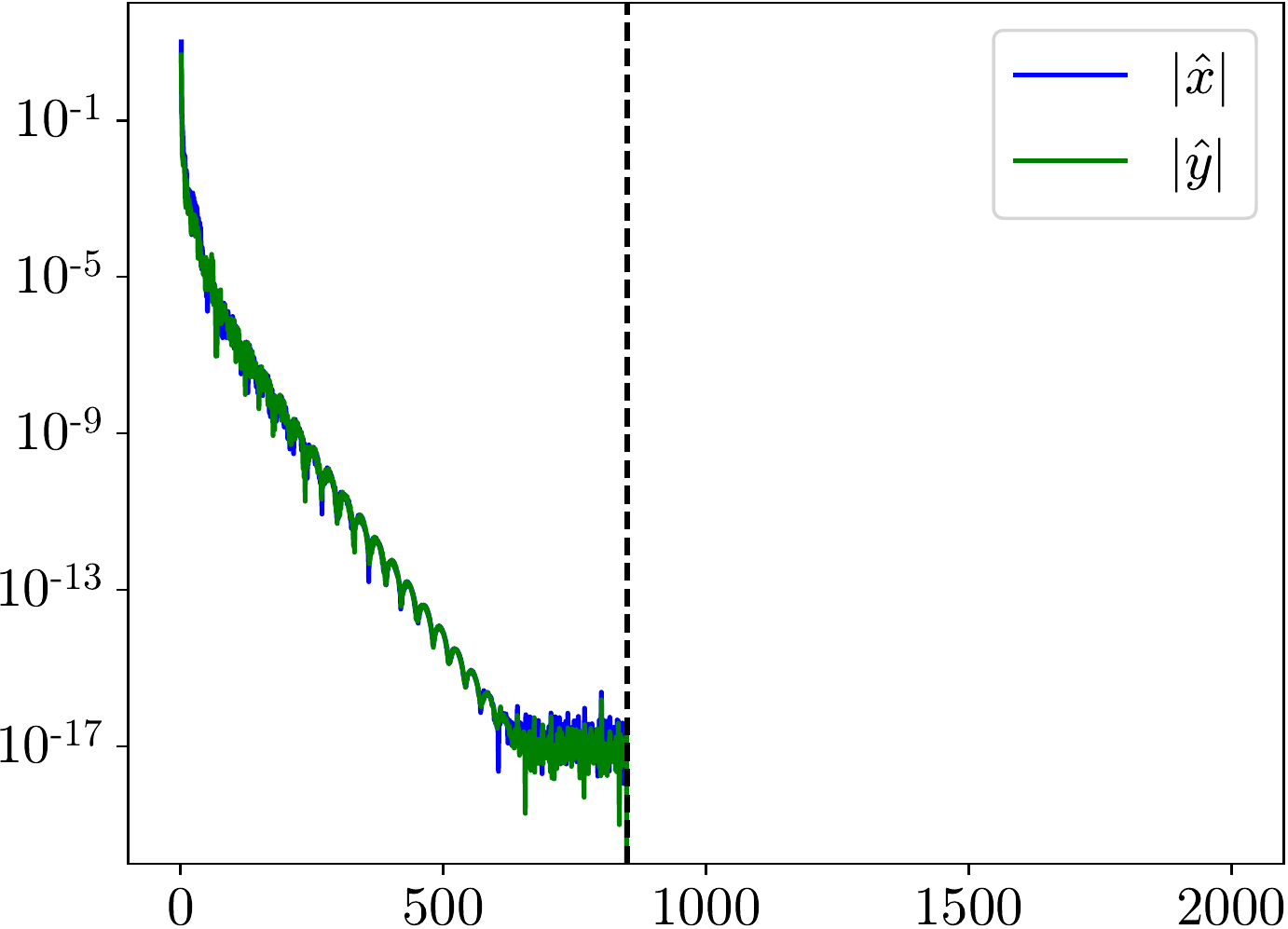}
    \caption{Fourier coefficients of the final curve}
  \end{subfigure}
  \caption{Fourier coefficients of $x(t)$ and $y(t)$ corresponding to
  Figure~\ref{fig:closed6}.
  The value of $n_{\text{coefs}}$ is indicated by a vertical dashed line.} 
  \label{fig:closed6cs2}
\end{figure}

\begin{figure}[!h]
  \centering
  \begin{subfigure}[b]{0.20\linewidth}
    \includegraphics[width=\linewidth]{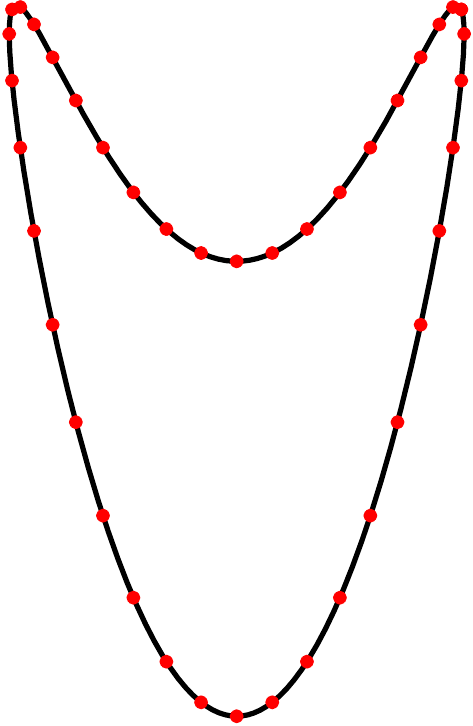}
    \caption{The curve before smoothing}
  \end{subfigure}\hspace{10mm}
  \begin{subfigure}[b]{0.20\linewidth}
    \includegraphics[width=\linewidth]{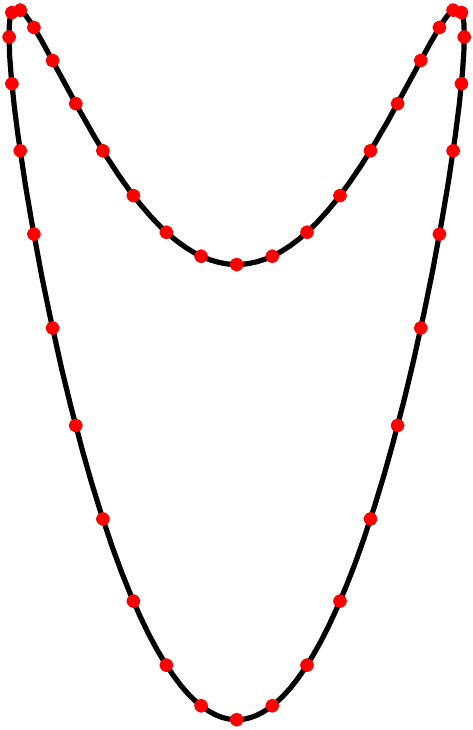}
    \caption{The curve after smoothing}
  \end{subfigure}
  \caption{The result of algorithm applied to Figure $4.3$ in~\cite{bandlimited}. The red dots 
 mark the sample points.}
  \label{fig:closed5}
\end{figure}
\begin{figure}[!h]
  \centering
  \includegraphics[scale=0.5]{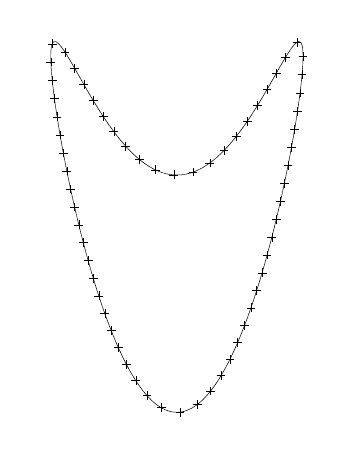}
  \caption{Figure $4.3$ in \cite{bandlimited}}
  \label{fig:closed5ex}
\end{figure}

\begin{figure}[!h]
  \centering
  \begin{subfigure}[b]{0.45\linewidth}
    \includegraphics[width=\linewidth]{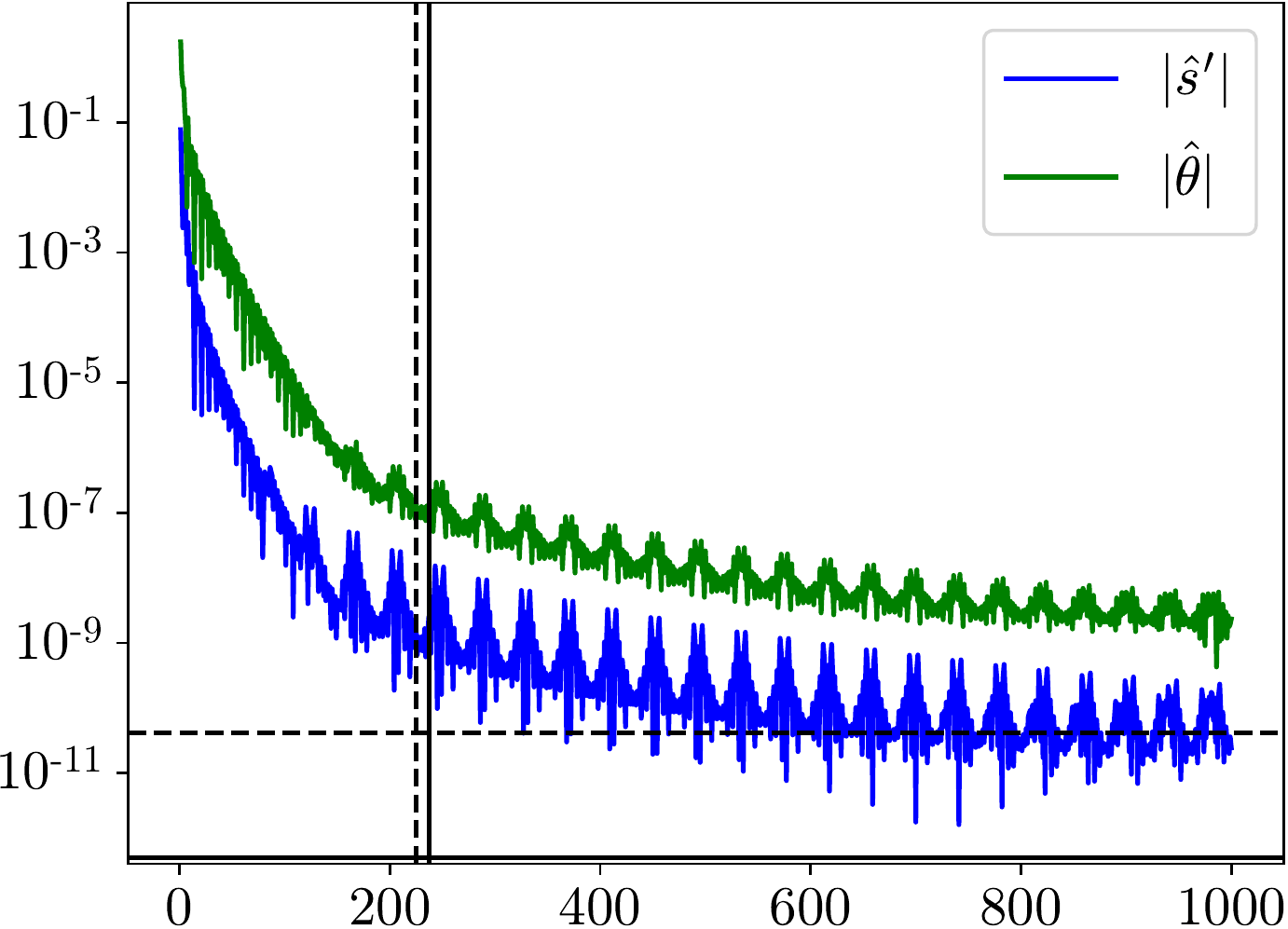}
    \caption{Before filtering}
  \end{subfigure}
  \begin{subfigure}[b]{0.45\linewidth}
    \includegraphics[width=\linewidth]{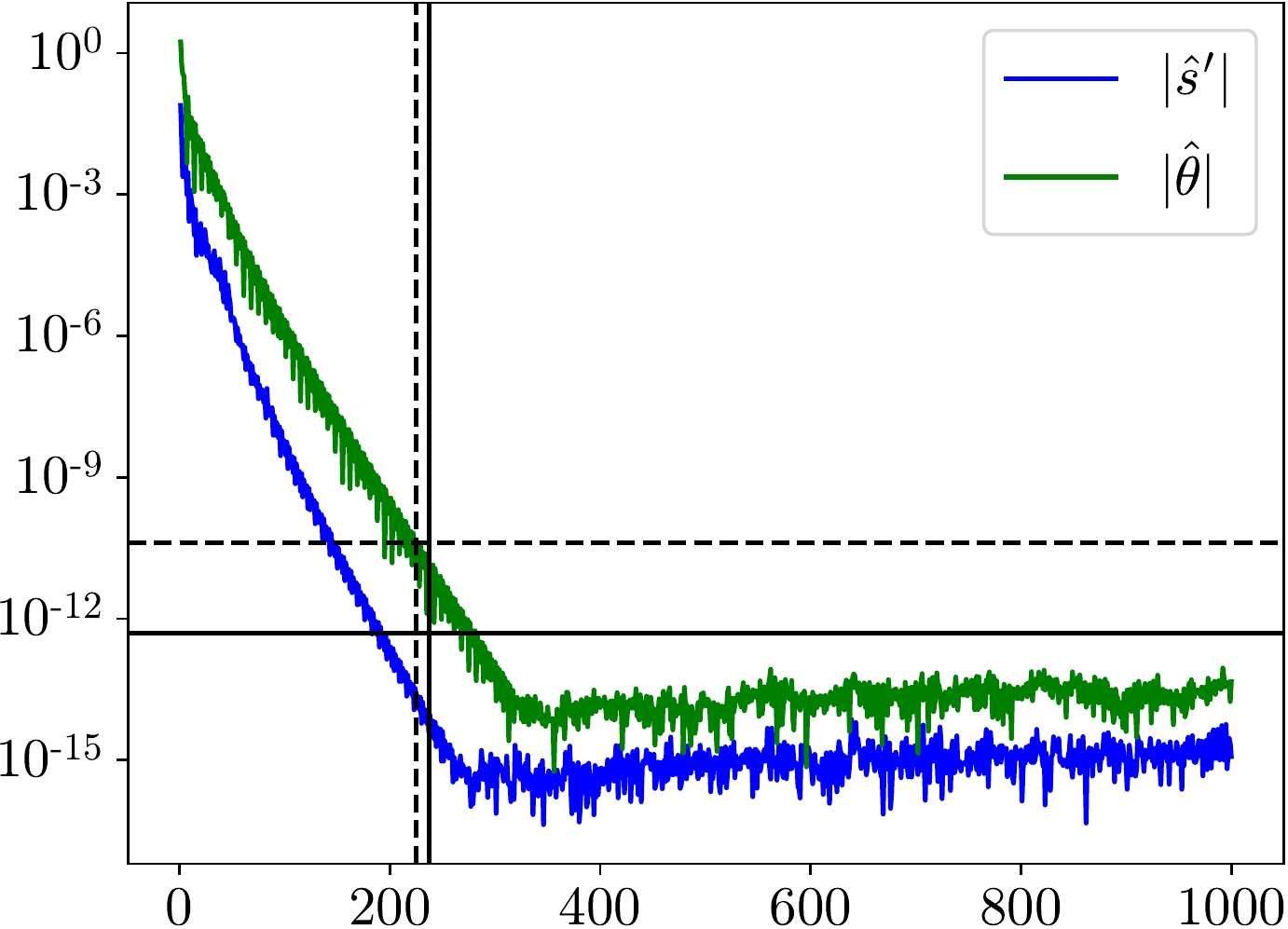}
    \caption{After filtering}
  \end{subfigure}
  \caption{Fourier coefficients of $s'(t)$ and $\theta(t)$ corresponding to
  Figure~\ref{fig:closed5}.
  The value of $\delta_{s'}$ is indicated by a horizontal solid line 
  and the value of $\delta_{\theta}$ is indicated by a horizontal dashed line.
  The $237$th coefficients of $s'(t)$ decays to $\delta_{s'}$, indicated by 
  a vertical solid line.
  The $225$th coefficients of $\theta(t)$ decays to $\delta_{\theta}$,
  indicated by a vertical dashed line.}
  \label{fig:closed5cs}
\end{figure}
\begin{figure}[!h]
  \centering
  \begin{subfigure}[b]{0.45\linewidth}
    \includegraphics[width=\linewidth]{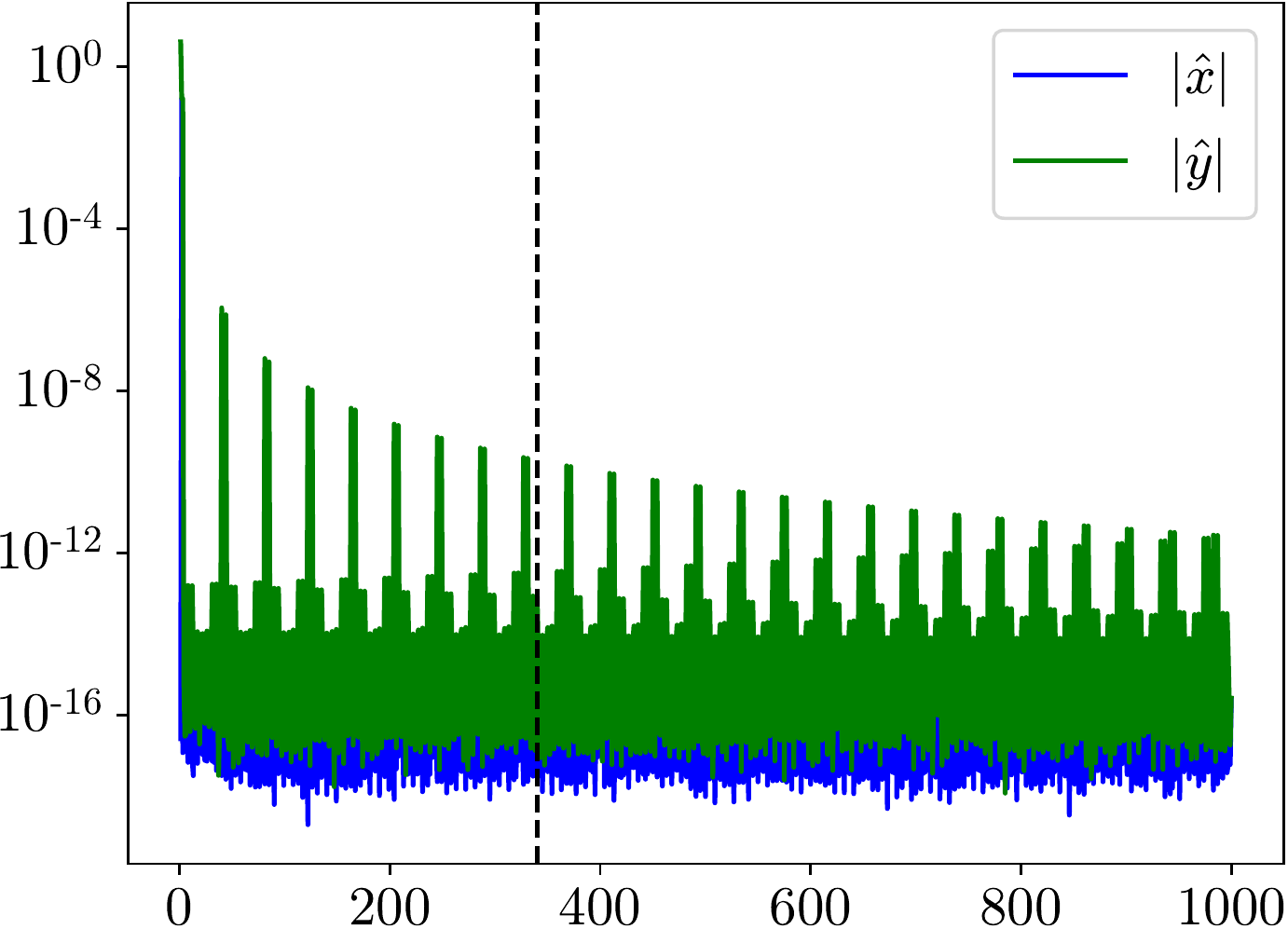}
    \caption{Fourier coefficients of the initial curve}
  \end{subfigure}
  \begin{subfigure}[b]{0.45\linewidth}
    \includegraphics[width=\linewidth]{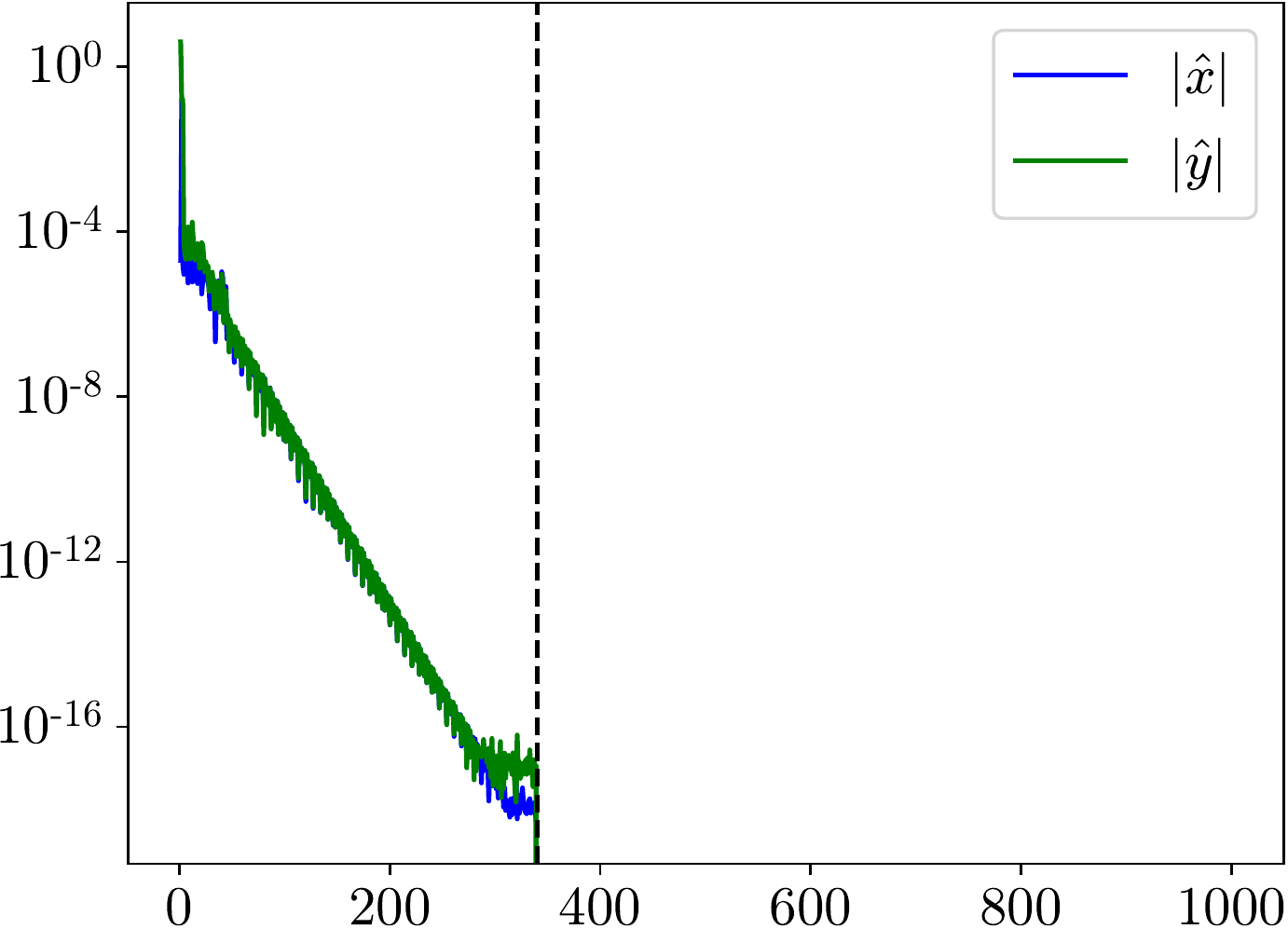}
    \caption{Fourier coefficients of the final curve}
  \end{subfigure}
  \caption{Fourier coefficients of $x(t)$ and $y(t)$ corresponding to
  Figure~\ref{fig:closed5}.
  The value of $n_{\text{coefs}}$ is indicated by a vertical dashed line.} 
  \label{fig:closed5cs2}
\end{figure}

\begin{table}[h]
  \centering
\begin{tabular}{cccccc}
Case & $N=1024$ & $N=2048$ & $N=4096$ & $N=8192$  \\
\hline
\T Figure~\ref{fig:closed1} &$0.44050\e{-03}$& $0.72667\e{-03}$ & $0.13412\e{-02}$& $0.26530\e{-02}$ \\
\T Figure~\ref{fig:closed2} &$0.36162\e{-03}$& $0.61837\e{-03}$ & $0.11924\e{-02}$& $0.24492\e{-02}$ \\
\end{tabular}
\caption{
  \label{tab:table2}
Average runtime per iteration, for the first two closed curves, calculated by determining 
the total runtime for $250$ iterations and  
dividing by the number of iterations.
}
\end{table}
\newpage
\section{Conclusion}
Our algorithm produces a bandlimited curve passing through a set of points, up to machine precision. 
It first constructs a $C^{2}$ B\'ezier spline passing through the points, and then
recursively applies a Gaussian filter to both the derivative of the arc length function and 
the tangential angle of the curve, to control the bandwidth
of the coefficients, followed by smooth corrections. The resulting curve can be represented by a small number of
coefficients, and resembles a smooth curve drawn naturally by hand, free of ringing artifacts.
The algorithm costs O($N\log{N}$) 
operations at each iteration, and the cost can be further reduced by calling
the FFT in the FFTW library~\cite{fftw},  
in which the speed of the FFT routines is optimized for inputs of certain sizes. 

One possible extension of this paper is to design an algorithm for curves and surfaces in $\mathbb{R}^{3}$.   
The main methodology is still applicable, if
we parametrize a curve in $\mathbb{R}^{3}$ by a
function $\gamma(t) \colon I \rightarrow \mathbb{R}^{3}$,
where $I \subset \mathbb{R}$, in terms of 
the same parameter $t$ as in this paper,
and a surface in $\mathbb{R}^{3}$ by a function
$\gamma(s,t) \colon I_1 \times I_2 \rightarrow \mathbb{R}$, 
where $I_1, I_2 \subset \mathbb{R}$,
in terms of both $s$ and $t$. We can apply the Chebyshev or the Fourier approximation 
in each parameter,
depending on whether the curve or surface is periodic in that parameter,
filter the coefficients and 
add smooth perturbations in a similar way.  
Another application is to implement the algorithm of this paper as a geometric primitive 
in CAD/CAM systems.
Since primitives are generally defined as level sets of polynomials (see Chapter $2$ of~\cite{cad}),
the techniques in this paper could be used for the 
constructions of more general $C^{\infty}$ shapes in CAD/CAM systems.

\newcommand{\sectionbreak}{\clearpage}


\end{document}